\documentclass[12pt]{article}
\usepackage{floatrow}
\usepackage{float}
\newfloatcommand{capbtabbox}{table}[][\FBwidth]

\usepackage{algorithm}
\usepackage{algpseudocode}

\usepackage{amsmath}
\usepackage{mathrsfs}
\usepackage{keyval}
\usepackage{multirow}
\usepackage{xcolor}
\usepackage{latexsym,amsmath,amsfonts,amscd,amsthm}
\usepackage{graphicx}
\usepackage{changebar}
\numberwithin{equation}{section}
\theoremstyle{plain}
\newtheorem{thm}{Theorem}[section]

\theoremstyle{definition}
\newtheorem{exam}[thm]{Example}

\begin{document}
\baselineskip=1.5pc

\vspace{.5in}

\begin{center}

{\large\bf High-order WENO-based semi-implicit Newton-type fast sweeping methods for static Hamilton-Jacobi equations}

\end{center}

\vspace{.1in}

\centerline{
Yuan Liu  
\footnote{Department of Mathematics, Statistics and Physics,
Wichita State University, Wichita, KS 67260, USA.
E-mail: yuan.liu@wichita.edu. Research is supported in part by NSF grant DMS-2213436.}, 
\qquad
Jianliang Qian\footnote{Department of Mathematics and Department of Computational Mathematics, Science and Engineering, Michigan State University, East Lansing, MI 48824, USA.
E-mail: jqian@msu.edu. Research is supported in part by NSF grants 2152011 and 2309534, and an MSU SPG grant.} 
}

\vspace{.2in}

\centerline{\bf Abstract}

%\vspace{.1in}
In this paper, we propose high-order weighted essentially non-oscillatory (WENO)-based semi-implicit Newton-type Gauss–Seidel Lax–Friedrichs fast sweeping methods for solving the generalized Eikonal equation arising in wave propagation through a moving fluid. Building upon the Newton-type framework of Li and Qian \cite{li2020newton}, which updates the solution line-wise using Newton's method with a tridiagonal and strictly diagonally dominant Jacobian, we extend the local solver to fifth-, seventh-, and ninth-order accuracy by incorporating high-order WENO approximations of the spatial derivatives into the numerical Hamiltonian. Three alternating sweeping strategies, column-wise, row-wise, and column-row-wise are considered in the Gauss–Seidel iteration framework. Numerical examples in both two and three spatial dimensions demonstrate the efficiency and accuracy of the proposed schemes.
%\bigskip

{\bf Key Words:}
Hamilton Jacobi, Numerical fluxes, Fast sweeping, line-wise sweeping, WENO reconstruction, Newton iteration, Gauss-Seidel iteration.
%Sparse grid; Discontinuous Galerkin; Reaction-diffusion.
%Hyperbolic conservation laws, maximum principle preserving, positivity preserving, unstructured meshes, finite volume schemes, WENO schemes, compressible Euler system.

%{\bf AMS(MOS) subject classification:} 65M99

\pagenumbering{arabic}

%section 1

% ===== BEGIN introduction.tex =====
\section{Introduction}
\label{sec1}
\setcounter{equation}{0}
\setcounter{figure}{0}
\setcounter{table}{0}

In this paper, we are interested in solving the following boundary value problem for static HJ equations  
\begin{align}\label{eq1}
\begin{cases}
H(\mathbf{x}, \nabla T(\mathbf{x})) = f(\mathbf{x}), \qquad \mathbf{x} \in \Omega,
\\
T(\mathbf{x}) = g(\mathbf{x}), \qquad \qquad \qquad \mathbf{x} \in \Gamma,
\end{cases}
\end{align}
where $\Omega \subset \mathcal{R}^n$ is a computational domain with boundary conditions prescribed on the subset $\Gamma$, $T(\mathbf{x})$ is the travel time of wave propagation, and $H$ is the Hamiltonian which is usually a nonlinear Lipschitz continuous function. Particularly, if the Hamiltonian 
\begin{align}\label{geikonal}
H(\mathbf{x}, \nabla T(\mathbf{x})) = || \nabla \text{T}(\mathbf{x}) || -\frac{|1-\mathbf{v(x)} \cdot \nabla \text{T}(\mathbf x) |}{F(\mathbf x)},
\end{align}
then (\ref{eq1}) is called the generalized Eikonal equation, which is frequently used to model the travel time $T(\mathbf{x})$ of wave propagation in an isotropic acoustic medium occupied by a moving fluid. Here, $F(\mathbf{x})$ denotes the speed of wave propagation and $\mathbf{v(x)}$ is the velocity of the background moving fluid. When $\mathbf{v(x)} = 0$, we obtain the further reduced case as the standard isotropic Eikonal equation with Hamiltonian  
\begin{align}\label{eikonal}
H(\mathbf{x}, \nabla T(\mathbf{x})) = || \nabla \text{T}(\mathbf{x}) ||.
\end{align}

The generalized Eikonal equation is a prototypical static Hamilton–Jacobi (HJ) equation of broad practical significance and the background velocity field $\mathbf{v}(\mathbf{x})$ plays a central role in determining the traveltime solution across a wide range of applications. In seismic wave propagation, the first-arrival traveltime $T(\mathbf{x})$ from a source to receivers is influenced by crustal deformation and mantle convection, whose effects are encoded in $\mathbf{v}(\mathbf{x})$; accurate recovery of this field is therefore essential for seismic imaging and traveltime tomography \cite{leuqia06,lileuqia14,leuqiahu21,siphiweiqia25}. In underwater acoustics, acoustic wavefronts propagating through the ocean must account for spatially varying current velocities, which can substantially distort the traveltime field \cite{brekhovskikh2003fundamentals, kornhauser1953ray}. Beyond these geophysical settings, the equation also arises in geometric optics modeling of steady wind flow over mountains \cite{tanushev2007mountain}, structural geology fold modeling \cite{hjelle2011hamilton}, and sonic boom propagation from cruising and supersonic aircraft \cite{blumrich2005meteorologically, loubeau2009effects}.

Numerically simulating the static HJ equations presents significant challenges due to the nonlinear nature of the Hamiltonian and the possible formation of singularities in the viscosity solution, and the numerical solution of the Eikonal and general static HJ equations has been pursued along two broad lines: Lagrangian and Eulerian approaches. Classical Lagrangian ray-tracing methods \cite{cerveny2001seismic, bois1972well} solve a system of ordinary differential equations along the characteristics of the HJ equation, tracing a bundle of rays along which the traveltime is evaluated. These methods are efficient in smooth media, but yield solutions on non-uniformly distributed point clouds that can leave shadow zones with few or no rays in complex heterogeneous environments, which would become problematic for applications that require a globally resolved, uniformly accurate traveltime field on a regular computational grid. Eulerian grid-based methods overcome this by solving the HJ equation directly on a mesh, and two efficient families have emerged: fast marching methods \cite{tsitsiklis1995efficient, sethian1996fast} and fast sweeping methods \cite{zhao2005fast, tsai2003fast}.

Fast marching methods propagate the solution causally in $O(N \log N)$ operations, but enforcing causality for the generalized eikonal equation, where the wavefront advances in both normal and tangential directions, requires careful algorithmic design \cite{dahiya2013characteristic, ho2019improved}. Fast sweeping methods instead achieve convergence through alternating-direction Gauss–Seidel iterations that simultaneously cover multiple characteristic directions, naturally handling non-convex Hamiltonians without strict causality when a Lax-Friedrichs type numerical flux is used. For the generalized Eikonal equation, \cite{kao2004lax} proposed the Lax–Friedrichs sweeping algorithm, which discretizes the Hamiltonian with a Lax–Friedrichs numerical flux, yielding a simple point-wise updating formula applicable to general static HJ equations. \cite{zhang2006high} extended this to high-order accuracy using weighted essentially non-oscillatory (WENO) approximations \cite{jiang2000weighted}, and \cite{glowinski2016operator} developed operator-splitting based fast sweeping methods for the generalized Eikonal equation by exploiting the additive structure of its Hamiltonian.

Several more high-order extensions have since been proposed. \cite{xiong2010fast} designed a fifth-order WENO fast sweeping method with accurate boundary treatment. \cite{ren2020high} proposed Hermite WENO fast sweeping methods achieving fifth-order accuracy with compact stencils. A significant convergence issue with classical high-order (fifth-order and above) WENO fast sweeping is that iteration residues may stagnate at truncation-error level rather than converging to machine zero, complicating the choice of stopping criteria. \cite{li2021absolutely} resolved this by introducing absolutely convergent fixed-point fast sweeping methods based on multi-resolution WENO (MR-WENO) local solvers, whose residues settle to round-off errors. \cite{hu2024high} subsequently extended this framework to eikonal and factored eikonal equations. See also a recent work which uses sparse grid implementations to reduce computational cost in higher dimensions \cite{miksis2024sparse}.

%, and physics-informed neural network approaches for eikonal-type equations %\cite{waheed2021pinneik, grubas2023neural}.

Despite such progress, one major limitation shared by all the above Lax–Friedrichs sweeping algorithms is that the Gauss–Seidel updates are carried out point by point: $T_{i,j}$ is updated one grid point at a time in scalar fashion. This leads to: (1) relatively slow convergence requiring many sweepings, especially on fine meshes and in 3D; and (2) poor utilization of modern hardware optimized for vectorized array operations.
To overcome these limitations, Li and Qian \cite{li2020newton} proposed the Newton-type Gauss–Seidel Lax–Friedrichs sweeping algorithm, which replaces point-wise updates with a line-wise sweeping strategy so that the entire solution column $\{T_{i,j},\, i=1,\ldots,N\}$ is updated simultaneously via Newton's method. The Jacobian of the resulting Newton system is tridiagonal and strictly diagonally dominant, so each Newton linear system is solved in $O(N)$ operations by the Thomas algorithm, preserving $O(N^2)$ per-sweep complexity in 2D while executing entirely through vector operations. This semi-implicit character accelerates convergence, and numerical experiments in \cite{li2020newton} show that the Newton-type algorithm requires roughly half as many sweepings as the point-wise scheme across all tested configurations, and achieves decisive wall-clock speedups in 3D and for high-order schemes. Li and Qian \cite{li2020newton} developed first-order and third-order (WENO3) versions of this framework. The present paper extends the Newton-type Gauss–Seidel Lax–Friedrichs sweeping framework to fifth- (WENO5), seventh- (WENO7), and ninth-order (WENO9) accuracy. The extension of the Newton-type framework requires non-trivial considerations, as wider WENO stencils necessitate appropriate initialization near source locations and consistent boundary treatment compared with WENO3. Moreover, we show that incorporating higher-order WENO approximations into the Jacobian $\nabla \mathbf{F}_j$ preserves its tridiagonal and strictly diagonally dominant structure at all orders considered, so the Thomas algorithm still applies and the $O(N)$ complexity per Newton step is maintained. 

The remainder of this paper is organized as follows. In Section 2, we will give a brief review of the explicit point-wise Gauss-Seidel fast sweeping methods. The proposed high-order Newton-type fast sweeping schemes are presented in Section 3. Numerical examples demonstrating the performance of the proposed methods are provided in Section 4. Finally, concluding remarks are given in Section 5.

%{\bf Remark by Qian: fill in the details here}

% ===== END introduction.tex =====

% section 2
% ===== BEGIN algorithm.tex =====
%\section{Numerical methods}

\section{Point-wise Gauss-Seidel fast sweeping methods}
\setcounter{equation}{0}
\setcounter{figure}{0}
\setcounter{table}{0}

%{\bf Q. You should input those detailed numerical schemes in the context.}
For notational simplicity, we present our numerical schemes for the 2D case when the Hamiltonian takes the following form,
\begin{equation}\label{h1}
H(u, v) = \sqrt{u^2+v^2} - \frac{|1-v_1 u -v_2 v|}{F(\mathbf{x})}
\end{equation}
with $\nabla T = (u, v)$ and $\mathbf{v}(\mathbf{x}) = (v_1, v_2)$. Following \cite{li2020newton}, the Lax-Friedrichs numerical Hamiltonian is applied to discretize the Hamiltonian (\ref{eq1}), which takes the form
\begin{equation}\label{h4}
\hat{H}^{\text{LF}}(u^+, u^-; v^+, v^-) = H(\frac{u^++u^-}{2}, \frac{v^++v^-}{2})-\alpha^x \cdot \frac{u^+-u^-}{2} -\alpha^y \cdot \frac{v^+-v^-}{2}, 
\end{equation}
where $u^{\pm}$ and $v^{\pm}$ are reconstructed approximations for the derivatives $T_x$ and $T_y$ from the right and left sides in the $x$- and $y$-direction, respectively. $\alpha^x$ and $\alpha^y$ are viscosity 
parameters, respectively, in the $x$- and $y$-direction, satisfying 
\begin{equation}
\alpha^x > \max \left |\frac{\partial H}{\partial u}\right |, \qquad \alpha^y > \max \left |\frac{\partial H}{\partial v}\right |.
\end{equation}
For the generalized Eikonal equation with Hamiltonian (\ref{h1}), we choose 
\begin{equation}
\alpha^x = \alpha \left (1+\max \left|\frac{v_1}{F}\right |\right ), \qquad \alpha^y = \alpha \left (1+\max \left |\frac{v_2}{F}\right |\right ), 
\end{equation}
with
\begin{equation}\label{h2}
\begin{aligned}
 \frac{\partial H}{\partial u} = \frac{u}{\sqrt{u^2+v^2}} + \frac{v_1}{F(\mathbf{x})} \cdot \text{sgn}(1-v_1 u-v_2 v), \\
 \frac{\partial H}{\partial v} = \frac{v}{\sqrt{u^2+v^2}} + \frac{v_2}{F(\mathbf{x})} \cdot \text{sgn}(1-v_1 u-v_2 v), 
 \end{aligned}
\end{equation}
where sgn($\cdot$) represents the signum function. Moreover, in the numerical simulation, we will add a small positive number $10^{-12}$ to $\sqrt{u^2+v^2}$ to avoid zero denominator in (\ref{h2}). 

We consider solving (\ref{eq1}) on the rectangular domain $\Omega$, and perform a shape regular partition
$\Omega = \cup_{i, j} R_{i,j}$, where $R_{i, j} = I_i \times J_j = [x_i, x_{i+1}] \times [y_j, y_{j+1}]$ for $i = 1, ..., M$ and $j = 1, ..., N$. The spatial mesh sizes are uniform in each direction, denoting as $\Delta x$ and $\Delta y$.
Then, a Lax-Friedrichs scheme reads as 
\begin{equation}\label{h5}
H \left (\frac{u^+_{ij}+u^-_{ij}}{2}, \frac{v^+_{ij}+v^-_{ij}}{2}\right )-\alpha^x \frac{u^+_{ij}-u^-_{ij}}{2} -\alpha^y \frac{v^+_{ij}-v^-_{ij}}{2} =f_{ij}
\end{equation}
with $f_{ij} = f(x_i, y_j)$. Further derivation will yield 
\begin{equation}\label{5.5}
F_{ij} (T) = 0
\end{equation}
where
\begin{align}\label{h6}
F_{ij}(T) ={}& H\left(\frac{T_{i+1,j}-T_{i-1,j}}{2\Delta x},
\frac{T_{i,j+1}-T_{i,j-1}}{2\Delta y}\right) \nonumber\\
&-\alpha^x \frac{T_{i+1,j}+T_{i-1,j}-2T_{i,j}}{2\Delta x}
-\alpha^y \frac{T_{i,j+1}+T_{i,j-1}-2T_{i,j}}{2\Delta y} - f_{ij}, 
\end{align}
and
\begin{align}\label{h7}
& T_{i-1,j} = T_{ij} - \Delta x (T_x)^-_{ij}\,;  \qquad T_{i+1,j} = T_{ij} + \Delta x (T_x)_{ij}^+\,; \\ \nonumber
& T_{i, j-1} = T_{ij} - \Delta y (T_y)^-_{ij}\,; \qquad T_{i, j+1} = T_{ij} +\Delta y (T_y)_{ij}^+\,.
\end{align}
In \cite{kao2004lax, zhang2006high}, point-wise implicit Lax-Friedrichs sweeping methods based on Gauss-Seidel iteration with alternating-direction sweepings are proposed to simulate the static HJ equation,
\begin{align}\label{h8}
T^{\text{new}}_{i,j} ={}& T^{\text{old}}_{i,j}
+ \left(\frac{1}{\frac{\alpha^x}{\Delta x}+\frac{\alpha^y}{\Delta y}}\right)
\Biggl[ f_{ij}
- H\left( \frac{T_{i+1,j}-T_{i-1,j}}{2\Delta x},
\frac{T_{i,j+1}-T_{i,j-1}}{2\Delta y} \right) \nonumber\\
&+\alpha^x \frac{T_{i+1,j}+T_{i-1,j}}{2\Delta x}
+\alpha^y \frac{T_{i,j+1}+T_{i,j-1}}{2\Delta y} \Biggr].
\end{align}
Particularly, first order approximation to $(T_x)^{\pm}_{ij}$ and $(T_y)^{\pm}_{ij}$ in \cite{kao2004lax} are employed while third order WENO reconstructed approximations are used in \cite{zhang2006high}. 

\section{Newton-type fast sweeping method}
Instead of explicit point-wise iteration, following the ideas initiated in \cite{zhang2006high,li2020newton}, we are interested in systematically developing higher-order WENO-based Newton-type Gauss-Seidel sweeping methods to solve (\ref{h5}) in the following semi-implicit manner, i.e., explicitly using higher-order WENO reconstructions to approximate the derivatives in numerical Hamiltonians \cite{zhang2006high} and implicitly solving the Newton linear system \cite{li2020newton}. One significant advantage of our semi-implicit treatment is that, we are able to bring high-order derivative information into the numerical Hamiltonian at the cost of only solving a tridiagonal linear system at each Newton step. 

\subsection{Line-wise Newton iteration}
We propose to solve the nonlinear equations (\ref{h6}) by Newton iteration in line-wise fashion. The general architecture of our Newton-type solver is based on the first order approximation to the derivative terms, i.e.,
\begin{align}\label{h9}
& (T_x)^-_{i,j} = \frac{T_{i,j}-T_{i-1,j}}{\Delta x}, \qquad (T_x)^+_{i,j} = \frac{T_{i+1,j}-T_{i,j}}{\Delta x}, \\ &
(T_y)^-_{i,j} =\frac{T_{i,j}-T_{i,j-1}}{\Delta y}, \qquad (T_y)^+_{i,j} = \frac{T_{i,j+1}-T_{i,j}}{\Delta y}.
\end{align}
First of all, we present our solver for (\ref{h6}) in column-wise fashion, i.e., we propose to solve the $j$-th column of $\{T_{ij}, i=1,2,..., M \}$ simultaneously while treating all the other columns of $T_{i,j}$ explicitly, i.e.,  (\ref{h6}) can be reformulated as
\begin{align}\label{h10_0}
F_{ij}(T^*_{i-1,j}, T^*_{ij}, T^*_{i+1,j}) = 0
\end{align}
and
\begin{align}\label{h10}
F_{ij}(T^*_{i-1,j}, T^*_{i,j}, T^*_{i+1,j}) ={}&
H\left(\frac{T^*_{i+1,j}-T^*_{i-1,j}}{2\Delta x},
\frac{T_{i,j+1}-T_{i,j-1}}{2\Delta y}\right) \nonumber\\
&\hspace{-8.5em}{}-
\alpha^x \frac{T^*_{i+1,j}+T^*_{i-1,j}}{2\Delta x}
-\alpha^y \frac{T_{i,j+1}+T_{i,j-1}}{2\Delta y}
-\left(\frac{\alpha^x}{\Delta x}+\frac{\alpha^y}{\Delta y}\right)T^*_{i,j}-f_{i,j}, 
\end{align}
where the superscript $*$ indicates the variables to be solved at the current iteration. For the interior points of the $j$-th column, $\mathbf{T_j} = (T_{2,j}, \cdots , T_{M-1, j})^T$, we solve the system of $(M-2)$ equations $\mathbf{F_j (T_j)} = 0$ with the following component-wise form, 
\begin{equation}\label{h11}
\begin{cases}
F_{2, j} (T_{1, j}, T^*_{2,j}, T^*_{3, j})  = 0;\\
F_{3, j}(T^*_{2,j}, T^*_{3,j}, T^*_{4,j})  = 0; \\
\qquad \vdots \\
F_{M-2, j} (T^*_{M-3, j}, T^*_{M-2, j}, T^*_{M-1,j}) = 0; \\
F_{M-1, j} (T^*_{M-2, j}, T^*_{M-1,j}, T_{M,j}) = 0.
\end{cases}
\end{equation}
Therefore, Newton's method yields 
\begin{equation}\label{h12}
\mathbf{T}^{\text{new}}_{\mathbf{j}} = \mathbf{T}^{\text{old}}_{\mathbf{j}} - (\mathbf{\nabla \mathbf{F_j})^{-1} \cdot \mathbf{F_j}}, 
\end{equation}
where $\nabla \mathbf{F_j}$ is the Jacobian matrix of the nonlinear function $\mathbf{F_j}$ with the following form 

\[
 \nabla \mathbf{F_j} = \begin{bmatrix}
  \frac{\partial F_{2, j}}{\partial T_{2, j}} &  \frac{\partial F_{2, j}}{\partial T_{3, j}}      \\ \\
  \frac{\partial F_{3, j}}{\partial T_{2,j}} & \frac{\partial F_{3, j}}{\partial T_{3,j}} & \frac{\partial F_{3, j}}{\partial T_{4,j}}     \\  
      \\
     & \ddots & \ddots  &  \ddots  \\
           \\
         & &\frac{\partial F_{M-2, j}}{\partial T_{M-3,j}} &\frac{\partial F_{M-2, j}}{\partial T_{M-2,j}} &  \frac{\partial F_{M-2, j}}{\partial T_{M-1,j}}  \\ \\ 
& & &\frac{\partial F_{M-1, j}}{\partial T_{M-2,j}} &  \frac{\partial F_{M-1, j}}{\partial T_{M-1,j}}
 \end{bmatrix}_{(M-2) \times (M-2)}   
\] 
and can be evaluated explicitly 
\[
 \nabla \mathbf{F_j} = \begin{bmatrix}
  a_2 &  c_2     \\ \\
  b_3 & a_3 & c_3    \\  
      \\
     & \ddots & \ddots  &  \ddots  \\
           \\
         & &b_{M-2} &a_{M-2} &  c_{M-2} \\ \\ 
& & &b_{M-1} & a_{M-1}
 \end{bmatrix}_{(M-2) \times (M-2)}   
\] 
with
\begin{align}\label{h13}
& a_i = \frac{\partial F_{i, j}}{\partial T_{i, j}} =\frac{\alpha^x}{\Delta x} +\frac{\alpha^y}{\Delta y},   && i=2, 3, ..., M-1;\\
& b_i = \frac{\partial F_{i,j}}{\partial T_{i-1,j}} = - \left (\frac{\partial H}{\partial u} \right)_{i,j} \cdot \frac{1}{2\Delta x} - \frac{\alpha^x}{2\Delta x};  && i = 3, 4, ..., M-1; \\
& c_i = \frac{\partial F_{i,j}}{\partial T_{i+1, j}} = \left( \frac{\partial H}{\partial u} \right)_{i, j} \cdot \frac{1}{2\Delta x} - \frac{\alpha^x}{2\Delta x}, && i=2, 3, ..., M-2.
\end{align}
As proved in \cite{li2020newton}, the Jacobian matrix $\nabla F_j$ is strictly diagonally dominant and thus non-singular. 

Similarly, we can also solve (\ref{h5}) in row-wise fashion by the implicit Newton iteration, 
\begin{equation}\label{h14}
\mathbf{T}^{\text{new}}_{\mathbf{i}} = \mathbf{T}^{\text{old}}_{\mathbf{i}} - (\mathbf{\nabla \mathbf{F_i})^{-1} \cdot \mathbf{F_i}},
\end{equation}
where $\mathbf{T_i} = (T_{i,2}, \cdots , T_{i,N-1})^T $ are the interior numerical unknowns of the $i$-th row, and 
\[
 \nabla \mathbf{F_i} = \begin{bmatrix}
  a_2 &  c_2     \\ \\
  b_3 & a_3 & c_3    \\  
      \\
     & \ddots & \ddots  &  \ddots  \\
           \\
         & &b_{N-2} &a_{N-2} &  c_{N-2} \\ \\ 
& & &b_{N-1} & a_{N-1}
 \end{bmatrix}_{(N-2) \times (N-2)}   
\] 
with
\begin{align}\label{h15}
& a_j = \frac{\partial F_{i, j}}{\partial T_{i, j}} =\frac{\alpha^x}{\Delta x} +\frac{\alpha^y}{\Delta y},   && j=2, 3, ..., N-1;\\
& b_j = \frac{\partial F_{i,j}}{\partial T_{i,j-1}} = - \left (\frac{\partial H}{\partial v} \right)_{i,j} \cdot \frac{1}{2\Delta y} - \frac{\alpha^x}{2\Delta y};  && j = 3, 4, ..., N-1; \\
& c_j = \frac{\partial F_{i,j}}{\partial T_{i, j+1}} = \left( \frac{\partial H}{\partial v} \right)_{i, j} \cdot \frac{1}{2\Delta y} - \frac{\alpha^x}{2\Delta y}, && j=2, 3, ..., N-2.
\end{align}

On the computational boundary, linear extrapolation is used, 
\begin{align}\label{h16}
& T_{1, j} = \max\{2T_{2,j}-T_{3,j}, T_{3,j}\}; \qquad  T_{M, j} = \max\{2T_{M-1,j}-T_{M-2,j}, T_{M-2,j}\}; \\
&  T_{i, 1} = \max\{2T_{i,2}-T_{i,3}, T_{i,3}\}; \qquad  T_{i, N} = \max\{2T_{i,N-1}-T_{i,N-2}, T_{i,N-2}\}.
\end{align}

\subsection{Gauss–Seidel iterations with line-wise alternating sweepings}
Moreover, we incorporate Gauss-Seidel iteration with alternating directional sweep into the Newton iteration (\ref{h12}) or (\ref{h14}). Specifically, we propose to solve (\ref{h5}) line-wisely with alternating directional sweepings. During each iteration, the whole domain will be swept twice for the one sided column-wise and row-wise sweepings, which is the same as the explicit point-wise Gauss-Seidel iteration. For the two-sided column-row-wise sweeping, the whole domain will be swept four times om pme oteration. Below we illustrate the pseudo codes based on one iteration from $\mathbf{T}^n$ to $\mathbf{T}^{n+1}$.

 \begin{algorithm}
 \caption{Column-wise sweeping Newton-type method}\label{alg:column}
 \begin{algorithmic}
 \State $\mathbf{T}^{\text{old}} = \mathbf{T}^n;$ 
 \For {j=2 to N-1} \\
 \quad Compute $\mathbf{T}^{\text{new}}_{\mathbf{j}} = (T_{2, j}, T_{3, j}, \ldots, T_{M-1, j})^T$ by Newton iteration (\ref{h12}); \\
 \quad Calculate $T^{\text{new}}_{1, j}$ and $T^{\text{new}}_{M, j}$ by boundary treatment (\ref{h16});
 \EndFor \\
 \quad Compute $(T^{\text{new}}_{1, j}, \ldots, T^{\text{new}}_{M,j})^T$ by boundary treatment (\ref{h16}) for $j=1$ and $j=N$;\\
 \quad $\mathbf{T}^{\text{old}} = \mathbf{T}^{\text{new}}$;
 \For {j=N-1 to 2} \\
 \quad Compute $\mathbf{T}^{\text{new}}_{\mathbf{j}} = (T_{2, j}, T_{3, j}, ... T_{M-1, j})^T$ by Newton iteration (\ref{h12}); \\
 \quad Calculate $T^{\text{new}}_{1, j}$ and $T^{\text{new}}_{M, j}$ by boundary treatment (\ref{h16});
 \EndFor \\
 \quad Compute $(T^{\text{new}}_{1, j}, \ldots, T^{\text{new}}_{M,j})^T$ by boundary treatment (\ref{h16}) for $j=1$ and $j=N$;\\
 \quad $\mathbf{T}^{n+1} = \mathbf{T}^{\text{new}}$. 
 \end{algorithmic}
 \end{algorithm}

 \begin{algorithm}
 \caption{Row-wise sweeping Newton-type method}\label{alg:row}
 \begin{algorithmic}
 \State $\mathbf{T}^{\text{old}} = \mathbf{T}^n;$ 
 \For {i=2 to M-1} \\
 \quad Compute $\mathbf{T^{\text{new}}_i} = (T_{i,2}, \cdots , T_{i,N-1})^T $ by Newton iteration (\ref{h14}); \\
 \quad Calculate $T^{\text{new}}_{i, 1}$ and $T^{\text{new}}_{i, N}$ by boundary treatment (\ref{h16});
 \EndFor \\
 \quad Compute $(T^{\text{new}}_{i, 1}, \ldots, T^{\text{new}}_{i,N})^T$ by boundary treatment (\ref{h16}) for $i=1$ and $i=M$;\\
 \quad $\mathbf{T}^{\text{old}} = \mathbf{T}^{\text{new}}$; 
 \For {i=M-1 to 2} \\
 \quad Compute $\mathbf{T^{\text{new}}_i} = (T_{i,2}, \cdots , T_{i,N-1})^T $ by Newton iteration (\ref{h14}); \\
 \quad Calculate $T^{\text{new}}_{i, 1}$ and $T^{\text{new}}_{i, N}$ by boundary treatment (\ref{h16});
 \EndFor \\
 \quad Compute $(T^{\text{new}}_{i, 1}, \ldots, T^{\text{new}}_{i,N})^T$ by boundary treatment (\ref{h16}) for $i=1$ and $i=M$;\\
 \quad $\mathbf{T}^{n+1} = \mathbf{T}^{\text{new}}$. 
 \end{algorithmic}
 \end{algorithm}

 \begin{algorithm}
 \caption{Column-row-wise sweeping Newton-type method}\label{alg:columnrow}
 \begin{algorithmic}
 \State $\mathbf{T}^{\text{old}} = \mathbf{T}^n;$ 
 \For {j=2 to N-1} \\
 \quad Compute $\mathbf{T}^{\text{new}}_{\mathbf{j}} = (T_{2, j}, T_{3, j}, \ldots, T_{M-1, j})^T$ by Newton iteration (\ref{h12}); \\
 \quad Calculate $T^{\text{new}}_{1, j}$ and $T^{\text{new}}_{M, j}$ by boundary treatment (\ref{h16});
 \EndFor \\
 \quad Compute $(T^{\text{new}}_{1, j}, \ldots, T^{\text{new}}_{M,j})^T$ by boundary treatment (\ref{h16}) for $j=1$ and $j=N$;\\
 \quad $\mathbf{T}^{\text{old}} = \mathbf{T}^{\text{new}}$;
 \For {j=N-1 to 2} \\
 \quad Compute $\mathbf{T}^{\text{new}}_{\mathbf{j}} = (T_{2, j}, T_{3, j}, \ldots, T_{M-1, j})^T$ by Newton iteration (\ref{h12}); \\
 \quad Calculate $T^{\text{new}}_{1, j}$ and $T^{\text{new}}_{M, j}$ by boundary treatment (\ref{h16});
 \EndFor \\
 \quad Compute $(T^{\text{new}}_{1, j}, \ldots, T^{\text{new}}_{M,j})^T$ by boundary treatment (\ref{h16}) for $j=1$ and $j=N$;\\
 \quad $\mathbf{T}^{\text{old}} = \mathbf{T}^{\text{new}}$;
 \For {i=2 to M-1} \\
 \quad Compute $\mathbf{T^{\text{new}}_i} = (T_{i,2}, \cdots , T_{i,N-1})^T $ by Newton iteration (\ref{h14}); \\
 \quad Calculate $T^{\text{new}}_{i, 1}$ and $T^{\text{new}}_{i, N}$ by boundary treatment (\ref{h16});
 \EndFor \\
 \quad Compute $(T^{\text{new}}_{i, 1}, ... T^{\text{new}}_{i,N})^T$ by boundary treatment (\ref{h16}) for $i=1$ and $i=M$;\\
 \quad $\mathbf{T}^{\text{old}} = \mathbf{T}^{\text{new}}$; 
 \For {i=M-1 to 2} \\
 \quad Compute $\mathbf{T^{\text{new}}_i} = (T_{i,2}, \cdots , T_{i,N-1})^T $ by Newton iteration (\ref{h14}); \\
 \quad Calculate $T^{\text{new}}_{i, 1}$ and $T^{\text{new}}_{i, N}$ by boundary treatment (\ref{h16});
 \EndFor \\
 \quad Compute $(T^{\text{new}}_{i, 1}, \ldots, T^{\text{new}}_{i,N})^T$ by boundary treatment (\ref{h16}) for $i=1$ and $i=M$;\\
 \quad $\mathbf{T}^{n+1} = \mathbf{T}^{\text{new}}$. 
 \end{algorithmic}
 \end{algorithm}

% \begin{description}
% \item [Column-wise sweep:] For $j-$th column, $\{ u_{2, j}, u_{3, j},... u_{M-1, j} \}$ will be solved by Newton iteration of the nonlinear system. The sweep is in two-directional $j = 2:N-1$ and $j=N-1:2:-1$.
% \item [Row-wise sweep:] For $i-$th row, $\{ u_{i, 2}, u_{i, 3},... u_{i, N-1} \}$ will be updated by Newton iteration of the nonlinear system. The sweep is in two-directional $i = 2:M-1$ and $i=M-1:2:-1$.
% \item [Column-row-wise: ] Column wise sweep first and then row wise.  The sweep in four directional $j = 2:N-1$, $j=N-1:2:-1$, $i = 2:M-1$ and $i=M-1:2:-1$.
% \end{description}

\newpage

\subsection{High-order schemes}
In order to achieve high-order accuracy, $T^{\pm}_x$ and $T^{\pm}_y$ in (\ref{h7}) will be approximated by high-order WENO reconstructions. In this section, we will present the WENO5, WENO7 and WENO9 reconstructions, respectively. 

\subsubsection{WENO5 reconstruction}
To approximate $(T^-_x)_{i, j}$, a left-biased stencil in the $x$-direction $\{T_{k, j}, k=i-3, ..., i+2\}$ will be used, from which three third-order approximations can be obtained,
\begin{align}\label{w1}
& (T_x)^{-,0}_{i, j} = \frac{1}{3} \frac{\triangle^+_x T_{i-3, j}}{\Delta x} -\frac{7}{6} \frac{\triangle^+_x T_{i-2, j}}{\Delta x} +\frac{11}{6}\frac{\triangle^+_x T_{i-1, j}}{\Delta x}; \\ \nonumber
& (T_x)^{-,1}_{i, j} = -\frac{1}{6} \frac{\triangle^+_x T_{i-2,j}}{\Delta x} + \frac{5}{6} \frac{\triangle^+_x T_{i-1,j}}{\Delta x} +\frac{1}{3} \frac{\triangle^+_x T_{i, j}}{\Delta x}; 
\\ \nonumber
& (T_x)^{-,2}_{i, j} = \frac{1}{3} \frac{\triangle^+_x T_{i-1,j}}{\Delta x} +\frac{5}{6} \frac{\triangle^+_x T_{i, j}}{\Delta x} -\frac{1}{6} \frac{\triangle^+_x T_{i+1, j}}{\Delta x}. 
\end{align}
Then, a fifth-order WENO approximation for $T_x(x_i, y_j)$ is a convex combination of $(T_x)^{-,0}_{i, j}$, $(T_x)^{-,1}_{i, j}$ and $(T_x)^{-,2}_{i, j}$, 
\begin{equation}\label{w2}
(T_x)^-_{i, j} = \omega_0 (T_x)^{-,0}_{i, j} +\omega_1 (T_x)^{-,1}_{i, j} + \omega_2 (T_x)^{-,2}_{i, j}, 
\end{equation}
where $\omega_k \ge 0$ are the nonlinear weights associated with the stencils given by
\begin{equation}\label{w3}
\omega_k = \frac{\alpha_k}{\sum^2_{k=0} \alpha_k}, \qquad \alpha_k = \frac{d_k}{(\epsilon_{\text{weno}} +\beta_k)^2}, \qquad k=0,1,2.
\end{equation}
The constants $d_k$ are optimal linear weights chosen as $d_0=0.1, d_1=0.6$ and $d_2=0.3$ such that $\sum^2_{k=0} d_k (T_x)^{-, k}_{i ,j}$ approximates $T_x(x_i, y_j)$ with full fifth-order accuracy. $\epsilon_{weno}$ is a small parameter added to avoid division by zero in the denominator. $\beta_k$ are called smoothness indicators which provide a measure of the smoothness of the solution over the corresponding stencils. We follow the definition of the classical smoothness indicators in \cite{jiang1996efficient}, which yields
\begin{align}\label{w4}
& \scalebox{0.92}{$\displaystyle \beta_0 = \frac{13}{12} \left( \frac{\triangle^+_x T_{i-3, j}}{\Delta x}-2\frac{\triangle^+_x T_{i-2, j}}{\Delta x}+\frac{\triangle^+_x T_{i-1,j}}{\Delta x} \right)^2 +\frac{1}{4} \left (\frac{\triangle^+_x T_{i-3,j}}{\Delta x}-4\frac{\triangle^+_x T_{i-2, j}}{\Delta x}+3\frac{\triangle^+_x T_{i-1, j}}{\Delta x} \right)^2$}, \\ \nonumber
& \beta_1 = \frac{13}{12} \left(\frac{\triangle^+_x T_{i-2,j}}{\Delta x} -2\frac{\triangle^+_x T_{i-1, j}}{\Delta x}+\frac{\triangle^+_x T_{i, j}}{\Delta x}\right)^2 + \frac{1}{4} \left( \frac{\triangle^+_x T_{i,j}}{\Delta x} -\frac{\triangle^+_x T_{i-2, j}}{\Delta x} \right)^2, \\ \nonumber
& \scalebox{0.95}{$\displaystyle \beta_2 = \frac{13}{12} \left(\frac{\triangle^+_x T_{i-1, j}}{\Delta x}-2\frac{\triangle^+_x T_{i, j}}{\Delta x}+\frac{\triangle^+_x T_{i+1, j}}{\Delta x} \right)^2 +\frac{1}{4} \left( 3\frac{\triangle^+_x T_{i-1, j}}{\Delta x}-4\frac{\triangle^+_x T_{i, j}}{\Delta x}+\frac{\triangle^+_x T_{i+1, j}}{\Delta x} \right)^2.$}
\end{align}
Similarly, $(T_x)^+_{i,j}$ can be reconstructed by a right-biased stencil $\{T_{k, j}, k=i-2, ..., i+3\}$, and the definitions for $(T_y)^{\pm}_{i,j}$ are analogous. 

\subsubsection{WENO7 reconstruction}
We present the WENO7 reconstruction for $(T_x)^-_{i,j}$ with a left-biased stencil $\{T_{k,j}, k=i-4, ..., i+3 \}$, with the understanding that the reconstructions of $(T_x)^+_{i,j}$ and $(T_y)^{\pm}_{i,j}$ are analogous. With the selection of small stencils, we can calculate four fourth-order approximations,
\begin{align}\label{w5}
& (T_x)^{-,0}_{i,j} = -\frac{1}{4} \frac{\triangle^+_x T_{i-4, j}}{\Delta x} +\frac{13}{12} \frac{\triangle^+_x T_{i-3, j}}{\Delta x} -\frac{23}{12} \frac{\triangle^+_x T_{i-2, j}}{\Delta x} +\frac{25}{12} \frac{\triangle^+_x T_{i-1, j}}{\Delta x}, \\ \nonumber
& (T_x)^{-,1}_{i,j} =\frac{1}{12} \frac{\triangle^+_x T_{i-3, j}}{\Delta x} -\frac{5}{12} \frac{\triangle^+_x T_{i-2, j}}{\Delta x} +\frac{13}{12} \frac{\triangle^+_x T_{i-1, j}}{\Delta x} +\frac{1}{4} \frac{\triangle^+_x T_{i, j}}{\Delta x}, \\ \nonumber
& (T_x)^{-,2}_{i,j} = -\frac{1}{12} \frac{\triangle^+_x T_{i-2, j}}{\Delta x} +\frac{7}{12} \frac{\triangle^+_x T_{i-1, j}}{\Delta x} +\frac{7}{12} \frac{\triangle^+_x T_{i, j}}{\Delta x} -\frac{1}{12} \frac{\triangle^+_x T_{i+1, j}}{\Delta x}, \\ \nonumber
& (T_x)^{-,3}_{i,j} = \frac{1}{4} \frac{\triangle^+_x T_{i-1, j}}{\Delta x} +\frac{13}{12}\frac{\triangle^+_x T_{i, j}}{\Delta x} -\frac{5}{12} \frac{\triangle^+_x T_{i+1, j}}{\Delta x} +\frac{1}{12} \frac{\triangle^+_x T_{i+2, j}}{\Delta x}, 
\end{align}
which will be used to construct a seventh-order approximation 
\begin{equation}\label{w5.5}
(T_x)^-_{i,j} = \sum^3_{k=0} \omega_k (T_x)^{-,k}_{i,j}. 
\end{equation}
The smoothness indicators are given by \cite{balsara2000monotonicity}
\begin{align}\label{w6}
\beta_0 ={}& \frac{\triangle^+_x T_{i-4,j}}{\Delta x}\left(547\frac{\triangle^+_x T_{i-4,j}}{\Delta x}-3882\frac{\triangle^+_x T_{i-3,j}}{\Delta x}+4642\frac{\triangle^+_x T_{i-2,j}}{\Delta x}-1854\frac{\triangle^+_x T_{i-1,j}}{\Delta x}\right) \nonumber\\
&+\frac{\triangle^+_x T_{i-3,j}}{\Delta x}\left(7043\frac{\triangle^+_x T_{i-3,j}}{\Delta x}-17246\frac{\triangle^+_x T_{i-2,j}}{\Delta x}+7042\frac{\triangle^+_x T_{i-1,j}}{\Delta x}\right) \nonumber\\
&+\frac{\triangle^+_x T_{i-2,j}}{\Delta x}\left(11003\frac{\triangle^+_x T_{i-2,j}}{\Delta x}-9402\frac{\triangle^+_x T_{i-1,j}}{\Delta x}\right)+2107\left(\frac{\triangle^+_x T_{i-1,j}}{\Delta x}\right)^2, \nonumber\\[0.4em]
\beta_1 ={}& \frac{\triangle^+_x T_{i-3,j}}{\Delta x}\left(267\frac{\triangle^+_x T_{i-3,j}}{\Delta x}-1642\frac{\triangle^+_x T_{i-2,j}}{\Delta x}+1602\frac{\triangle^+_x T_{i-1,j}}{\Delta x}-494\frac{\triangle^+_x T_{i,j}}{\Delta x}\right) \nonumber\\
&+\frac{\triangle^+_x T_{i-2,j}}{\Delta x}\left(2843\frac{\triangle^+_x T_{i-2,j}}{\Delta x}-5966\frac{\triangle^+_x T_{i-1,j}}{\Delta x}+1922\frac{\triangle^+_x T_{i,j}}{\Delta x}\right) \nonumber\\
&+\frac{\triangle^+_x T_{i-1,j}}{\Delta x}\left(3443\frac{\triangle^+_x T_{i-1,j}}{\Delta x}-2522\frac{\triangle^+_x T_{i,j}}{\Delta x}\right)+547\left(\frac{\triangle^+_x T_{i,j}}{\Delta x}\right)^2, 
\end{align}
\begin{align}
\beta_2 ={}& \frac{\triangle^+_x T_{i-2,j}}{\Delta x}\left(547\frac{\triangle^+_x T_{i-2,j}}{\Delta x}-2522\frac{\triangle^+_x T_{i-1,j}}{\Delta x}+1922\frac{\triangle^+_x T_{i,j}}{\Delta x}-494\frac{\triangle^+_x T_{i+1,j}}{\Delta x}\right) \nonumber\\
&+\frac{\triangle^+_x T_{i-1,j}}{\Delta x}\left(3433\frac{\triangle^+_x T_{i-1,j}}{\Delta x}-5966\frac{\triangle^+_x T_{i,j}}{\Delta x}+1602\frac{\triangle^+_x T_{i+1,j}}{\Delta x}\right) \nonumber\\
&+\frac{\triangle^+_x T_{i,j}}{\Delta x}\left(2843\frac{\triangle^+_x T_{i,j}}{\Delta x}-1642\frac{\triangle^+_x T_{i+1,j}}{\Delta x}\right)+267\left(\frac{\triangle^+_x T_{i+1,j}}{\Delta x}\right)^2, \nonumber\\[0.4em]
\beta_3 ={}& \frac{\triangle^+_x T_{i-1,j}}{\Delta x}\left(2107\frac{\triangle^+_x T_{i-1,j}}{\Delta x}-9402\frac{\triangle^+_x T_{i,j}}{\Delta x}+7042\frac{\triangle^+_x T_{i+1,j}}{\Delta x}-1854\frac{\triangle^+_x T_{i+2,j}}{\Delta x}\right) \nonumber\\
&+\frac{\triangle^+_x T_{i,j}}{\Delta x}\left(11003\frac{\triangle^+_x T_{i,j}}{\Delta x}-17246\frac{\triangle^+_x T_{i+1,j}}{\Delta x}+4642\frac{\triangle^+_x T_{i+2,j}}{\Delta x}\right) \nonumber\\
&+\frac{\triangle^+_x T_{i+1,j}}{\Delta x}\left(7043\frac{\triangle^+_x T_{i+1,j}}{\Delta x}-3882\frac{\triangle^+_x T_{i+2,j}}{\Delta x}\right)+547\left(\frac{\triangle^+_x T_{i+2,j}}{\Delta x}\right)^2, 
\end{align}
and the optimal linear weights are $d_0 = \frac{1}{35}, d_1 = \frac{12}{35}, d_2 = \frac{18}{35}$ and $d_3 = \frac{4}{35}$.

\subsubsection{WENO9 reconstruction}
Similar to the previous sections, we present the WENO9 reconstruction of $(T_x)^-_{i,j}$ on the left-biased stencil $\{T_{k,j}, k=i-5, ..., i+4 \}$, where the lower-order approximations on small stencils are
\begin{align}\label{w7} \nonumber
& (T_x)^{-,0}_{i,j} =\frac{1}{5} \frac{\triangle^+_x T_{i-5,j}}{\Delta x}-\frac{21}{20} \frac{\triangle^+_x T_{i-4,j}}{\Delta x} +\frac{137}{60} \frac{\triangle^+_x T_{i-3,j}}{\Delta x} - \frac{163}{60} \frac{\triangle^+_x T_{i-2,j}}{\Delta x} +\frac{137}{60} \frac{\triangle^+_x T_{i-1,j}}{\Delta x}, \\ \nonumber
& (T_x)^{-,1}_{i,j} = -\frac{1}{20} \frac{\triangle^+_x T_{i-4,j}}{\Delta x} +\frac{17}{60} \frac{\triangle^+_x T_{i-3,j}}{\Delta x}-\frac{43}{60} \frac{\triangle^+_x T_{i-2,j}}{\Delta x} +\frac{77}{60} \frac{\triangle^+_x T_{i-1,j}}{\Delta x}+\frac{1}{5}\frac{\triangle^+_x T_{i,j}}{\Delta x}, \\ \nonumber
& (T_x)^{-,2}_{i,j} = \frac{1}{30}\frac{\triangle^+_x T_{i-3,j}}{\Delta x} -\frac{13}{60}\frac{\triangle^+_x T_{i-2,j}}{\Delta x}+\frac{47}{60}\frac{\triangle^+_x T_{i-1,j}}{\Delta x} +\frac{9}{20} \frac{\triangle^+_x T_{i,j}}{\Delta x}-\frac{1}{20} \frac{\triangle^+_x T_{i+1,j}}{\Delta x}, \\ \nonumber
& (T_x)^{-,3}_{i,j} = -\frac{1}{20} \frac{\triangle^+_x T_{i-2,j}}{\Delta x} +\frac{9}{20} \frac{\triangle^+_x T_{i-1,j}}{\Delta x} +\frac{47}{60} \frac{\triangle^+_x T_{i,j}}{\Delta x} -\frac{13}{60} \frac{\triangle^+_x T_{i+1,j}}{\Delta x} +\frac{1}{30} \frac{\triangle^+_x T_{i+2,j}}{\Delta x}, \\ \nonumber
& (T_x)^{-,4}_{i,j} = \frac{1}{5} \frac{\triangle^+_x T_{i-1,j}}{\Delta x} +\frac{77}{60} \frac{\triangle^+_x T_{i,j}}{\Delta x} -\frac{43}{60} \frac{\triangle^+_x T_{i+1,j}}{\Delta x} +\frac{17}{60} \frac{\triangle^+_x T_{i+2,j}}{\Delta x}-\frac{1}{20} \frac{\triangle^+_x T_{i+3,j}}{\Delta x}, 
\end{align}
and then
\begin{equation}\label{w8}
(T_x)^-_{i,j} = \sum^4_{k=0} \omega_k (T_x)^{-,k}_{i,j}.
\end{equation}
The optimal linear weights are $d_0 = \frac{1}{126}, d_1 = \frac{10}{63}, d_2 = \frac{10}{21}, d_3 = \frac{20}{63}$ and $d_4 = \frac{5}{126}$. The smoothness indicators are

\begin{align*}
%\label{w9}\nonumber
\beta_0 & = 22658\left(\frac{\triangle^+_x T_{i-5, j}}{\Delta x}\right)^2 +482963 \left(\frac{\triangle^+_x T_{i-4,j}}{\Delta x}\right)^2 +1521393\left(\frac{\triangle^+_x T_{i-3,j}}{\Delta x}\right)^2 \\ \nonumber & +1020563 \left(\frac{\triangle^+_x T_{i-2,j}}{\Delta x}\right)^2  + 107918\left(\frac{\triangle^+_x T_{i-1,j}}{\Delta x}\right)^2 -2462076\frac{\triangle^+_x T_{i-3,j}}{\Delta x}\frac{\triangle^+_x T_{i-2,j}}{\Delta x} \\ \nonumber & + 758823 \frac{\triangle^+_x T_{i-3,j}}{\Delta x} \frac{\triangle^+_x T_{i-1,j}}{\Delta x} 
 -649501 \frac{\triangle^+_x T_{i-2,j}}{\Delta x}\frac{\triangle^+_x T_{i-1,j}}{\Delta x} \\ \nonumber & + \frac{\triangle^+_x T_{i-4,j}}{\Delta x}\left(-1704396 \frac{\triangle^+_x T_{i-3,j}}{\Delta x}+1358458 \frac{\triangle^+_x T_{i-2,j}}{\Delta x}-411487 \frac{\triangle^+_x T_{i-1,j}}{\Delta x} \right) \\ \nonumber 
& + \frac{\triangle^+_x T_{i-5,j}}{\Delta x}\left(-208501 \frac{\triangle^+_x T_{i-4,j}}{\Delta x}+364863 \frac{\triangle^+_x T_{i-3,j}}{\Delta x}-288007 \frac{\triangle^+_x T_{i-2,j}}{\Delta x}\right.\\ \nonumber 
  & \left.+86329\frac{\triangle^+_x T_{i-1,j}}{\Delta x}\right),
  \end{align*}
  %\vspace{-1.3em}
  \begin{equation*}
\scalebox{0.95}{$\begin{aligned}
\beta_1 & = 6908 \left(\frac{\triangle^+_x T_{i-4,j}}{\Delta x}\right)^2+138563 \left(\frac{\triangle^+_x T_{i-3,j}}{\Delta x}\right)^2+406293 \left(\frac{\triangle^+_x T_{i-2,j}}{\Delta x}\right)^2\\
  & +242723 \left(\frac{\triangle^+_x T_{i-1,j}}{\Delta x}\right)^2 +22658 \left(\frac{\triangle^+_x T_{i,j}}{\Delta x}\right)^2-611976 \frac{\triangle^+_x T_{i-2,j}}{\Delta x} \frac{\triangle^+_x T_{i-1,j}}{\Delta x} \\
  & +165153\frac{\triangle^+_x T_{i-2,j}}{\Delta x}\frac{\triangle^+_x T_{i,j}}{\Delta x} -140251 \frac{\triangle^+_x T_{i-1,j}}{\Delta x} \frac{\triangle^+_x T_{i,j}}{\Delta x}  \\
  & + \frac{\triangle^+_x T_{i-3,j}}{\Delta x}\left(-464976 \frac{\triangle^+_x T_{i-2,j}}{\Delta x} +337018\frac{\triangle^+_x T_{i-1,j}}{\Delta x}-88297 \frac{\triangle^+_x T_{i,j}}{\Delta x}\right) \\
  & + \frac{\triangle^+_x T_{i-4,j}}{\Delta x}\left(-60871 \frac{\triangle^+_x T_{i-3,j}}{\Delta x}+99213 \frac{\triangle^+_x T_{i-2,j}}{\Delta x}-70237\frac{\triangle^+_x T_{i-1,j}}{\Delta x}+18079 \frac{\triangle^+_x T_{i,j}}{\Delta x}\right),
\end{aligned}$}
\end{equation*}
\vspace{0.8em}
\begin{equation*}
\scalebox{0.95}{$\begin{aligned}
\beta_2 & = 6908\left(\frac{\triangle^+_x T_{i-3,j}}{\Delta x}\right)^2+104963\left(\frac{\triangle^+_x T_{i-2,j}}{\Delta x}\right)^2 +231153 \left(\frac{\triangle^+_x T_{i-1,j}}{\Delta x}\right)^2 \\
  & +104963 \left(\frac{\triangle^+_x T_{i,j}}{\Delta x}\right)^2  +6908 \left(\frac{\triangle^+_x T_{i+1,j}}{\Delta x}\right)^2  -299076 \frac{\triangle^+_x T_{i-1,j}}{\Delta x}\frac{\triangle^+_x T_{i,j}}{\Delta x} \\
  &+67923 \frac{\triangle^+_x T_{i-1,j}}{\Delta x} \frac{\triangle^+_x T_{i+1,j}}{\Delta x}-51001 \frac{\triangle^+_x T_{i,j}}{\Delta x}\frac{\triangle^+_x T_{i+1,j}}{\Delta x} \\ & + \frac{\triangle^+_x T_{i-2,j}}{\Delta x}\left(-299076 \frac{\triangle^+_x T_{i-1,j}}{\Delta x}+179098 \frac{\triangle^+_x T_{i, j}}{\Delta x}-38947 \frac{\triangle^+_x T_{i+1, j}}{\Delta x} \right)\\ & +\frac{\triangle^+_x T_{i-3, j}}{\Delta x}\left(-51001 \frac{\triangle^+_x T_{i-2,j}}{\Delta x}+67923\frac{\triangle^+_x T_{i-1, j}}{\Delta x}-38947 \frac{\triangle^+_x T_{i, j}}{\Delta x}+8209\frac{\triangle^+_x T_{i+1,j}}{\Delta x}\right),
\end{aligned}$}
\end{equation*}
 % \vspace{-1.3em}
  \begin{equation*}
\scalebox{0.93}{$\begin{aligned}
\beta_3 & = 22658 \left(\frac{\triangle^+_x T_{i-2, j}}{\Delta x}\right)^2+242723 \left(\frac{\triangle^+_x T_{i-1, j}}{\Delta x}\right)^2+406293 \left(\frac{\triangle^+_x T_{i, j}}{\Delta x}\right)^2 \\ & +138563\left(\frac{\triangle^+_x T_{i+1, j}}{\Delta x}\right)^2
+6908\left(\frac{\triangle^+_x T_{i+2, j}}{\Delta x}\right)^2 -464976 \frac{\triangle^+_x T_{i, j}}{\Delta x} \frac{\triangle^+_x T_{i+1, j}}{\Delta x} \\ & +99213 \frac{\triangle^+_x T_{i, j}}{\Delta x}\frac{\triangle^+_x T_{i+2, j}}{\Delta x}-60871 \frac{\triangle^+_x T_{i+1, j}}{\Delta x} \frac{\triangle^+_x T_{i+2, j}}{\Delta x} \\
&+\frac{\triangle^+_x T_{i-1,j}}{\Delta x}\left(-611976 \frac{\triangle^+_x T_{i, j}}{\Delta x}+337018 \frac{\triangle^+_x T_{i+1,j}}{\Delta x}-70237\frac{\triangle^+_x T_{i+2,j}}{\Delta x}\right) \\ 
& + \frac{\triangle^+_x T_{i-2,j}}{\Delta x}\left(-140251 \frac{\triangle^+_x T_{i-1,j}}{\Delta x}+165153 \frac{\triangle^+_x T_{i,j}}{\Delta x}-88297 \frac{\triangle^+_x T_{i+1,j}}{\Delta x}+18079 \frac{\triangle^+_x T_{i+2,j}}{\Delta x}\right), 
\end{aligned}$}
\end{equation*}
%\vspace{-1.3em}
\begin{equation*}
\scalebox{0.92}{$\begin{aligned}
%\label{w11}\nonumber
\beta_4 & = 107918 \left(\frac{\triangle^+_x T_{i-1,j}}{\Delta x}\right)^2+1020563 \left(\frac{\triangle^+_x T_{i, j}}{\Delta x}\right)^2+1521393 \left(\frac{\triangle^+_x T_{i+1,j}}{\Delta x}\right)^2 \\ &+482963 \left(\frac{\triangle^+_x T_{i+2,j}}{\Delta x}\right)^2  +22658 \left(\frac{\triangle^+_x T_{i+3,j}}{\Delta x}\right)^2-1704396 \frac{\triangle^+_x T_{i+1,j}}{\Delta x} \frac{\triangle^+_x T_{i+2,j}}{\Delta x}\\ & +364863 \frac{\triangle^+_x T_{i+1,j}}{\Delta x} \frac{\triangle^+_x T_{i+3,j}}{\Delta x} -208501\frac{\triangle^+_x T_{i+2,j}}{\Delta x} \frac{\triangle^+_x T_{i+3,j}}{\Delta x} \\ & +\frac{\triangle^+_x T_{i,j}}{\Delta x}\left(-2462076 \frac{\triangle^+_x T_{i+1,j}}{\Delta x}+1358458\frac{\triangle^+_x T_{i+2,j}}{\Delta x}-288007\frac{\triangle^+_x T_{i+3,j}}{\Delta x}\right)\\ 
&+ \frac{\triangle^+_x T_{i-1,j}}{\Delta x}\left(-649501 \frac{\triangle^+_x T_{i,j}}{\Delta x}+758823 \frac{\triangle^+_x T_{i+1,j}}{\Delta x}-411487 \frac{\triangle^+_x T_{i+2,j}}{\Delta x}+86329 \frac{\triangle^+_x T_{i+3,j}}{\Delta x}\right).
\end{aligned}$}
\end{equation*}

%\newpage
% ===== END algorithm.tex =====

% section 3
% ===== BEGIN numerical_staticHJ.tex =====

%\pagebreak
\section{Numerical examples}
\label{numer}

In this section, we conduct numerical experiments to validate the convergence of the proposed high-order implicit Newton-type WENO schemes. For all the high-order methods, the first-order point-wise Gaussian sweeping method \cite{kao2004lax} is used to generate the initial conditions. The iterative sweeping is stopped with the convergence criterion $||T^{n+1}-T^n|| \le \epsilon$. Moreover, the small parameter in the WENO nonlinear weights  is denoted as $\epsilon_{\text{WENO}}$. For the boundary treatment, a first-order linear interpolation is adopted. The numerical errors in the examples are measured in the interior region with a margin of 0.2 from the computational domain. 
%%%%%%%%%%%%%%%%%%%%%%%%%%%%%%%%%%%%%%%%%%%%%%t
%\subsection{Static HJ equation}
\subsection{Isotropic Eikonal equation}

\begin{exam}\label{2d_static_eikonal_1}
In this example, we first consider the reduced case (\ref{geikonal}) when there is no background moving medium, i.e., $v_1=v_2=0.0$. The computational domain is chosen as $\Omega = [-1, 1] \times [-1, 1]$, $\Gamma = \{(0, 0)\}$, and $g(0, 0) = -2$. $F(\textbf{x})$ is computed by assuming the exact solution is 
\begin{equation}\label{eq10}
\text{T}(x,y) = -\cos(\frac{\pi}{2}x) - \cos(\frac{\pi}{2}y).
\end{equation}
The values around $\Gamma$ are assigned based on the expression of the exact solution.  
\end{exam}
It is known that there are only isotropic waves. We solve this example by the proposed high-order implicit Newton-type WENO schemes, taking $\alpha = 1.01$, and the convergence criterion is chosen as  $\epsilon=10^{-8}$. The numerical errors and orders are given in Tables \ref{2d_ex1_table1} - \ref{2d_ex1_table3}. It is easy to observe that the proposed Newton-type WENO schemes are all convergent. Moreover, the desired high order of accuracy are achieved with fifth-, seventh- and ninth-order WENO reconstructions with all the three sweeping strategies. Moreover, we can see that the numerical errors are almost the same for column-row-wise, column-wise and row-wise sweep approaches. The iteration numbers of the one-directional sweep are identical which is about half of those of the two-directional sweep. This matches our expectation as the simulation setup and the exact solution are symmetric.

\begin{table}[H]
\begin{center}
\caption {Accuracy study for Example \ref{2d_static_eikonal_1}. Newton-type WENO5 method. $\epsilon = 10^{-8}$, $\epsilon_{\text{WENO}} = 10^{-9}$, $\alpha=1.01$.  }
\begin{tabular}{|c|c|c c|cc|cc|c|}
%\cline{3-8}
\hline
 & N  & $L^1$ error & order   & $L^2$ error  & order  & $L^{\infty}$  & order & sweep \\\hline
\multicolumn{1}{|c|}{\multirow{4}{*}{Column-row-wise}}
& 60 & 1.05E-07&	--&	2.85E-09	&--&	1.66E-06&	--&48
\\\cline{2-9}  
& 80 &2.14E-08	&5.51	&3.65E-10	&7.15	&1.45E-07	&8.46&54
\\\cline{2-9}  
& 100 & 7.27E-09	&4.84	&9.72E-11	&5.93	&2.04E-08	&8.80&61
\\\cline{2-9}  
& 120 & 2.95E-09	&4.95	&3.27E-11	&5.98	&7.33E-09	&5.62&67
\\\cline{2-9}  
& 140 & 1.38E-09	&4.92	&1.31E-11	&5.93	&3.45E-09	&4.88&74
\\\cline{2-9}  
& 160 & 7.14E-10	&4.94	&5.92E-12	&5.95	&1.79E-09	&4.94&81
\\\cline{2-9}  
& 180 & 3.99E-10	&4.95	&2.93E-12	&5.96	&9.97E-10	&4.96&88
\\\hline
\multicolumn{1}{|c|}{\multirow{4}{*}{Column-wise}}
& 60 & 1.05E-07	&--	&2.85E-09	&--	&1.65E-06	&--&108
\\\cline{2-9}  
& 80 & 2.14E-08	&5.51	&3.65E-10	&7.14	&1.45E-07	&8.45&124
\\\cline{2-9}  
& 100 & 7.28E-09	&4.84	&9.73E-11	&5.92	&2.11E-08	&8.64&139
\\\cline{2-9}  
& 120 & 2.95E-09	&4.95	&3.27E-11	&5.98	&7.89E-09	&5.38&153
\\\cline{2-9}  
& 140 &1.38E-09	&4.92	&1.31E-11	&5.93	&4.25E-09	&4.01&166
\\\cline{2-9}  
& 160 & 7.14E-10	&4.94	&5.92E-12	&5.96	&1.93E-09	&5.93&182
\\\cline{2-9}  
& 180 & 3.99E-10	&4.94	&2.94E-12	&5.95	&1.07E-09	&4.98&198
\\\hline
\multicolumn{1}{|c|}{\multirow{4}{*}{Row-wise}}
& 60 & 1.05E-07	&--	&2.85E-09	&--	&1.65E-06	&--&108
\\\cline{2-9}  
& 80 & 2.14E-08	&5.51	&3.65E-10	&7.14	&1.45E-07	&8.45&124
\\\cline{2-9}  
& 100 & 7.28E-09	&4.84	&9.73E-11	&5.92	&2.11E-08	&8.64&139
\\\cline{2-9}  
& 120 & 2.95E-09	&4.95	&3.27E-11	&5.98	&7.89E-09	&5.38&153
\\\cline{2-9}  
& 140 &1.38E-09	&4.92	&1.31E-11	&5.93	&4.25E-09	&4.01&166
\\\cline{2-9}  
& 160 & 7.14E-10	&4.94	&5.92E-12	&5.96	&1.93E-09	&5.93&182
\\\cline{2-9}  
& 180 & 3.99E-10	&4.94	&2.94E-12	&5.95	&1.07E-09	&4.98&198
\\\hline
\end{tabular}
\label{2d_ex1_table1}
\end{center}
\end{table}

\begin{table}[H]
\begin{center}
\caption {Accuracy study for Example \ref{2d_static_eikonal_1}. Newton-type WENO7 method. $\epsilon = 10^{-8}$, $\epsilon_{\text{WENO}} = 10^{-9}$, $\alpha=1.01$. }
\begin{tabular}{|c|c|c c|cc|cc|c|}
%\cline{3-8}
\hline
 & N  & $L^1$ error & order   & $L^2$ error  & order  & $L^{\infty}$  & order & sweep \\\hline
\multicolumn{1}{|c|}{\multirow{4}{*}{ Column-row-wise }}
& 40 & 3.18E-06	&--&	2.81E-07&--&		1.04E-04&--&46	
\\\cline{2-9}  
& 60 & 3.63E-08	&11.0	&2.02E-09	&12.2	&1.10E-06	&11.2&51
\\\cline{2-9}  
& 80 & 2.21E-09	&9.72	&1.26E-10	&9.65	&1.48E-07	&6.97&54
\\\cline{2-9}  
& 100 &1.41E-10	&12.4	&4.76E-12	&14.7	&9.47E-09	&12.3&59
\\\cline{2-9}  
& 120 & 2.97E-11	&8.53	&6.51E-13	&10.9	&1.04E-09	&12.1&64
%\\\cline{2-9}  
%& 140 & 8.71E-12	&7.96	&1.55E-13	&9.29	&2.58E-10	&9.03&70
\\\hline
\multicolumn{1}{|c|}{\multirow{4}{*}{Column-wise }}
& 40 & 3.18E-06&	--	&2.81E-07	&--	&1.04E-04	&--&101
\\\cline{2-9}  
& 60 &3.66E-08	&11.0	&2.03E-09	&12.2	&1.09E-06	&11.2&103
\\\cline{2-9}  
& 80 & 2.42E-09	&9.44	&1.26E-10	&9.67	&1.47E-07	&6.96&107
\\\cline{2-9}  
& 100 & 1.87E-10	&11.5	&5.09E-12	&14.4	&8.88E-09	&12.6&119
\\\cline{2-9}  
& 120 & 2.68E-11	&10.7	&6.57E-13	&11.2	&1.82E-09	&8.69&135
\\\hline
\multicolumn{1}{|c|}{\multirow{4}{*}{Row-wise  }}
& 40 & 3.18E-06&	--	&2.81E-07	&--	&1.04E-04	&--&101
\\\cline{2-9}  
& 60 &3.66E-08	&11.0	&2.03E-09	&12.2	&1.09E-06	&11.2&103
\\\cline{2-9}  
& 80 & 2.42E-09	&9.44	&1.26E-10	&9.67	&1.47E-07	&6.96&107
\\\cline{2-9}  
& 100 & 1.87E-10	&11.5	&5.09E-12	&14.4	&8.88E-09	&12.6&119
\\\cline{2-9}  
& 120 & 2.68E-11	&10.7	&6.57E-13	&11.2	&1.82E-09	&8.69&135
\\\hline
\end{tabular}
\label{2d_ex1_table2}
\end{center}
\end{table}

\begin{table}[H]
\begin{center}
\caption {Accuracy study for Example \ref{2d_static_eikonal_1}. Newton-type WENO9 method. $\epsilon = 10^{-8}$, $\epsilon_{\text{WENO}} = 10^{-9}$, $\alpha=1.01$.  }
\begin{tabular}{|c|c|c c|cc|cc|c|}
%\cline{3-8}
\hline
 & N  & $L^1$ error & order   & $L^2$ error  & order  & $L^{\infty}$  & order & sweep \\\hline
\multicolumn{1}{|c|}{\multirow{4}{*}{ Column-row-wise }}
%& 20 & 2.60E-06&--&		2.08E-07&--&		7.59E-05&--	&55
%\\\cline{2-9}  
& 40 & 7.08E-08	&8.88	&3.92E-09	&9.80	&1.22E-06	&10.2&59
\\\cline{2-9}  
& 60 & 3.43E-09	&10.5	&1.55E-10	&11.2	&1.29E-07	&7.81&61
\\\cline{2-9}  
& 80 & 3.26E-10	&10.5	&1.12E-11	&11.8	&1.27E-08	&10.4&65
\\\cline{2-9}  
& 100 & 1.24E-10	&5.30	&2.52E-12	&8.19	&2.39E-09	&9.15&66
% \\\cline{2-9}  
% & 120 & 4.28E-11	&6.91E+00	&7.72E-13	&7.67E+00	&9.44E-10	&6.03E+00&72
% \\\cline{2-9}  
% & 140 & 1.29E-11	&8.99E+00	&2.09E-13	&9.77E+00	&2.99E-10	&8.61E+00&79
\\\hline
\multicolumn{1}{|c|}{\multirow{4}{*}{Column-wise }}
& 40 &2.60E-06	&--	&2.08E-07	&--	&7.59E-05	&--&116
\\\cline{2-9}  
& 60 & 7.11E-08&	8.87	&3.93E-09	&9.80	&1.22E-06	&10.2&121
\\\cline{2-9}  
& 80 & 3.71E-09	&10.3	&1.55E-10	&11.2	&1.31E-07	&7.77&123
\\\cline{2-9}  
& 100 & 6.58E-10	&7.75	&1.52E-11	&10.4	&1.44E-08	&9.88&125
\\\hline
\multicolumn{1}{|c|}{\multirow{4}{*}{Row-wise  }}
& 40 & 2.60E-06	&--	&2.08E-07	&--	&7.59E-05	&--&115
\\\cline{2-9}  
& 60 & 7.11E-08&	8.87	&3.93E-09	&9.80	&1.22E-06	&10.2&121
\\\cline{2-9}  
& 80 &3.71E-09	&10.3	&1.55E-10	&11.2	&1.31E-07	&7.77&123
\\\cline{2-9}  
& 100 & 6.58E-10	&7.75	&1.52E-11	&10.4	&1.44E-08	&9.88&125
\\\hline
\end{tabular}
\label{2d_ex1_table3}
\end{center}
\end{table}

%-----------------------------------------------------------------
\newpage
\newpage
\begin{exam}\label{2d_static_eikonal_2}
In this example, we consider the two-dimensional point source distance function problem. For the general static Eikonal equation (\ref{geikonal}), the velocity of the moving fluid is set to be $v_1 = v_2 =0$ on the computational domain  $\Omega = [-1, 1] \times [-1, 1]$,  $\Gamma = \{(0,0) \}$, $f(x,y) = 1$, and $g|_{\Gamma} = 0$. The exact solution can be calculated as
\begin{equation}\label{ex2_1}
\text{T}(x, y) = \sqrt{x^2+y^2}.
\end{equation}
The values on the boundary cells in the domain $[-0.1, 0.1] \times [-0.1, 0.1]$ are assigned based on the exact solution. 
\end{exam}
We simulate the problem using our proposed high-order implicit Newton-type WENO schemes with the column-row-wise sweep approach. The convergence criterion tolerance is $\epsilon=10^{-9}$ and $\alpha =2$. The WENO parameters are $\epsilon_{\text{WENO}} = 10^{-3}$ for WENO5 and $\epsilon_{\text{WENO}} = 10^{-9}$ for WENO7 and WENO9. The numerical errors and orders for the high-order implicit Newton-type WENO schemes are presented in Table \ref{2d_ex2_table1}, in which we can clearly see that the designed high orders of accuracy are achieved. Meanwhile, we observe that the WENO5 scheme needs more refined mesh to achieve the convergence and designed order of accuracy compared with seventh- and ninth-order schemes.  
%++++++++++++++++++++++++++++++++++++++++++++++++++

\begin{table}[H]
\begin{center}
\caption {Accuracy study for Example \ref{2d_static_eikonal_2}. Column-row-wise sweep.  }
\begin{tabular}{|c|c|c c|cc|cc|c|}
%\cline{3-8}
\hline
 & N  & $L^1$ error & order   & $L^2$ error  & order  & $L^{\infty}$  & order & sweep \\\hline
\multicolumn{1}{|c|}{\multirow{4}{*}{ WENO5 }}
%& 300 &1.03E-07&	--	&4.37E-10	&--	&2.09E-07	&-- &107
%\\\cline{2-9}
%& 320 & 4.87E-08	&11.5	&2.01E-10	&12.1	&9.89E-08	&11.6&112 
& 320 & 4.87E-08	&--	&2.01E-10	&--	&9.89E-08	&-- &112 
\\\cline{2-9}  
& 340 &2.52E-08	&10.9	&1.03E-10	&11.0	&6.30E-08	&7.46&115 
\\\cline{2-9}  
& 360 &1.48E-08	&9.33	&6.06E-11	&9.34	&4.46E-08	&6.02&121
\\\cline{2-9}  
& 400 & 6.67E-09	&7.56	&2.65E-11	&7.86	&2.45E-08	&5.68&130  
\\\cline{2-9}  
& 500 &1.84E-09	&5.78	&5.85E-12	&6.77	&7.13E-09	&5.54&155 
\\\cline{2-9}  
& 640 & 4.95E-10	&5.31	&1.18E-12	&6.48	&1.83E-09	&5.51&187
\\\hline
\multicolumn{1}{|c|}{\multirow{4}{*}{WENO7  }}
%& 80 & -- & --  & -- & --  &  --  &-- & --  
%\\\cline{2-9}  
%& 100 &  -- & --  & -- & --  &  --  &-- & --  
%\\\cline{2-9}  
& 120 &9.35E-05	&5.74	&9.78E-07	&6.58	&2.26E-04	&5.29&67
\\\cline{2-9}  
& 140 & 2.65E-05	&8.19	&2.53E-07	&8.77	&7.18E-05	&7.43&72
\\\cline{2-9}  
& 160 & 8.64E-06	&8.38	&7.44E-08	&9.17	&2.47E-05	&7.98&78
\\\cline{2-9}  
& 180 & 2.77E-06	&9.67	&2.19E-08	&10.4	&7.81E-06	&9.79&84
\\\cline{2-9}  
& 200 &  9.57E-07	&10.1	&6.96E-09	&10.9	&2.72E-06	&10.0&89 
\\\cline{2-9} 
& 220 &3.30E-07	&11.2	&2.25E-09	&11.8	&1.02E-06	&10.3&93  
\\\cline{2-9}
&240 &1.19E-07	&11.8	&7.60E-10	&12.5	&4.02E-07	&10.7&99 
\\\cline{2-9}  
& 260 &4.35E-08	&12.5	&2.65E-10	&13.2	&1.60E-07	&11.5&105
%\\\cline{2-9}  
%& 280 &1.67E-08	&12.9	&9.61E-11	&13.7	&6.42E-08	&12.3&111  
%\\\cline{2-9}  
%& 300 &6.49E-09	&13.7	&3.55E-11	&14.5	&2.63E-08	&12.9&115  
%\\\cline{2-9}  
%& 320 &2.64E-09	&13.9	&1.35E-11	&14.9	&1.09E-08	&13.7&123
\\\hline
\multicolumn{1}{|c|}{\multirow{4}{*}{WENO9  }}
%& 20 & -- & --  & -- & --  &  --  &-- & --  
%\\\cline{2-9}
%& 40 & -- & --  & -- & --  &  --  &-- & --  
%\\\cline{2-9}  
%& 60 &-- & --  & -- & --  &  --  &-- & --  
%\\\cline{2-9}  
%& 80 & -- & --  & -- & --  &  --  &-- & --  
%\\\cline{2-9}  
%& 100 &  -- & --  & -- & --  &  --  &-- & --  
%\\\cline{2-9}  
%& 120 & -- & --  & -- & --  &  --  &-- & --  
%\\\cline{2-9}  
& 140 & 9.83E-06	&-- &1.01E-07	&--	&2.99E-05	&-- &75  
\\\cline{2-9}  
& 160 & 3.34E-06	&8.08	&3.10E-08	&8.85	&1.06E-05	&7.78&82 
\\\cline{2-9}  
& 180 & 1.15E-06	&9.09	&9.87E-09	&9.72	&3.79E-06	&8.71&89
\\\cline{2-9}  
& 200 &  4.14E-07	&9.66	&3.29E-09	&10.4	&1.39E-06	&9.49&95
\\\cline{2-9}
& 220 & 1.45E-07	&11.0	&1.08E-09	&11.7	&5.32E-07	&10.1&102  
\\\cline{2-9}
&240 &5.24E-08	&11.7	&3.67E-10	&12.5	&2.07E-07	&10.8&109
\\\cline{2-9}  
& 260 &1.90E-08	&12.7	&1.27E-10	&13.3	&8.26E-08	&11.5&115 
%\\\cline{2-9}  
%& 280 &7.17E-09	&13.2	&4.52E-11	&13.9	&3.31E-08	&12.3&122  
%\\\cline{2-9}  
%& 300 &2.70E-09	&14.2	&1.63E-11	&14.8	&1.33E-08	&13.2&131
%\\\cline{2-9}  
%& 320 &1.05E-09	&14.7	&5.99E-12	&15.5	&7.42E-09	&9.00&135 
%\\\cline{2-9}  
%& 200 &  a & --  & a & --  &  a & a&
\\\hline
\end{tabular}
\label{2d_ex2_table1}
\end{center}
\end{table}

%++++++++++++++++++++++++++++++++++++++++++++++++++

\begin{exam}\label{2d_static_eikonal_4}
We consider the two-dimensional general static Eikonal equation (\ref{geikonal}) on the computational domain $\Omega = [-1, 1] \times [-1, 1]$ and the boundary condition is $T(0, 0)=1$ at the point source $(0, 0)$. 
$\mathbf{v(x)} = (v_1, v_2)^{\text{T}} = (0, 0)^\text{T}$. The speed function $F(\mathbf{x})$ is computed analytically as
\begin{equation}
F(x,y) =\frac{1}{\sqrt{4x^2(1+y^2)^2+4y^2(1+x^2)^2}}
\end{equation}
such that the exact solution is $T(x, y)=(1+x^2)(1+y^2)$.

\end{exam}
%++++++++++++++++++++++++++++++++++++++++++++++++++

We solve this problem by the proposed Newton-type WENO schemes using column-row-wise sweeping. $\epsilon_{\text{WENO}}=10^{-9}$ and $\alpha=2$ for all the WENO reconstructions. For WENO5 and WENO7, we take the convergence criterion $\epsilon = 10^{-12}$, while $\epsilon=10^{-9}$ for WENO9.  The numerical errors and orders are presented in Table \ref{2d_ex4_table1}, from which we can clearly observe the convergence and the high-order accuracy of the proposed schemes.

\begin{table}[H]
\begin{center}
\caption {Accuracy study for Example \ref{2d_static_eikonal_4}. Column-row-wise sweep.  }
\begin{tabular}{|c|c|c c|cc|cc|c|}
%\cline{3-8}
\hline
 & N  & $L^1$ error & order   & $L^2$ error  & order  & $L^{\infty}$  & order & sweep \\\hline
\multicolumn{1}{|c|}{\multirow{4}{*}{ WENO5 }}
& 20 &1.39E-03	&--	&1.83E-04	&--	&2.34E-02	&--&39 
\\\cline{2-9}
& 40 &  2.15E-05	&6.01	&1.71E-06	&6.74	&5.51E-04	&5.41&	90
\\\cline{2-9}  
& 60 &  4.63E-07	&9.47	&2.76E-08	&10.2	&2.09E-05	&8.07	&	116  
\\\cline{2-9}  
& 80 &   1.99E-08	&10.9	&1.06E-09	&11.3	&1.50E-06	&9.17		&150
\\\cline{2-9}  
& 100 &   1.49E-09	&11.6	&6.65E-11	&12.4	&6.19E-08	&14.3		&179  
\\\cline{2-9}  
& 120 &  1.35E-10	&13.2	&7.07E-12	&12.3	&1.87E-08	&6.58		&195  
\\\cline{2-9}  
& 140 &   9.23E-12	&17.4	&4.22E-13	&18.3	&1.69E-09	&15.6		&231  
\\\cline{2-9}  
& 160 &  6.71E-13	&19.6	&2.50E-14	&21.2	&7.31E-11	&23.5	&268  
\\\cline{2-9}  
& 180 &   8.99E-14	&17.1	&4.03E-15	&15.5	&2.03E-11	&10.9		&308  
\\\cline{2-9}  
& 200 &   1.11E-14	&19.8	&3.59E-16	&23.0	&2.75E-12	&18.9		&350  
\\\hline
\multicolumn{1}{|c|}{\multirow{4}{*}{WENO7  }}
%& 20 &  & --  &  & --  &   &.    &   
%\\\cline{2-9}
& 40 & 1.23E-05	&--	&8.27E-07	&--	&2.95E-04	&--	&227  
\\\cline{2-9}  
& 60 &  6.63E-07	&7.20	&3.25E-08	&7.99	&1.97E-05	&6.67		&294 
\\\cline{2-9}  
& 80 &   1.12E-07	&6.18	&5.55E-09	&6.14	&6.23E-06	&4.01		&338
\\\cline{2-9}  
& 100 &   1.32E-08	&9.59	&6.20E-10	&9.82	&1.24E-06	&7.21		&381 
\\\cline{2-9}  
& 120 &   1.58E-09	&11.6	&6.14E-11	&12.7	&1.56E-07	&11.4		&402  
\\\cline{2-9}  
& 140 &  2.28E-10	&12.6	&7.27E-12	&13.8	&1.84E-08	&13.8		&449  
\\\cline{2-9}  
& 160 & 4.42E-11	&12.3	&1.67E-12	&11.0	&6.17E-09	&8.19		&429  
\\\cline{2-9}  
& 180 &   8.05E-12	&14.4	&2.95E-13	&14.7	&1.64E-09	&11.3		&470  
\\\cline{2-9}  
& 200 & 1.17E-12	&18.3	&3.59E-14	&20.0	&2.10E-10	&19.5		&564 
\\\hline
\multicolumn{1}{|c|}{\multirow{4}{*}{WENO9  }}
%& 20 &  & --  &  & --  &   &.    &   
%\\\cline{2-9}
& 40 & 9.70E-06	&--	&6.68E-07	&--	&2.55E-04	&--	&308   
\\\cline{2-9}  
& 60 & 5.76E-07	&6.96	&2.65E-08	&7.96	&1.33E-05	&7.30	&415  
\\\cline{2-9}  
& 80 &   8.05E-08	&6.84	&3.99E-09	&6.59	&5.23E-06	&3.23		&498  
\\\cline{2-9}  
& 100 &  8.10E-09	&10.3	&2.89E-10	&11.8	&3.86E-07	&11.7		&635  
\\\cline{2-9}  
& 120 &  1.62E-09	&8.81	&6.53E-11	&8.16	&1.54E-07	&5.06		&698   
\\\cline{2-9}  
& 140 &   2.32E-10	&12.6	&7.95E-12	&13.7	&2.46E-08	&11.9		&832   
\\\cline{2-9}  
& 160 &  4.21E-11	&12.8	&1.31E-12	&13.5	&2.92E-09	&16.0		&919  
\\\cline{2-9}  
& 180 &   9.00E-12	&13.1	&3.11E-13	&12.2	&1.60E-09	&5.08		&1171  
\\\cline{2-9}  
& 200 & 1.34E-12	&18.1	&3.60E-14	&20.5	&1.79E-10	&20.8	&1280 
\\\hline
\end{tabular}
\label{2d_ex4_table1}
\end{center}
\end{table}

\subsection{Anisotropic Eikonal equation}
In this subsection, we test anisotropic Eikonal problems with nonzero background flow and assess the robustness of the implicit Newton-type WENO sweeping methods under directional bias. 
%-----------------------------------------------------------------

\newpage
\begin{exam}\label{2d_static_eikonal_3}
In this example, we investigate the anisotropic scenarios of the Eikonal equation (\ref{eq1}). On the computational domain $\Omega = [0, 1] \times [0, 1]$, the velocity of the moving fluid is set to be $\mathbf{v(x)} = (v_1, v_2)^{\text{T}} = (0.1, 0.2)^\text{T}$. The speed function $F(\mathbf{x})$ is computed analytically as
\begin{equation}
F(x,y) = \frac{|1-v_1 \cos(x)\sin(y) -v_2 \sin(x) \cos(y)|}{\sqrt{\cos^2(x)\sin^2(y) +\sin^2(x) \cos^2(y)}}
\end{equation}
such that the exact solution is $T(x, y)=\sin(x) \sin(y)$. The exact boundary condition $T=0$ is imposed on two sides of the domain along $x=0$ and $y=0$. 
\end{exam}
We simulate the problem by the proposed implicit Newton-type WENO schemes based on column-row-wise sweeping. The convergence criterion is set to $\epsilon=10^{-12}$. The numerical errors and orders are given in Table \ref{2d_ex3_table1_case1}, in which we can observe that the numerical schemes achieve high-order accuracy and the schemes are all convergent. 

\begin{table}[H]
\begin{center}
\caption {Numerical errors and orders for Example \ref{2d_static_eikonal_3}. Column-row-wise sweep. $\alpha = 1.01$.}
\begin{tabular}{|c|c|c c|cc|cc|c|}
%\cline{3-8}
\hline
 & N  & $L^1$ error & order   & $L^2$ error  & order  & $L^{\infty}$  & order & sweep \\\hline
\multicolumn{1}{|c|}{\multirow{4}{*}{ WENO5 }}
%& 20 &1.19E-05	&--&	2.16E-06&	--	&2.21E-04	&--	&	21
%\\\cline{2-9}
%& 40 & 2.42E-07&	5.61	&3.06E-08	&6.14	&1.01E-05	&4.45	&	25
%\\\cline{2-9}  
%& 60 & 6.62E-09	&8.88&	5.36E-10	&9.98	&4.22E-07&	7.83	&	31
& 60 & 6.62E-09	&--&	5.36E-10	&--	&4.22E-07&	--	&	31
\\\cline{2-9}  
& 80 &  3.97E-10	&9.78&	2.92E-11&	10.1&	3.91E-08&	8.27	&	39
\\\cline{2-9}  
& 100 &  3.31E-11	&11.1	&2.30E-12	&11.4	&5.17E-09	&9.06		&46
\\\cline{2-9}  
& 120 &  3.47E-12	&12.4	&2.35E-13	&12.5&	8.21E-10&	10.1	&	53
\\\cline{2-9}  
& 140 & 4.92E-13&	12.7	&2.80E-14&	13.8	&1.35E-10	&11.7	&	60
\\\cline{2-9}  
& 160 &  1.12E-13&	11.1	&3.78E-15	&15.0&	2.33E-11	&13.2	&	67
%\\\cline{2-9}  
%& 180 &  4.28E-14	&8.16	&6.63E-16	&14.8	&4.13E-12&	14.7	&	73
%\\\cline{2-9}  
%& 200 &  2.23E-14&	6.19	&2.26E-16	&10.2&	7.46E-13&	16.2		&80
\\\hline
\multicolumn{1}{|c|}{\multirow{4}{*}{WENO7  }}
%& 20 & 1.03E-05&	-- &	1.79E-06	&--&	1.72E-04	&--&		26
%\\\cline{2-9}
%& 40 & 2.24E-07	&5.53	&2.80E-08	&6.00	&1.00E-05	&4.10	&	32
%\\\cline{2-9}  
%& 60 & 1.02E-08&	7.62&	8.89E-10&	8.51&	5.79E-07&	7.03	&	37
& 60 & 1.02E-08&	--&	8.89E-10&	--&	5.79E-07&	--	&	37
\\\cline{2-9}  
& 80 &  8.22E-10&	8.75&	6.01E-11&	9.37&	9.02E-08	&6.46	&	42
\\\cline{2-9}  
& 100 &  8.72E-11	&10.1&	6.11E-12&	10.2&	1.72E-08	&7.43	&	50
\\\cline{2-9}  
& 120 &  1.09E-11	&11.4	&7.64E-13	&11.4&	3.48E-09&	8.75	&	58
\\\cline{2-9}  
& 140 &  1.55E-12&	12.7&	1.09E-13&	12.6	&7.34E-10&	10.1	&	65
\\\cline{2-9}  
& 160 &  2.41E-13	&13.9	&1.72E-14&	13.9	&1.59E-10&	11.4		&73
%\\\cline{2-9}  
%& 180 & 4.20E-14	&14.8&	2.90E-15	&15.1	&3.53E-11	 &12.8	&	80
%\\\cline{2-9}  
%& 200 & 8.57E-15&	15.1	&5.14E-16	&16.4&	7.96E-12&	14.1	&	88
\\\hline
\multicolumn{1}{|c|}{\multirow{4}{*}{WENO9  }}
& 60 &1.14E-08	&--&	9.17E-10	&--&	6.34E-07&--&		51
\\\cline{2-9}  
& 80 &  1.40E-09	&7.30&	1.02E-10	&7.64&	1.45E-07&	5.12	&	53
\\\cline{2-9}  
& 100 & 2.20E-10&	8.29&	1.55E-11	&8.44	&3.80E-08	&6.01	&	57
\\\cline{2-9}  
& 120 &  4.13E-11&	9.17&	2.88E-12	&9.23&	1.08E-08	&6.90	&	63
\\\cline{2-9}  
& 140 &  8.70E-12&	10.1&	6.07E-13&	10.1	&3.24E-09&	7.82	&	70
\\\cline{2-9}  
& 160 &  2.00E-12&	11.0&	1.40E-13&	11.0&	1.00E-09&	8.76	&	78
%\\\cline{2-9}  
%& 180 &  4.91E-13&	11.9	&3.43E-14&	11.9	&3.21E-10&	9.70&	86
%\\\cline{2-9}  
%& 200 &  1.27E-13	&12.8	&8.85E-15&	12.8	&1.04E-10	&10.6	&	94
\\\hline
\end{tabular}
\label{2d_ex3_table1_case1}
\end{center}
\end{table}

%%%%%%%%%%%%%%%%%%%%%%%%%%%%%%%%%%%%%%%%%%%%%%%%%%%%
\newpage
\begin{exam}\label{2d_static_eikonal_3_case3}
We consider the two-dimensional general static Eikonal equation 
\begin{equation}
|| \nabla \text{T}(\mathbf{x}) || -\frac{|1-\mathbf{v(x)} \cdot \nabla \text{T}(\mathbf{x}) |}{F(\mathbf{x})} = 0.
\end{equation}
On the computational domain $\Omega = [0, 1] \times [0, 1]$, the velocity of the moving fluid is set to be $\mathbf{v(x)} = (v_1, v_2)^{\text{T}} =(1, 0)^{\text{T}}$. The speed function $F(\mathbf{x})$ is computed analytically as
\begin{equation}\label{FFF}
F(x,y) = \frac{|1-v_1 \cos(x)\sin(y) -v_2 \sin(x) \cos(y)|}{\sqrt{\cos^2(x)\sin^2(y) +\sin^2(x) \cos^2(y)}}
\end{equation}
such that the exact solution is $T(x, y)=\sin(x) \sin(y)+1$. The exact boundary condition $T=1$ is imposed on two sides of the domain along $x=0$ and $y=0$. 
The numerical parameters are chosen as $\epsilon_{\text{iteration}} = 10^{-9}$, $\epsilon_{\text{WENO}} = 10^{-9}$ and $\alpha=1.1$. We first simulate the problem using column-row-wise sweep, and the numerical errors and orders are given in Table \ref{2dex3_table1_case3_0}, from which we can clearly observe that the designed high-order accuracies are all achieved.

We also notice that in this example $v_1=1$ and $v_2=0$, meaning that the background velocity field is entirely aligned with the $x$-direction. Motivated by this directional bias, we further simulate the problem using column-wise and row-wise sweeping only, and the numerical errors and orders are presented in Table \ref{2dex3_table2_case3_1}-\ref{2d_ex3_table3_case3_2}, where high-order accuracy is again achieved with comparable iteration counts.

\end{exam}

%%%%%%%%%

 \begin{table}[H]
 \begin{center}
 \caption {Accuracy study for Example \ref{2d_static_eikonal_3_case3}. Column-row-wise sweep.  }
 \begin{tabular}{|c|c|c c|cc|cc|c|}
 %\cline{3-8}
 \hline
  & N  & $L^1$ error & order   & $L^2$ error  & order  & $L^{\infty}$  & order & sweep \\\hline
 \multicolumn{1}{|c|}{\multirow{4}{*}{ WENO5 }}
 &10		&1.11E-04	&--&	3.88E-05&--&		2.16E-03&--&		17 \\\cline{2-9}
 &20		&1.26E-06&	6.46&	2.49E-07	&7.28&	5.46E-05&	5.30&	24\\\cline{2-9}
 &40		&6.29E-09&	7.64&	7.95E-10	&8.29&	3.04E-07&	7.49&	38 \\\cline{2-9}
 &60		&1.05E-10&	10.1&	1.09E-11	&10.6&	8.79E-09&	8.74&	46 \\\cline{2-9}
 &80		&6.53E-12&	9.65&	3.01E-13	&12.5&	2.89E-10&	11.9&	54 
 %\\\cline{2-9}
% &100	&3.48E-12&	2.82&	1.38E-13	&3.49&	2.01E-10&	1.63&	62
 \\\hline
 \multicolumn{1}{|c|}{\multirow{1}{*}{WENO7  }}
 %&  \multicolumn{8}{|c|} {No convergence}
 &10		&9.05E-05&--		&3.14E-05&--		&1.77E-03	&--&	23 \\\cline{2-9}
 &20		&9.43E-07&	6.58	&2.15E-07	&7.20	&5.75E-05	&4.94	&33\\\cline{2-9}
 &40		&5.24E-09&	7.49&	6.70E-10	&8.32	&2.79E-07	&7.69	&42 \\\cline{2-9}
 &60		&9.29E-11	&      9.94&	9.37E-12	&10.5	&9.09E-09	&8.44	&52\\\cline{2-9}
 &80		&1.54E-11& 	6.24&	6.34E-13	&9.36	&5.15E-10	&9.98	&60 
 %\\\cline{2-9}
% &100		&7.90E-12	&3.01	&2.83E-13	&3.61	&3.11E-10	&2.26	&70
\\\hline
 \multicolumn{1}{|c|}{\multirow{4}{*}{WENO9  }}
 &10		&8.27E-05&	--	        &2.71E-05	&--	&                1.47E-03&--		&34 \\\cline{2-9}
 &20		&5.77E-07&	7.16	&1.38E-07&	7.61	&4.26E-05	&5.11&	42 \\\cline{2-9}
 &40		&3.71E-09&	7.28	&4.33E-10&	8.32	&1.79E-07	&7.90&	48\\\cline{2-9}
 &60		&8.40E-11& 	9.35	&7.93E-12&	9.86	&6.28E-09	&8.26&	60 \\\cline{2-9}
 &80		&8.40E-12&	8.00	&4.12E-13&	10.3	&4.73E-10	&8.99&	73
 % \\\cline{2-9}
% &100	&4.06E-12&	3.25	&1.99E-13&	3.25	&4.43E-10	&0.23&	85
 \\\hline
 \end{tabular}
 \label{2dex3_table1_case3_0}
 \end{center}
 \end{table}

\begin{table}[H]
\begin{center}
\caption {Accuracy study for Example \ref{2d_static_eikonal_3_case3}. Column-wise sweep. }
\begin{tabular}{|c|c|c c|cc|cc|c|}
%\cline{3-8}
\hline
 & N  & $L^1$ error & order   & $L^2$ error  & order  & $L^{\infty}$  & order & sweep \\\hline
\multicolumn{1}{|c|}{\multirow{4}{*}{ WENO5 }}
&10		&1.10E-04 & --		&3.88E-05 & --		&2.16E-03 &--		&29 \\\cline{2-9}  
&20		& 1.26E-06&	6.46	&2.49E-07	&7.28	&5.46E-05	&5.30&	42 \\\cline{2-9}  
&40		&6.30E-09 &	7.64&7.95E-10	&8.29	&3.04E-07	&7.49		&61 \\\cline{2-9}  
&60		&1.08E-10&	10.0	&1.09E-11 &	10.6 &	8.79E-09 &	8.74&	71\\\cline{2-9}  
&80		&1.24E-11	&7.54	&6.31E-13	&9.90	&8.22E-10&	8.24	&73 
%\\\cline{2-9}    
%&100		&9.54E-12	&1.16	&3.47E-13	&2.67	&5.03E-10	&2.21	&78
\\\hline
\multicolumn{1}{|c|}{\multirow{4}{*}{WENO7  }}
&10		&9.05E-05&--		&3.14E-05&--		&1.77E-03&--		&43\\\cline{2-9}
 &20		&9.43E-07&	6.58&	2.15E-07	&7.20&	5.75E-05	&4.94	&54 \\\cline{2-9}
 &40		&5.24E-09&	7.49&	6.70E-10	&8.32&	2.79E-07	&7.69	&76 \\\cline{2-9}
 &60		&1.01E-10&	9.75&	9.50E-12&	10.5&	9.09E-09&	8.44&	89\\\cline{2-9}
 &80		&1.75E-11	&  6.08&	1.08E-12&	7.57&	1.33E-09	&6.69	&95 
 %\\\cline{2-9}
 %&100		&1.19E-11	&1.73	&5.42E-13	&3.08	&9.96E-10	&1.29	&99
\\\hline
\multicolumn{1}{|c|}{\multirow{4}{*}{WENO9  }}
&10		&8.28E-05&--		&2.71E-05&--&		1.48E-03&--	&	69 \\\cline{2-9}
 &20		&5.77E-07&	7.17	&1.38E-07	&7.62	&4.27E-05	&5.11	&75 \\\cline{2-9}
 &40		&3.72E-09&	7.28	&4.32E-10	&8.32	&1.79E-07	&7.90	&90 \\\cline{2-9}
 &60		&9.23E-11	&  9.11	&8.05E-12	&9.82	&6.33E-09	&8.24	&112 \\\cline{2-9}
 &80		&2.30E-12&	12.8	&1.96E-13	&12.9	&2.38E-10	&11.4	&284
 %\\\cline{2-9}
 %&100		&1.44E-11&	-8.23&	9.81E-13&	-7.21&	1.84E-09	&-9.17	&5000
\\\hline
\end{tabular}
\label{2dex3_table2_case3_1}
\end{center}
\end{table}

 \begin{table}[H]
 \begin{center}
 \caption {Accuracy study for Example \ref{2d_static_eikonal_3_case3}. Row-wise sweep. }
 \begin{tabular}{|c|c|c c|cc|cc|c|}
 %\cline{3-8}
 \hline
  & N  & $L^1$ error & order   & $L^2$ error  & order  & $L^{\infty}$  & order & sweep \\\hline
 \multicolumn{1}{|c|}{\multirow{4}{*}{ WENO5 }}
 &10		&1.11E-04&--		&3.89E-05&--		&2.16E-03&--		&38 \\\cline{2-9}
 &20		&1.26E-06&	6.46&	2.49E-07&	7.28	&5.46E-05	&5.31&	66\\\cline{2-9}
 &40		&6.31E-09&	7.64&	7.95E-10&	8.29&	3.04E-07	&7.49&	116 \\\cline{2-9}
 &60		&1.09E-10&	10.0&	1.09E-11&	10.6&	8.77E-09	&8.74&	169\\\cline{2-9}
 &80		&5.46E-12&	10.4&	2.48E-13&	13.2&	2.83E-10	&11.9& 234
 %\\\cline{2-9}
 %&100		&1.28E-12&	6.49&	2.54E-14&	10.2&	1.57E-11&	13.0&	307
 \\\hline
 \multicolumn{1}{|c|}{\multirow{4}{*}{WENO7  }}
 &10		&9.06E-05&--		&3.15E-05&--		&1.77E-03&--		&46 \\\cline{2-9}
 &20		&9.43E-07&	6.59	&2.14E-07	&7.20&	5.75E-05	&4.95&	73 \\\cline{2-9}
 &40		&5.25E-09&	7.49	&6.70E-10	&8.32&	2.79E-07	&7.69&	123  \\\cline{2-9}
 &60		&9.14E-11	&      9.99	&9.37E-12	&10.5&	9.13E-09	&8.43&	157  \\\cline{2-9}
 &80		&3.14E-12&	11.7	&2.07E-13	&13.3&	2.89E-10	&12.0&	204 
 %\\\cline{2-9}
 %&100	&4.41E-13&	8.80	&1.04E-14	&13.4&	2.07E-11	&11.8&	253
  \\\hline
 \multicolumn{1}{|c|}{\multirow{4}{*}{WENO9  }}
 &10		&8.25E-05&--		&2.70E-05&--&	 1.46E-03&--&		71 \\\cline{2-9}
 &20		&5.77E-07&	7.16&	1.38E-07&	7.61&	4.27E-05	&5.10	&94 \\\cline{2-9}
 &40		&3.71E-09&	7.28	&4.33E-10&	8.32	&1.78E-07	&7.90	&158 \\\cline{2-9}
 &60		&9.26E-11	&9.11&	8.22E-12&	9.78	&6.54E-09&	8.15	&161 \\\cline{2-9}
 &80		&3.34E-12&	11.5&	2.18E-13&	12.6&	2.47E-10	&11.4	&207 
 %\\\cline{2-9}
% &100		&1.80E-13&	13.1&	9.97E-15&	13.8	&2.17E-11&	10.9	&257
 \\\hline
 \end{tabular}
 \label{2d_ex3_table3_case3_2}
 \end{center}
 \end{table}

%----------------------------------------------------

\newpage
% %++++++++++++++++++++++++++++++++++++++++++++++++++
\begin{exam}\label{3d_static_eikonal_7}
We consider the three-dimensional general static Eikonal equation (\ref{geikonal}) on the computational domain $\Omega = [0, 1] \times [0, 1] \times [0, 1]$. The background velocity is chosen as $\mathbf{v(x)} = (v_1, v_2, v_3)^{\text{T}} = (0.1, 0.2, 0.3)^\text{T}$. The boundary condition $T(x, y, z)=0$ is imposed along the three sides $x = 0$, $y=0$ and $z=0$. Similarly, the speed function $F(\mathbf{x})$ is computed analytically from 
\begin{equation}
F(x,y,z) = \frac{|1 - \mathbf{v} \cdot \nabla T(x, y, z)|}{||\nabla T(x, y, z)||}
\end{equation}
by assuming that the exact solution is $T(x, y, z)=\sin(x) \sin(y) \sin(z)$.

\end{exam}

\begin{figure}[p]
  %\centering
 % \vspace{-1.0in}
  \begin{minipage}[b]{0.23\textwidth}
    \centerline{
    \includegraphics[width=1.5in,angle=0,scale=1.05]{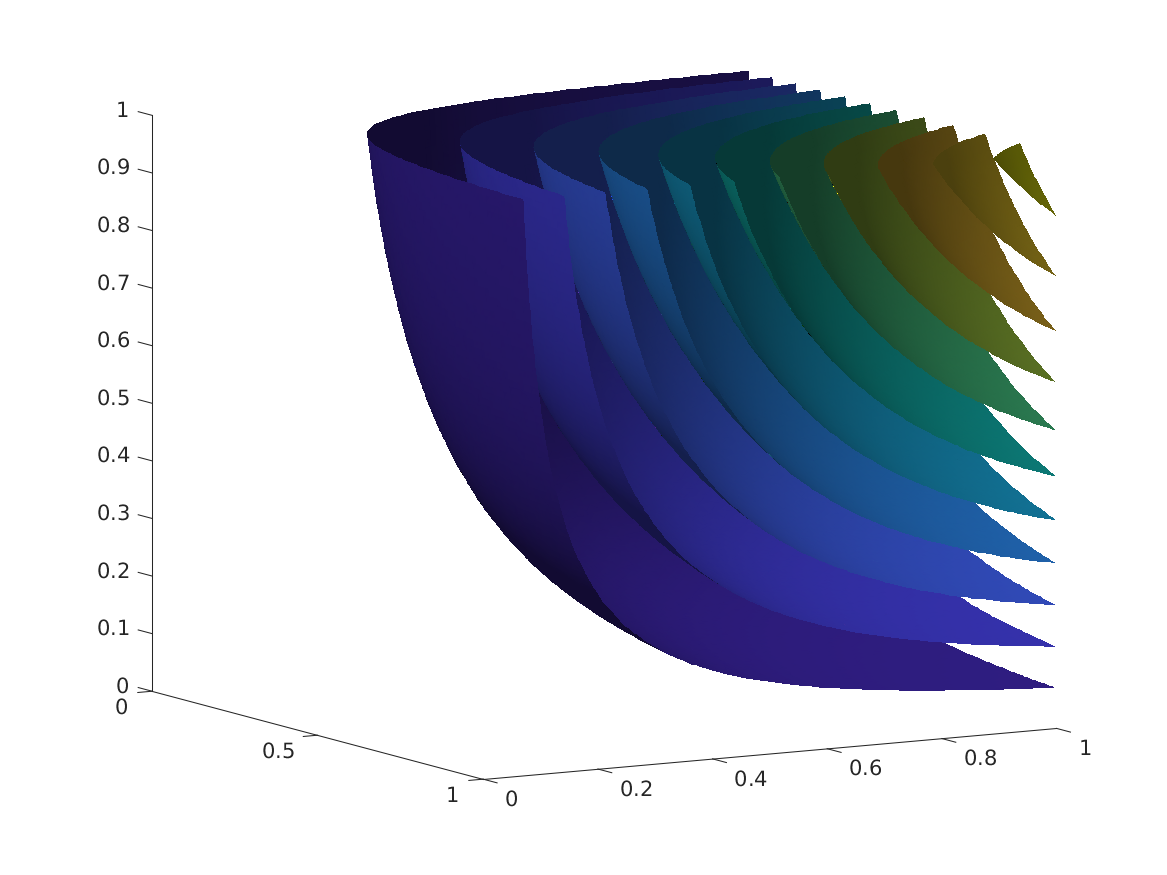}}
   % \caption{Caption 2}
 %  \vspace{-1.2in}
  \qquad     (a1) 
 %  \medskip
  \end{minipage}
%  \hspace{0.2\textwidth}
 %   \vspace{-1.0in}
  \begin{minipage}[b]{0.23\textwidth}
    \centerline{
    \includegraphics[width=1.5in,angle=0,scale=1.05]{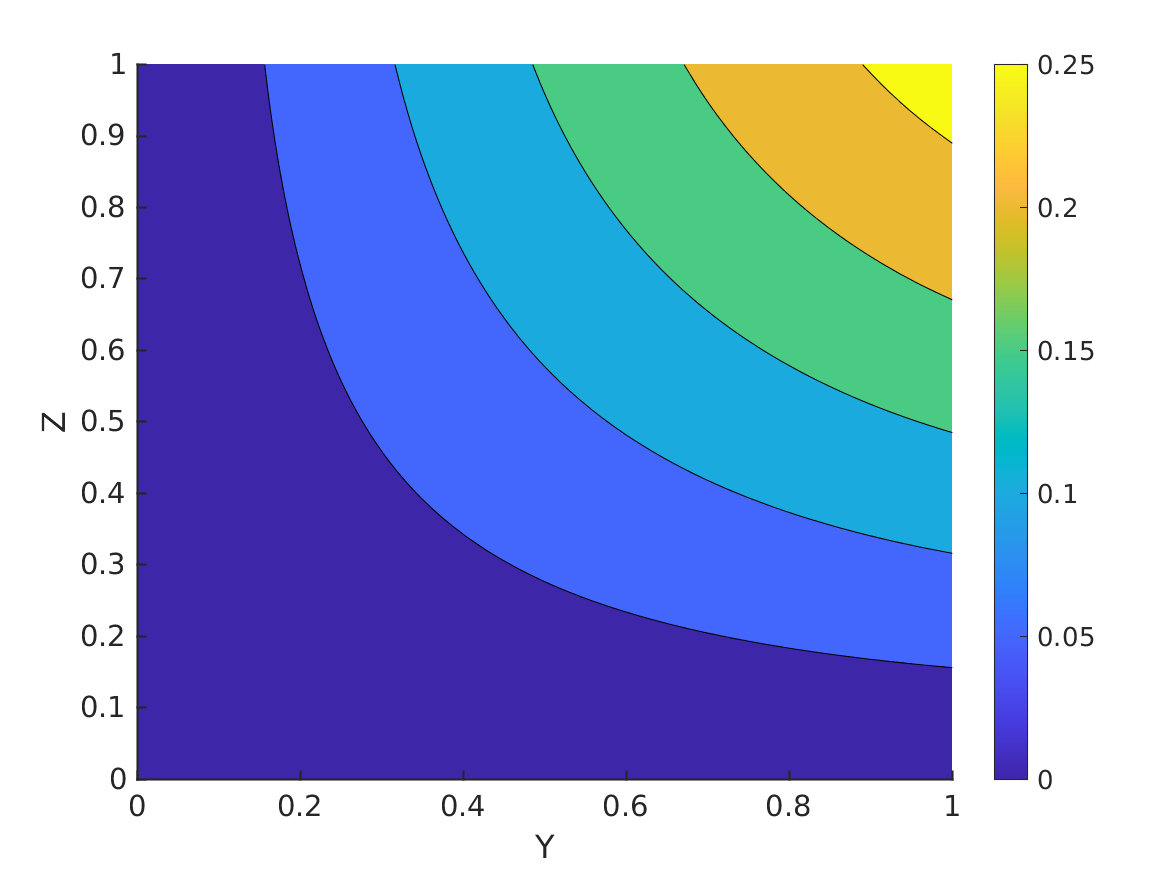}}
    %\caption{Caption 3}
   %  \vspace{-1.2in}
 \qquad   (a2) 
%    \medskip
  \end{minipage}
\begin{minipage}[b]{0.23\textwidth}
    \centerline{
    \includegraphics[width=1.5in,angle=0,scale=1.05]{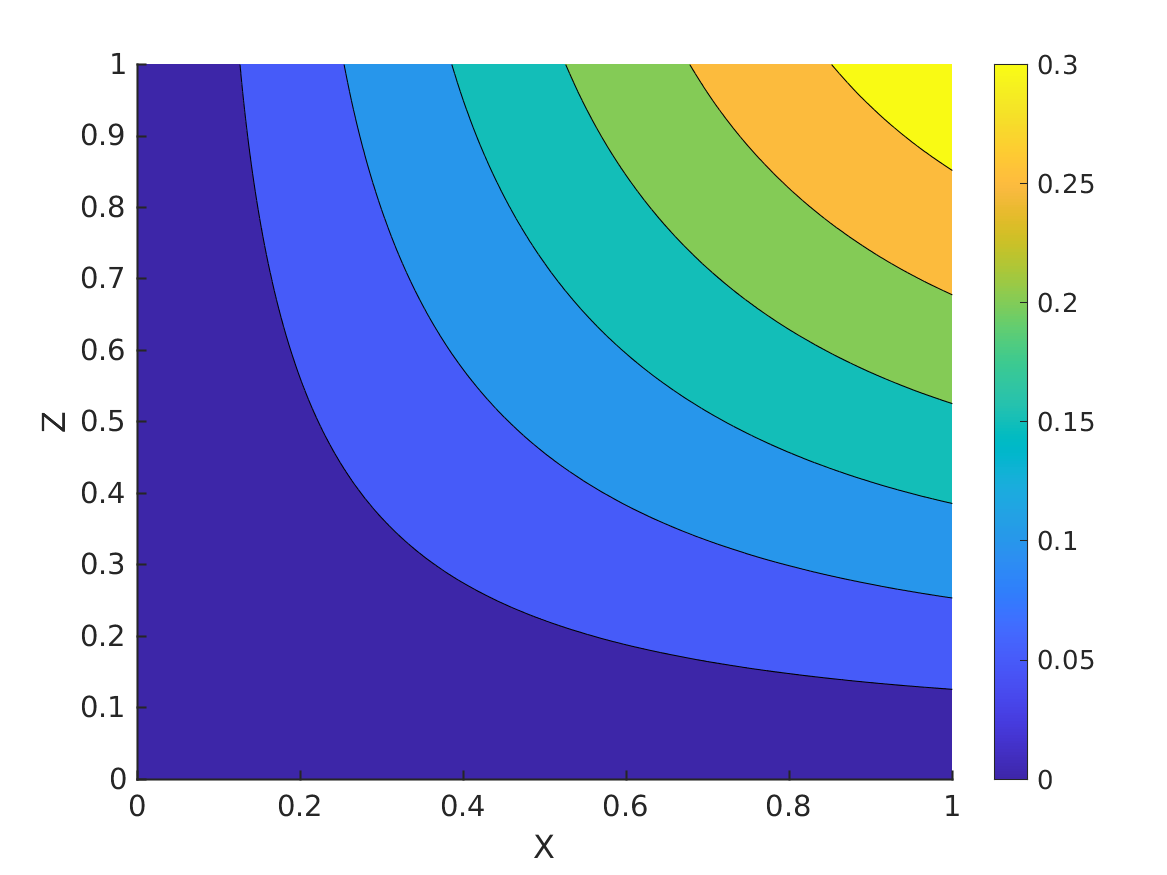}}
   % \caption{Caption 2}
 %  \vspace{-1.2in}
  \qquad      (a3) 
 %  \medskip
  \end{minipage}
 % \hspace{0.2\textwidth}
 %   \vspace{-1.0in}
  \begin{minipage}[b]{0.23\textwidth}
    \centerline{
    \includegraphics[width=1.5in,angle=0,scale=1.05]{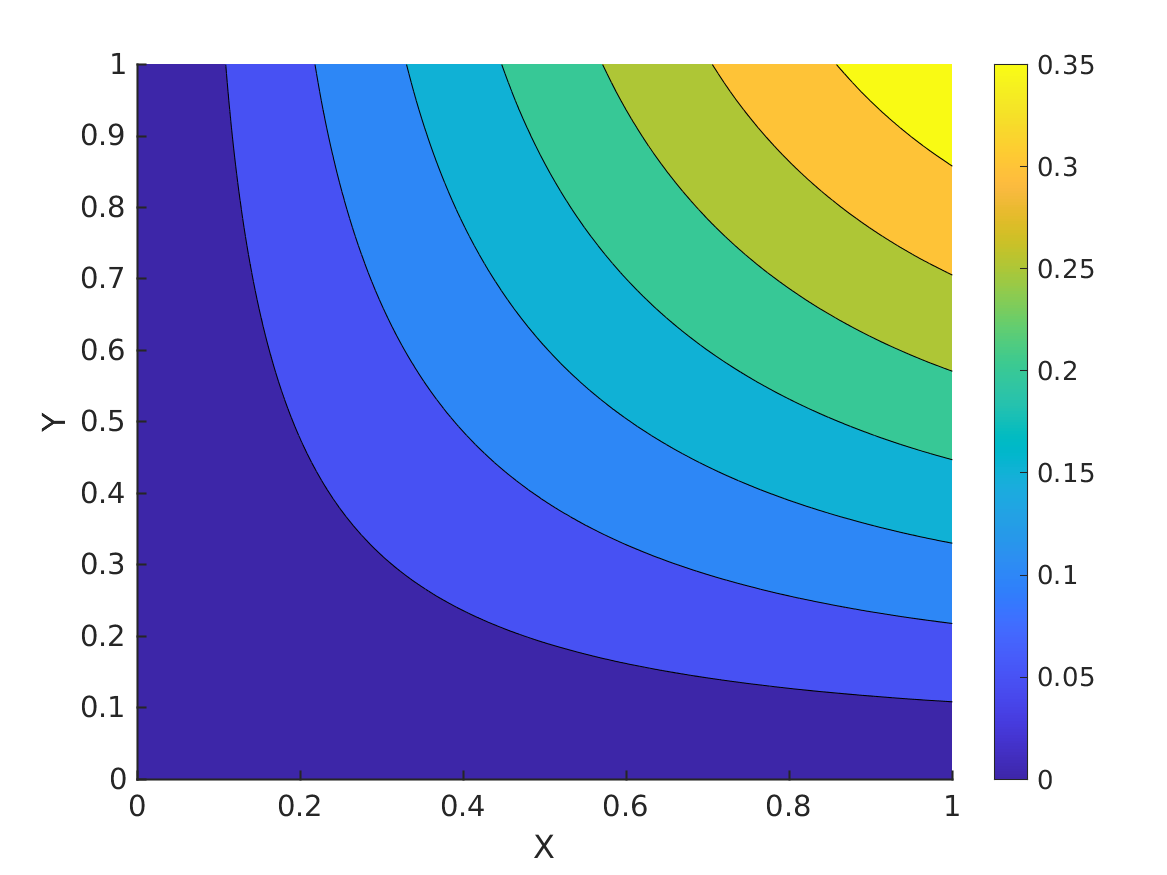}}
    %\caption{Caption 3}
   %  \vspace{-1.2in}
  \qquad      (a4)  
%    \medskip
  \end{minipage}  
 \begin{minipage}[b]{0.23\textwidth}
    \centerline{
    \includegraphics[width=1.5in,angle=0,scale=1.05]{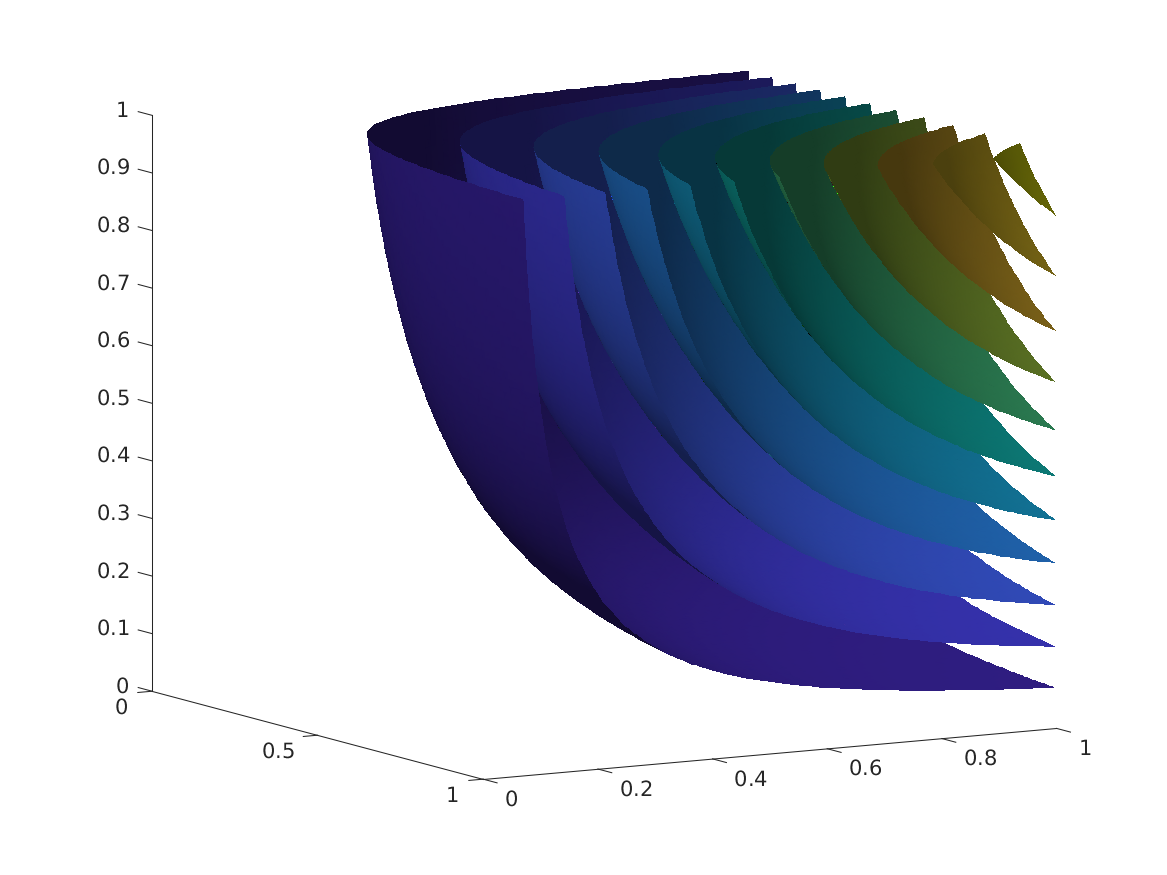}}
   % \caption{Caption 2}
 %  \vspace{-1.2in}
  \qquad     (b1) 
 %  \medskip
  \end{minipage}
%  \hspace{0.2\textwidth}
 %   \vspace{-1.0in}
  \begin{minipage}[b]{0.23\textwidth}
    \centerline{
    \includegraphics[width=1.5in,angle=0,scale=1.05]{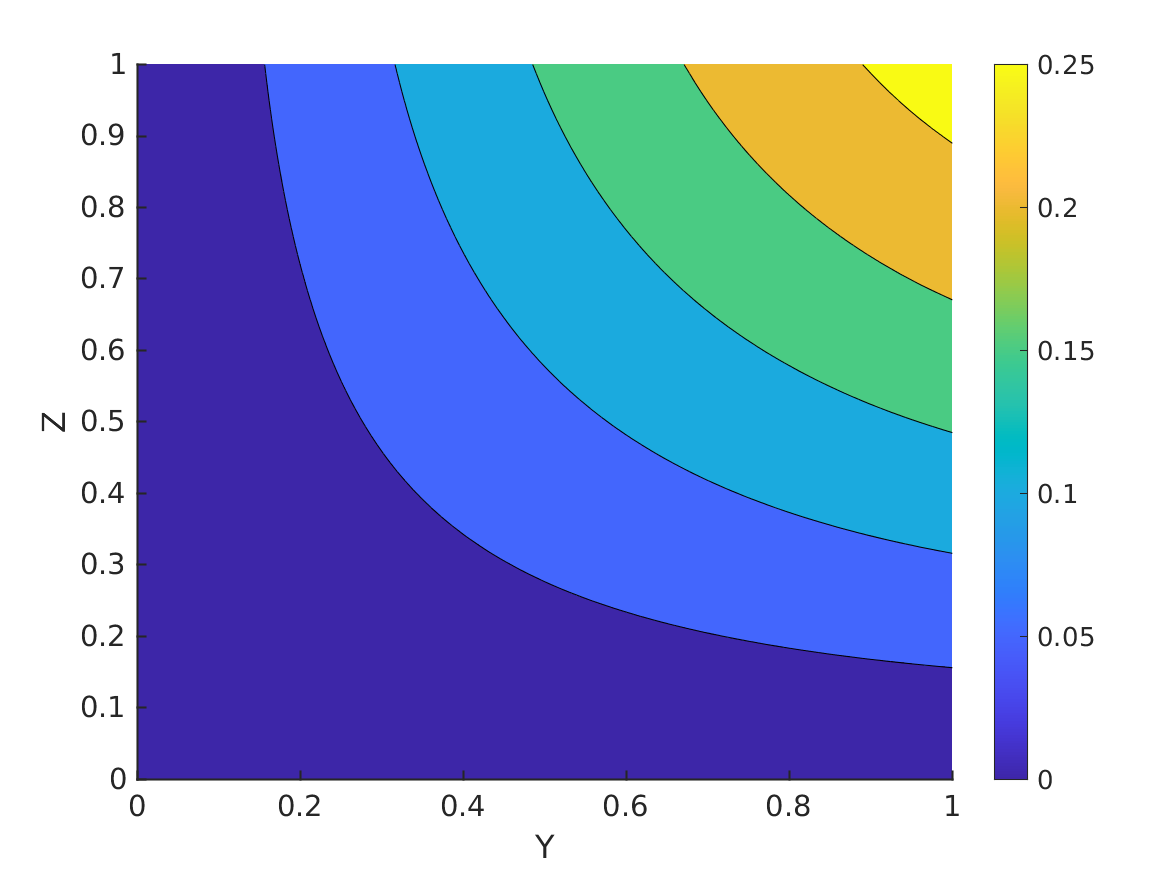}}
    %\caption{Caption 3}
   %  \vspace{-1.2in}
 \qquad   (b2) 
%    \medskip
  \end{minipage}
\begin{minipage}[b]{0.23\textwidth}
    \centerline{
    \includegraphics[width=1.5in,angle=0,scale=1.05]{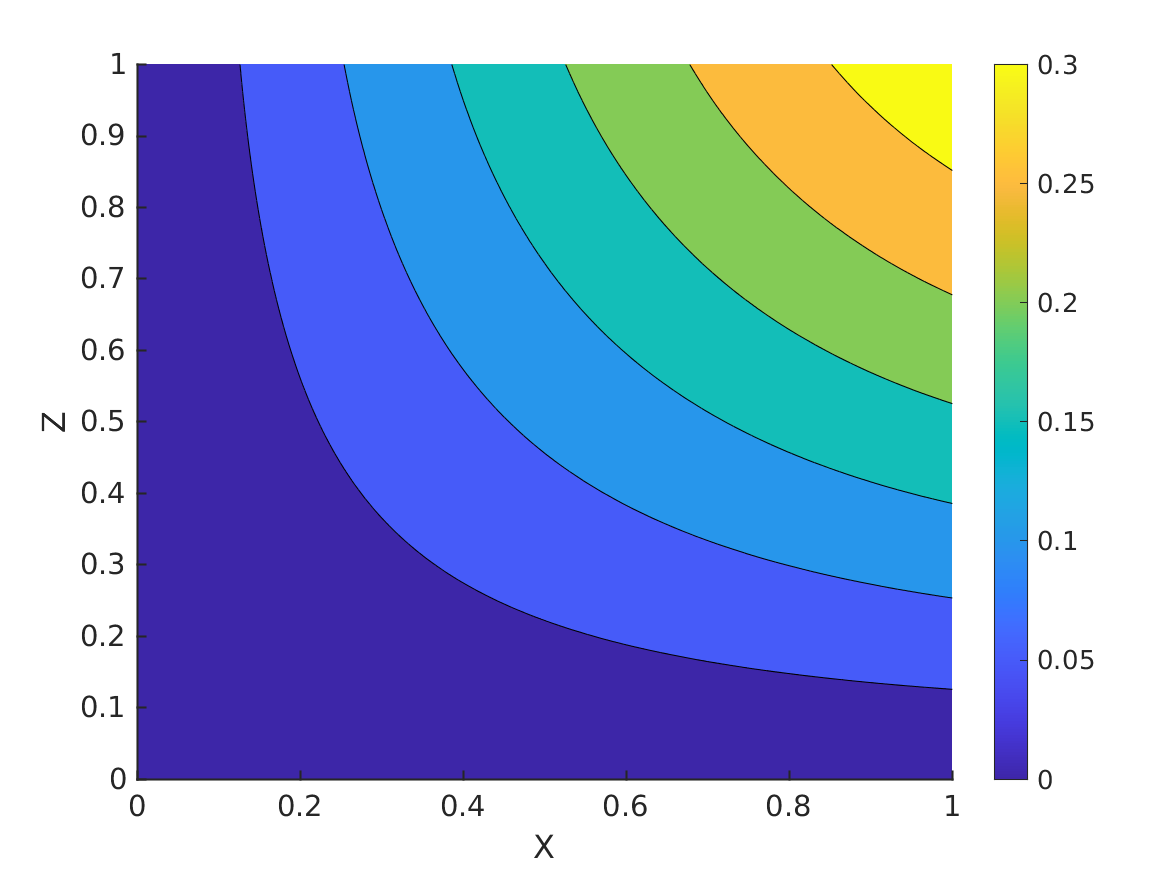}}
   % \caption{Caption 2}
 %  \vspace{-1.2in}
  \qquad      (b3) 
 %  \medskip
  \end{minipage}
 % \hspace{0.2\textwidth}
 %   \vspace{-1.0in}
  \begin{minipage}[b]{0.23\textwidth}
    \centerline{
    \includegraphics[width=1.5in,angle=0,scale=1.05]{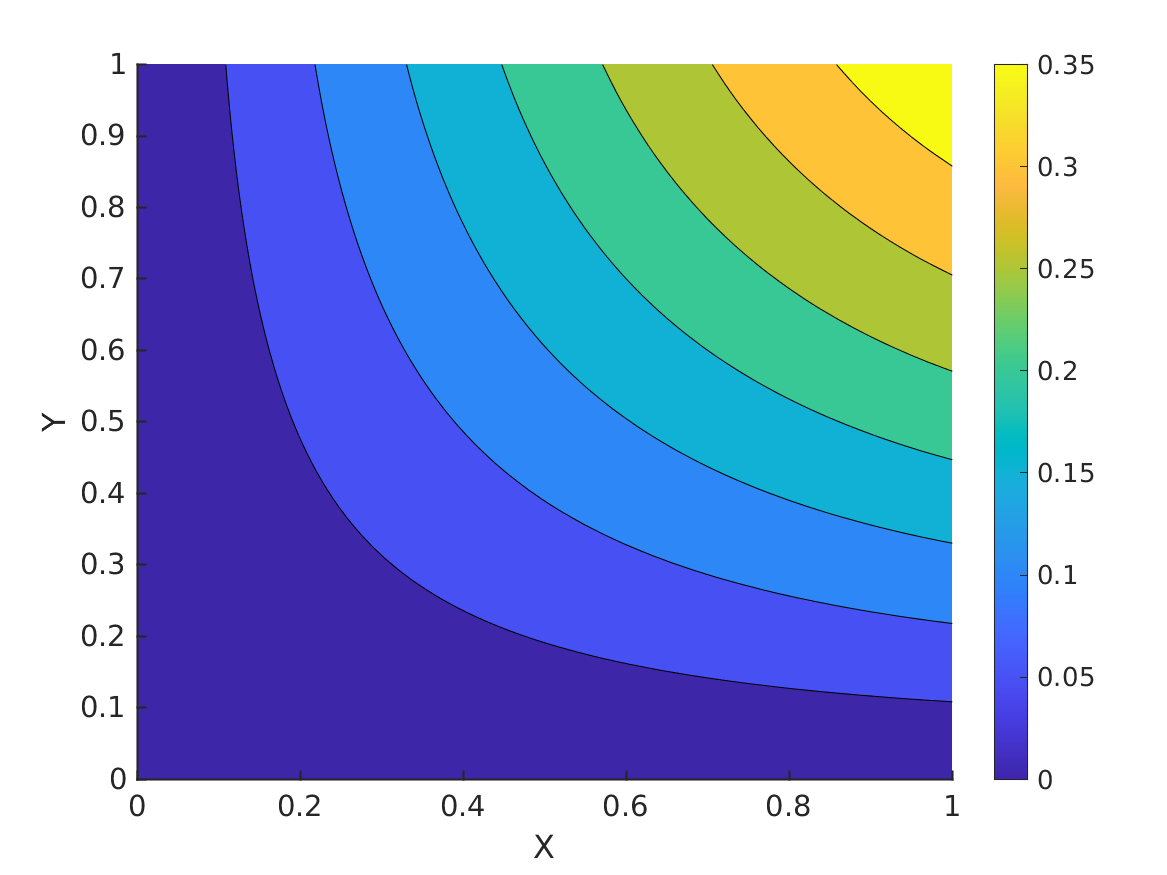}}
    %\caption{Caption 3}
   %  \vspace{-1.2in}
  \qquad      (b4)  
%    \medskip
  \end{minipage}  
   \begin{minipage}[b]{0.23\textwidth}
    \centerline{
    \includegraphics[width=1.5in,angle=0,scale=1.05]{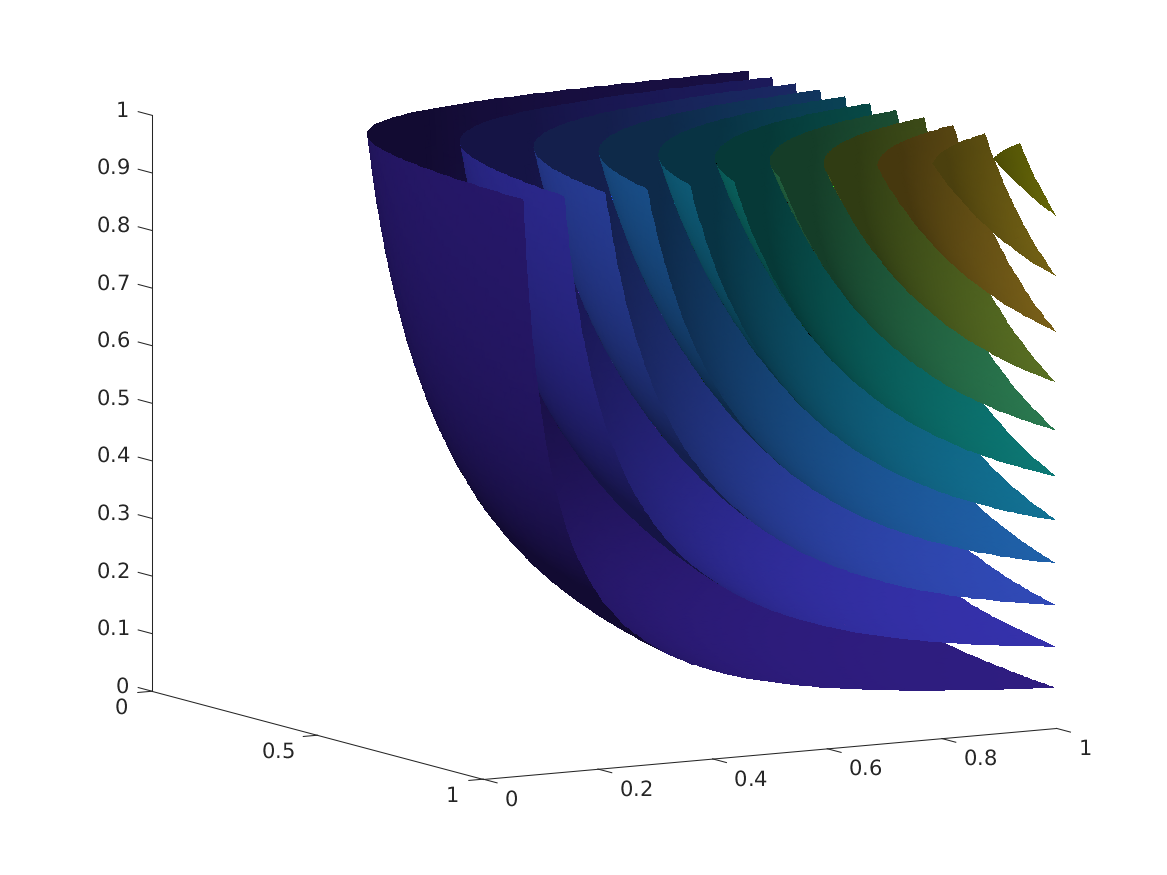}}
   % \caption{Caption 2}
 %  \vspace{-1.2in}
  \qquad     (c1) 
 %  \medskip
  \end{minipage}
%  \hspace{0.2\textwidth}
 %   \vspace{-1.0in}
  \begin{minipage}[b]{0.23\textwidth}
    \centerline{
    \includegraphics[width=1.5in,angle=0,scale=1.05]{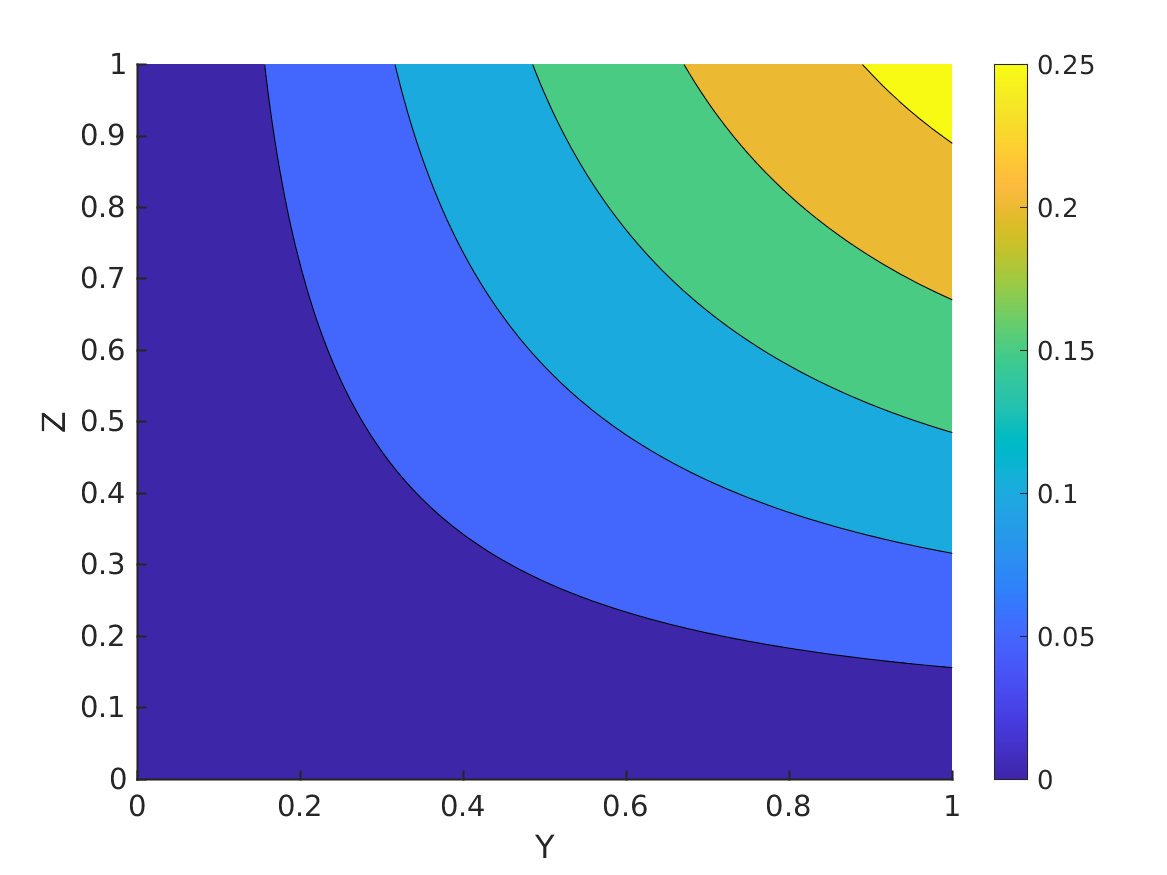}}
    %\caption{Caption 3}
   %  \vspace{-1.2in}
 \qquad   (c2) 
%    \medskip
  \end{minipage}
\begin{minipage}[b]{0.23\textwidth}
    \centerline{
    \includegraphics[width=1.5in,angle=0,scale=1.05]{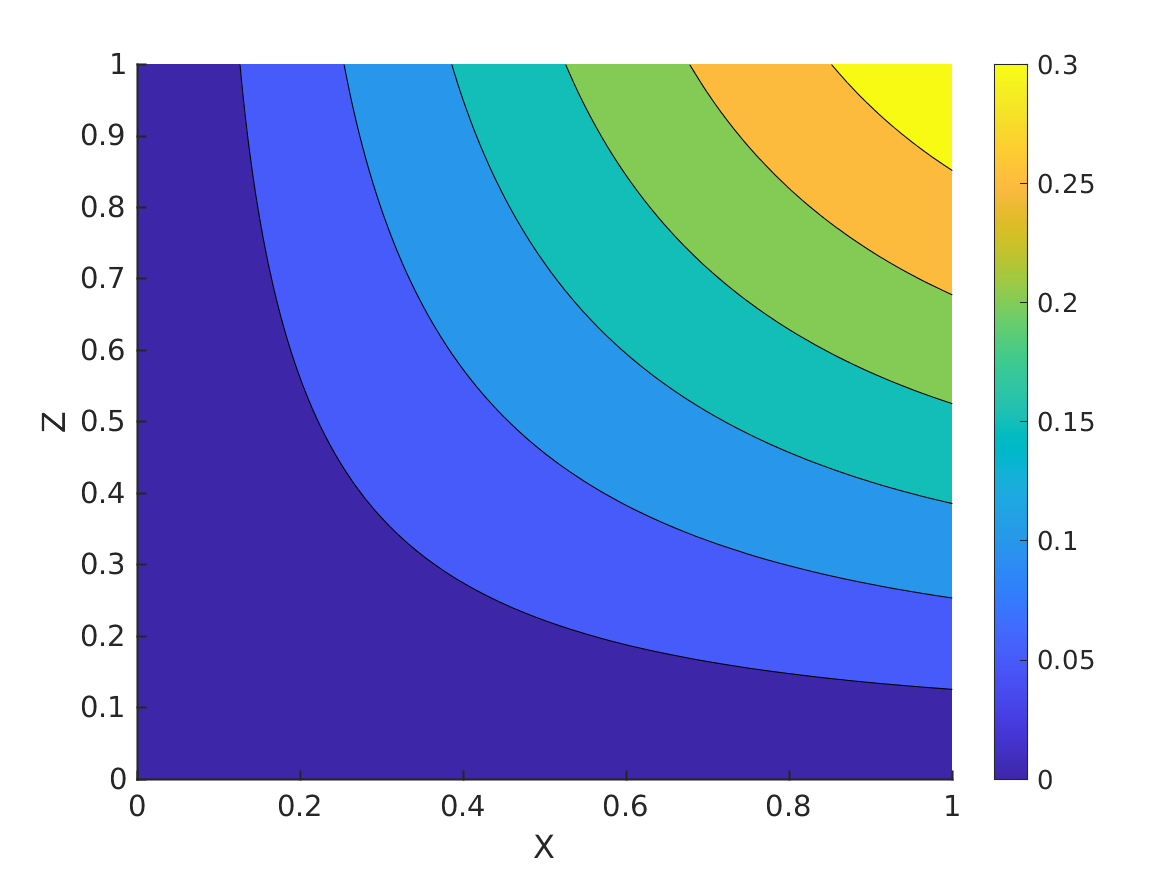}}
   % \caption{Caption 2}
 %  \vspace{-1.2in}
  \qquad      (c3) 
 %  \medskip
  \end{minipage}
 % \hspace{0.2\textwidth}
 %   \vspace{-1.0in}
  \begin{minipage}[b]{0.23\textwidth}
    \centerline{
    \includegraphics[width=1.5in,angle=0,scale=1.05]{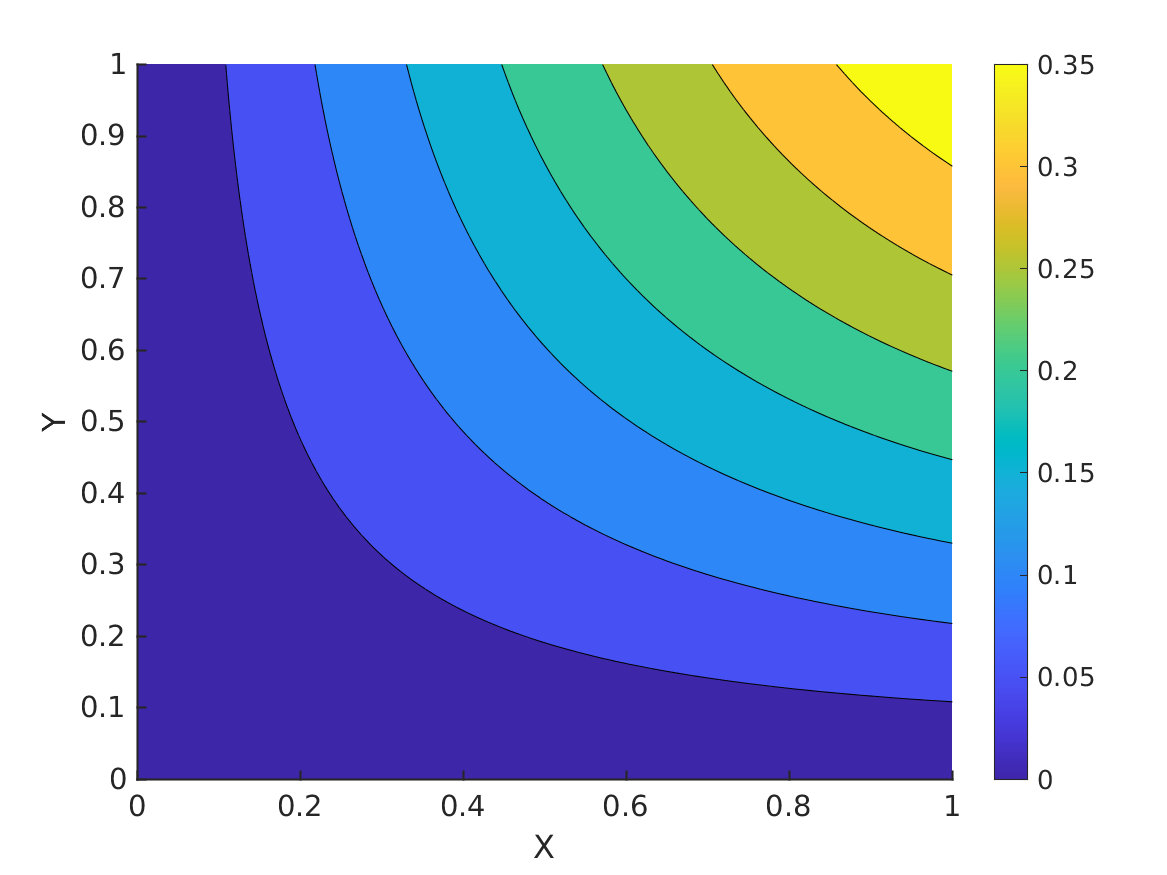}}
    %\caption{Caption 3}
   %  \vspace{-1.2in}
  \qquad      (c4)  
%    \medskip
  \end{minipage}  
    \caption{Numerical solutions of Example \ref{3d_static_eikonal_7}. $N_x =N_y=N_z = 80$. (a1)-(a4) Newton-type WENO5; (b1)-(b4) Newton-type WENO7; (c1)-(c4) Newton-type WENO9. (a1)-(c1) Isosurfaces of the solution from T = 0.05 to T = 0.55 with an increment of 0.05. (a2)-(c2) Cross-section of the solution along $x = 0.4$. (a3)-(c3) Cross-section along $y = 0.5$. (a4)-(c4) Cross-section along $z = 0.6$. } \label{figure_ex7_3d_accuracy}
\end{figure}

Let the partition be $x_i = (i-1)\Delta x$,  $y_j = (j-1)\Delta y$ and $z_k = (k-1)\Delta z$, we follow the approach in \cite{li2020newton} and simulate the problem using line-wise in $i$-th direction only. At each iteration, the numerical quantities \{$T(2, j, k), T(3, j, k), ..., T(Nx-1,j,k)$\} will be updated simultaneously, while the global solver sweeps the 3D computational domain in four directions,
\begin{eqnarray}
&& (1)\; k = 1: Nz, \quad j = 1: Ny; \qquad (2)\; k = 1:Nz, j = Ny:-1:1; \\ \nonumber
&& (3)\; k = Nz:-1:1, j=1:Ny; \quad (4)\; k = Nz:-1:1, j = Ny:-1:1.        
\end{eqnarray}
The convergence iteration criterion is $\epsilon = 10^{-9}$ and $\epsilon_{\text{WENO}} = 10^{-9}$ for all cases. There are differences in the numerical boundary condition setup. For WENO5, the boundary condition is enforced on the computational comain. For WENO7, we enforce exact solution boundary condition on the physical domain boundary and  4 neighboring points. For WENO9, exact solution boundary values are enforced with 1 neighboring point for $N=10$ and 5 neighboring points for other mesh partitions. $\alpha =2$ for WENO5 and WENO7, while $\alpha=1$ for WENO9. The numerical errors and orders are presented in Table \ref{3d_ex1_accuracy} in which we can clearly see the high-order accuracies are all achieved.  The isosurfaces plots of the travel time and the cross-sectional solutions given in Figure \ref{figure_ex7_3d_accuracy} demonstrate that the proposed high-order reconstruction could produce comparable numerical results on relatively coarse grids.

\begin{table} [htb]
\begin{center}
\caption {Accuracy study for Example \ref{3d_static_eikonal_7}.}
\begin{tabular}{|c|c|c c|cc|cc|c|}
%\cline{3-8}
\hline
 & N  & $L^1$ error & order   & $L^2$ error  & order  & $L^{\infty}$  & order & sweep \\\hline
\multicolumn{1}{|c|}{\multirow{4}{*}{ WENO5 }}
& 10 & 7.47E-05	&--	&1.00E-05	&-- &4.40E-03&-- &159
\\\cline{2-9}  
& 20 &6.56E-07	&6.83	&3.11E-08	&8.33 &4.68e-05	&6.56 &139
\\\cline{2-9}  
& 40 & 7.68E-09	&6.42	&1.48E-10	&7.71 &6.52E-07	&6.16 &72
\\\cline{2-9}  
%& 60 & --	&--	&--	&-- &--	&-- &--
%\\\cline{2-9}  
& 80 & 8.33E-12	&9.85	&5.06E-14	&11.5 &8.12E-10	&9.65 &186
\\\hline
\multicolumn{1}{|c|}{\multirow{4}{*}{WENO7}}
& 10 & 2.04E-05	&--	&5.72E-06	&-- &3.48E-03	&-- &25
\\\cline{2-9}  
& 20 & 3.86E-07	&5.72	&1.95E-08	&8.20 &3.83E-05	&6.51 &44
\\\cline{2-9}  
& 40 & 5.77E-09	&6.06	&1.09E-10	&7.48 &5.31E-07	&6.17 &71
\\\cline{2-9}  
%& 60 & --	&--	&--	&-- &--	&-- &--
%\\\cline{2-9}  
& 80 & 6.35E-12	&9.83	&5.69E-14	&10.9 &9.02E-10	&9.20 &129
\\\hline
\multicolumn{1}{|c|}{\multirow{4}{*}{WENO9  }}
& 10 & 2.4191e-05	&--	&2.9761e-06	&-- &9.5208e-04	&-- &65
\\\cline{2-9}  
& 20 & 2.1367e-07	&6.82	&1.8204e-08	&7.35 &4.7217e-05	&4.33 &47
\\\cline{2-9}  
& 40 & 1.2292e-09	&7.44	&3.4159e-11	&9.06 &6.3626e-07	&6.21 &75
\\\cline{2-9}  
%& 60 & --	&--	&--	&-- &--	&-- &--
%\\\cline{2-9}  
& 80 & 6.8720e-13	&10.8	&7.3425e-15	&12.2 &7.8340e-10	&9.67 &146
\\\hline
\end{tabular}
\label{3d_ex1_accuracy}
\end{center}
\end{table}

\begin{figure}[!htb]
  %\centering
 % \vspace{-1.0in}
  \begin{minipage}[b]{0.3\textwidth}
    \centerline{
    \includegraphics[width=2.1in,angle=0,scale=1.6]{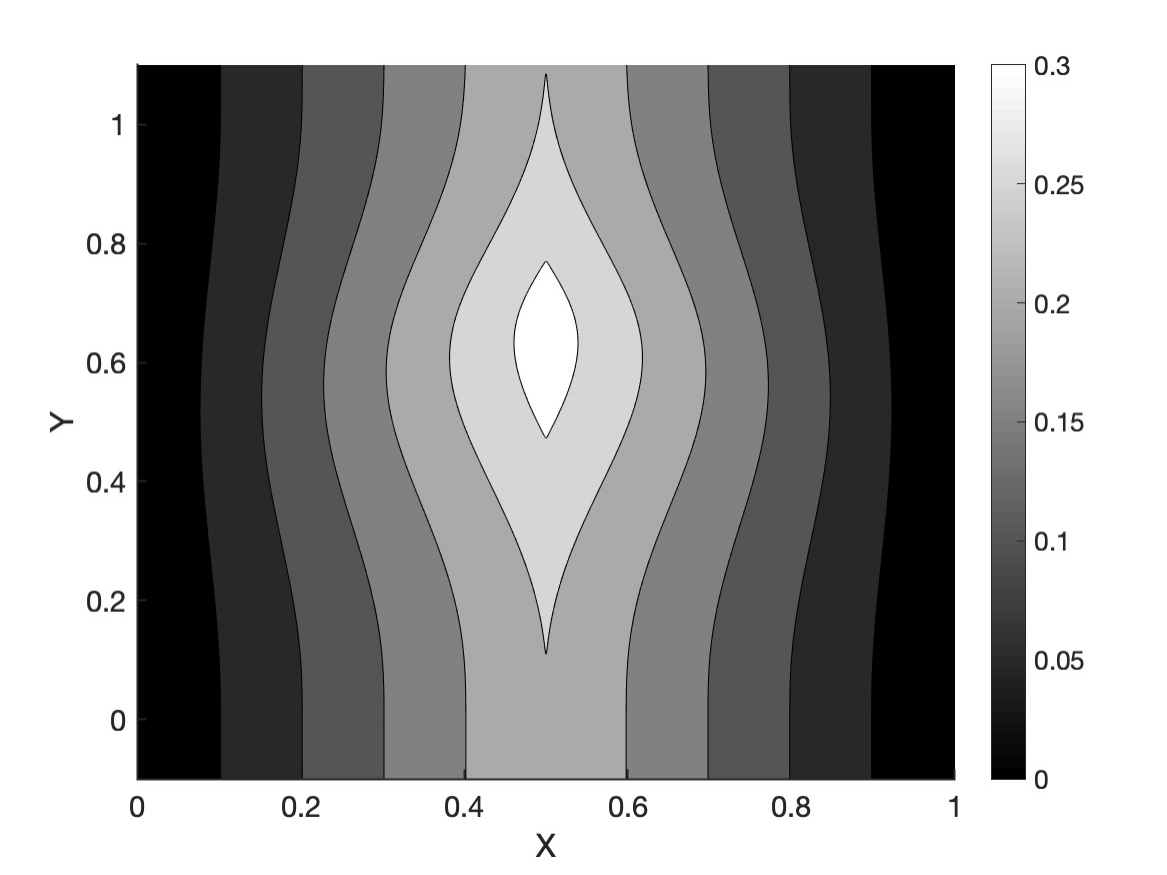}}
   % \caption{Caption 2}
 %  \vspace{-1.2in}
  \qquad     (a1) 
 %  \medskip
  \end{minipage}
  \hspace{0.2\textwidth}
 %   \vspace{-1.0in}
  \begin{minipage}[b]{0.3\textwidth}
    \centerline{
    \includegraphics[width=2.1in,angle=0,scale=1.6]{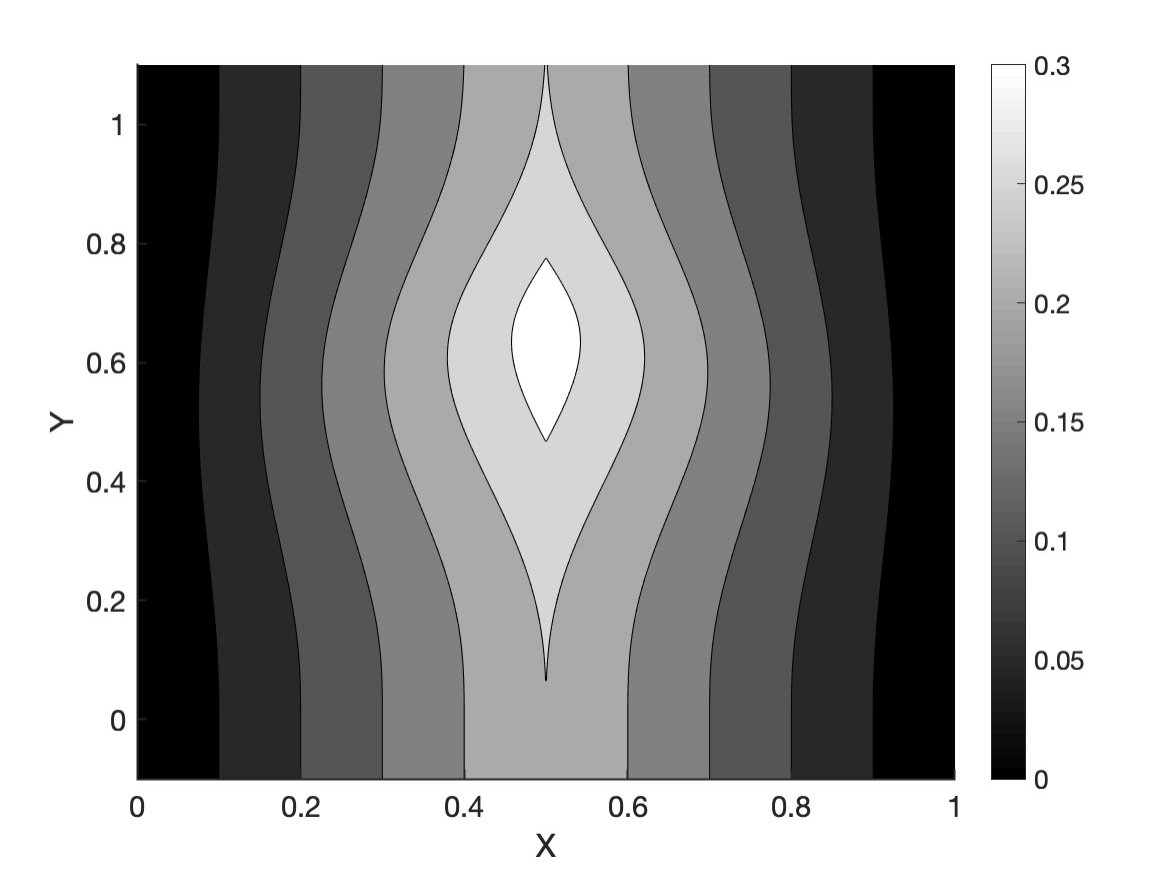}}
    %\caption{Caption 3}
   %  \vspace{-1.2in}
 \qquad   (a2) 
%    \medskip
  \end{minipage}
\begin{minipage}[b]{0.3\textwidth}
    \centerline{
    \includegraphics[width=2.1in,angle=0,scale=1.6]{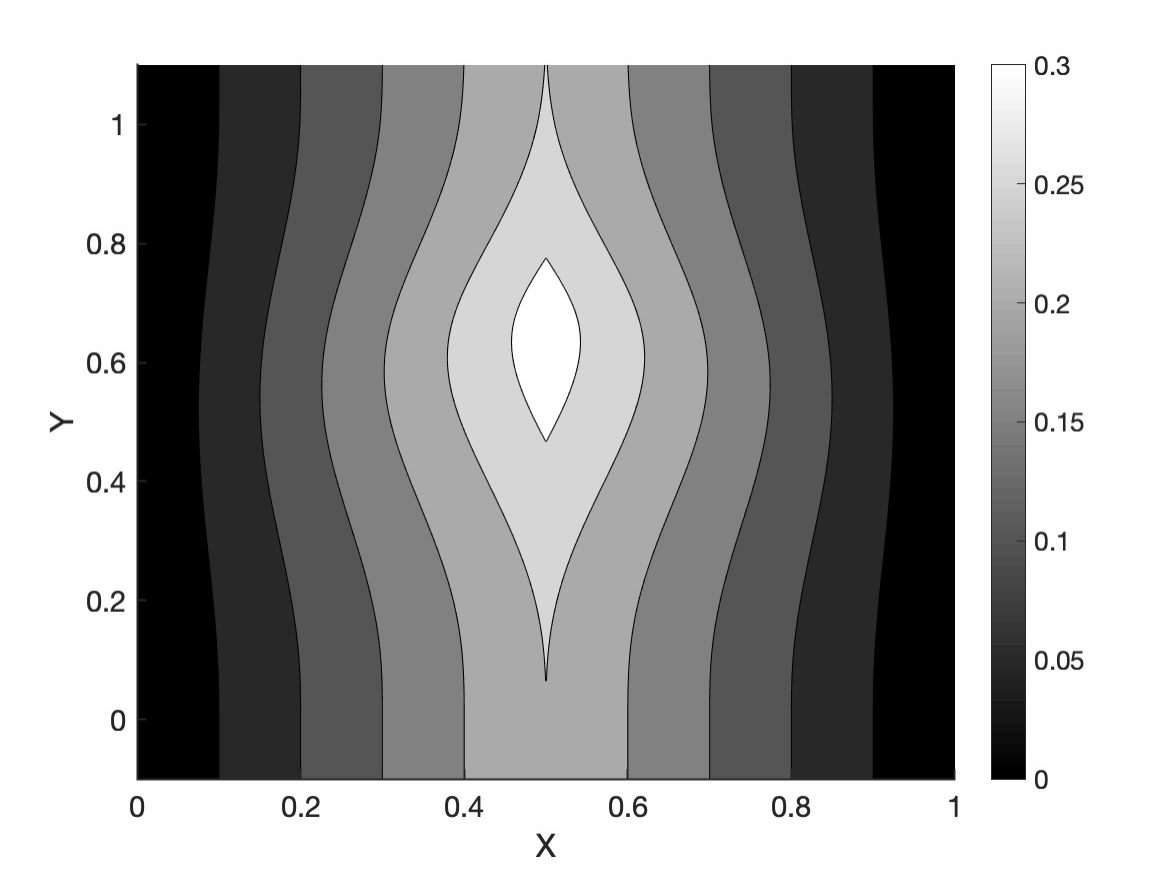}}
   % \caption{Caption 2}
 %  \vspace{-1.2in}
  \qquad      (b1) 
 %  \medskip
  \end{minipage}
  \hspace{0.2\textwidth}
 %   \vspace{-1.0in}
  \begin{minipage}[b]{0.3\textwidth}
    \centerline{
    \includegraphics[width=2.1in,angle=0,scale=1.6]{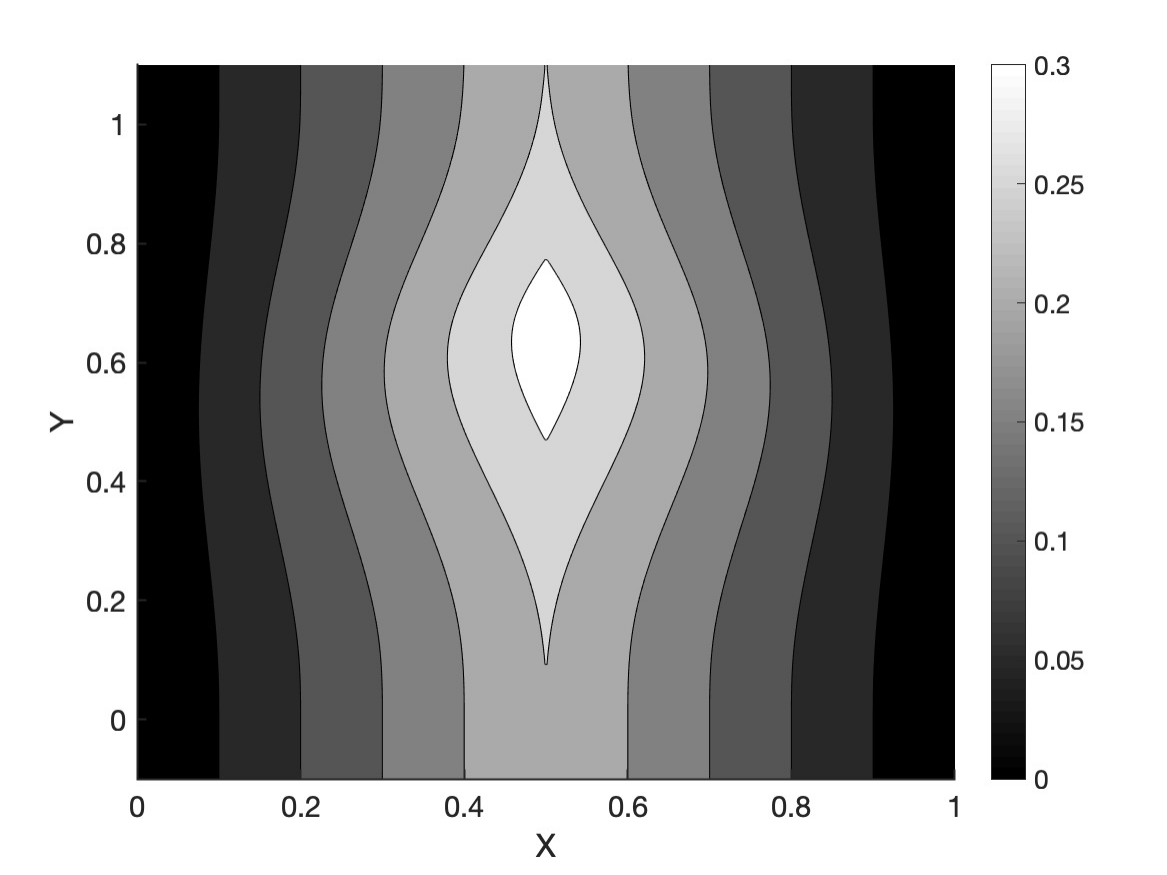}}
    %\caption{Caption 3}
   %  \vspace{-1.2in}
  \qquad      (b2)  
%    \medskip
  \end{minipage}  
    \caption{Numerical solution of Example \ref{2d_static_eikonal_6}.  (a1) First order point-wise sweeping, $\alpha=1.5$, $h=0.001$;  (a2) Newton-type sweep with WENO5, $\alpha=1.5$, $h=0.001$; (b1) Newton-type sweep with WENO7, $\alpha=1.5$, $h=0.001$; (b2) Newton-type sweep with WENO9, $\alpha=5.0$, $h=0.0025$. } \label{figure_ex6}
\end{figure}

\newpage
%++++++++++++++++++++++++++++++++++++++++++++++++++
\begin{exam}\label{2d_static_eikonal_6}
We consider an example with two incident planar wavefronts. For the general Eikonal equation (\ref{geikonal}) on the computational domain $\Omega = [0, 1] \times [-0.1, 1.1]$, the wave speed $F$ is taken as 
 
\begin{equation}\label{ex6_F}
F(x,y) = \left\{
\begin{aligned}
& 2-\frac{1}{2}\cos^2(\pi(y-\frac{1}{2})), \quad \text{if} \quad 0 \le y \le 1 \\
& 2, \qquad \qquad \qquad \qquad \qquad \text{otherwise.}
\end{aligned}
\right.
\end{equation}
The ambient velocity is chosen as $\mathbf{v(x)} = (v_1, v_2)^{\text{T}} = (0, 0.8)^\text{T}$, while the boundary condition $T(x, y)=0$ is imposed along the sides $x = 0$ and $x=1$ such that there are two planar wavefronts normal to the x-axis entering the computational domain. The numerical parameters for each method are as follows. $\epsilon = 10^{-13}$ for first-order explicit Gauss-Seidel Lax-Friedrichs sweep method and  $\epsilon = 10^{-12}$ for implicit Newton-type WENO sweeping methods. For WENO5 and WENO7, we set $\alpha = 1.5$ and $\epsilon_{\text{WENO}} = 10^{-9}$, while $\alpha = 5.0$ and $\epsilon_{\text{WENO}} = 10^{-5}$ for WENO9. The mesh size is $h=0.001$ for first order, WENO5 and WENO7 schemes, and $h=0.025$ for WENO9. From the numerical solutions in Figure \ref{figure_ex6}, we can clearly see that all the numerical methods can produce comparable results and the wavefronts are symmetric with respect to $x=0.5$. 

\end{exam}

%++++++++++++++++++++++++++++++++++++++++++++++++++
\pagebreak
\newpage
\begin{exam}\label{3d_application_example}

In this example, we consider a Gaussian wave speed, where
\begin{equation}\label{3d_application}
F(x,y,z) = 2.7 -0.5y-0.6e^{-10(x-0.3)^2-15(y-0.4)^2-20(z-0.3)^2}
\end{equation}
with a constant ambient velocity $\mathbf{v} = (0.4, 0.8, 1.2)$. The source condition $T(x,y,z) = 0$ is imposed at the point source $(x_s, y_s, z_s) = (0.5, 0.5, 0.5)$ in the computational domain $[0, 1]^3$. We simulate the problem with a uniform mesh partition $h=0.01$. We use the first-order implicit Newton-type WENO sweeping method to initialize the higher-order schemes. The iteration stopping criterions are taken as $\epsilon = 10^{-8}$ for the first-order scheme, while $\epsilon=10^{-9}$ for higher-order schemes. We choose $\alpha=1$ for all the schemes, $\epsilon_{\text{WENO}} = 10^{-5}$ for WENO5, $\epsilon_{\text{WENO}} = 10^{-1}$ for WENO7 and $\epsilon_{\text{WENO}} = 2$ for WENO9.

\end{exam}

\begin{figure}[!htb]
  %\centering
 % \vspace{-1.0in}
  \begin{minipage}[b]{0.3\textwidth}
    \centerline{
    \includegraphics[width=2.1in,angle=0,scale=1.4]{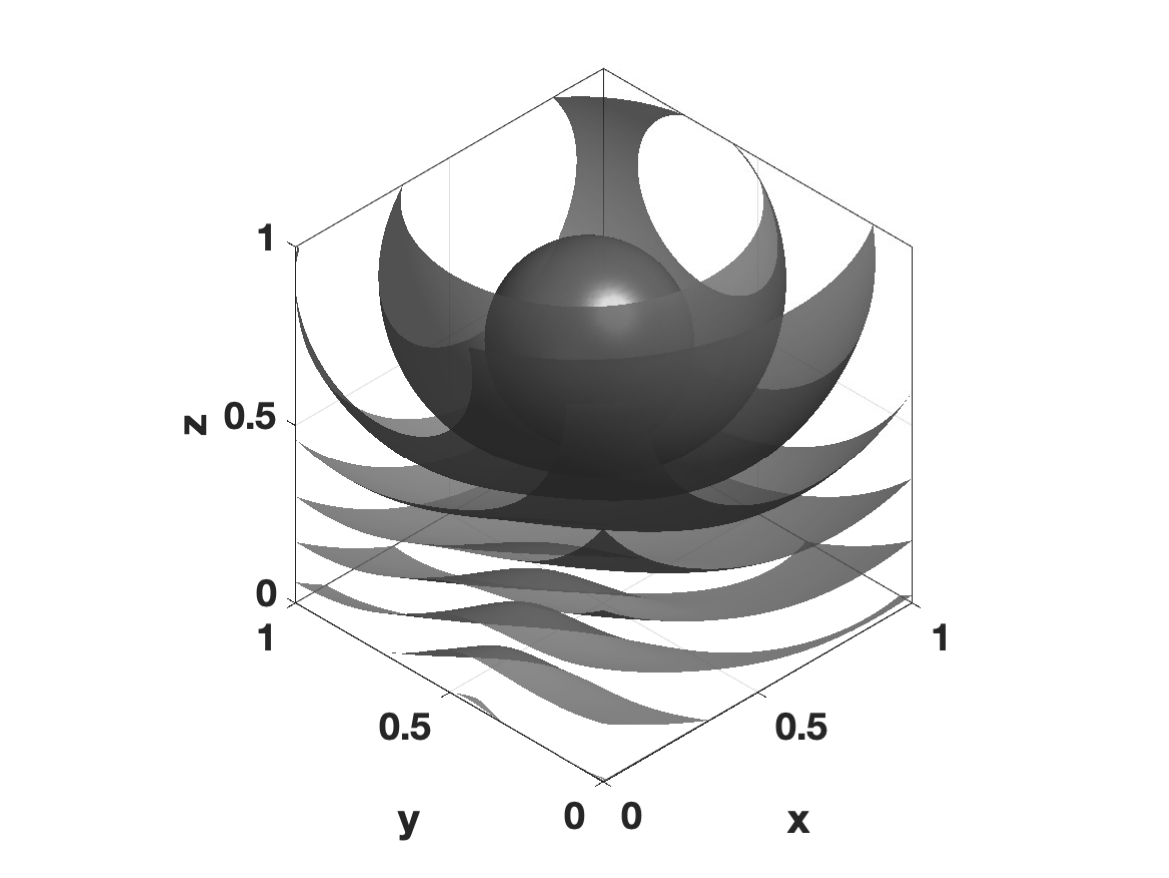}}
   % \caption{Caption 2}
 %  \vspace{-1.2in}
  \qquad     (a1) 
 %  \medskip
  \end{minipage}
  \hspace{0.2\textwidth}
 %   \vspace{-1.0in}
  \begin{minipage}[b]{0.3\textwidth}
    \centerline{
    \includegraphics[width=2.1in,angle=0,scale=1.4]{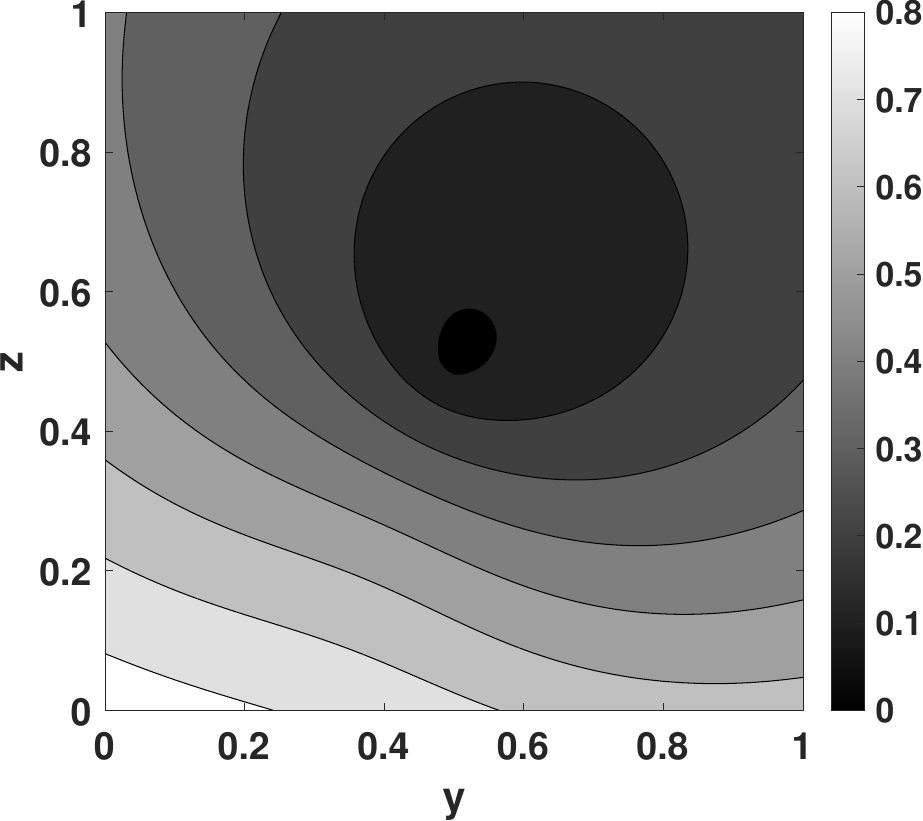}}
    %\caption{Caption 3}
   %  \vspace{-1.2in}
 \qquad   (a2) 
%    \medskip
  \end{minipage}
\begin{minipage}[b]{0.3\textwidth}
    \centerline{
    \includegraphics[width=2.1in,angle=0,scale=1.4]{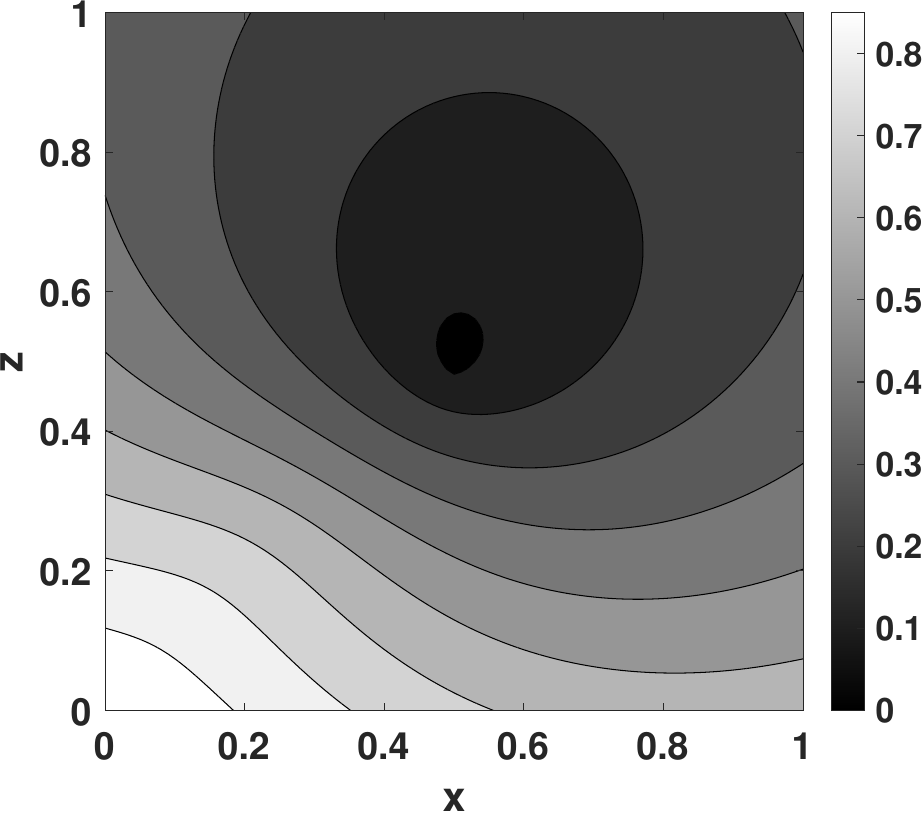}}
   % \caption{Caption 2}
 %  \vspace{-1.2in}
  \qquad      (b1) 
 %  \medskip
  \end{minipage}
  \hspace{0.2\textwidth}
 %   \vspace{-1.0in}
  \begin{minipage}[b]{0.3\textwidth}
    \centerline{
    \includegraphics[width=2.1in,angle=0,scale=1.4]{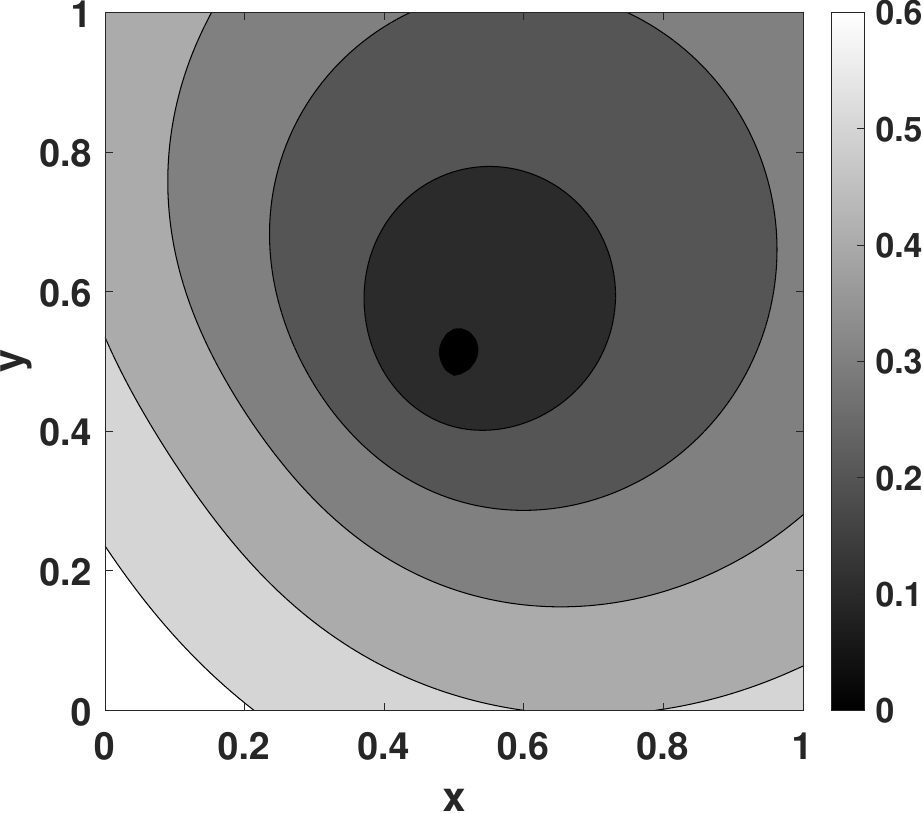}}
    %\caption{Caption 3}
   %  \vspace{-1.2in}
  \qquad      (b2)  
%    \medskip
  \end{minipage}  
    \caption{Numerical solution of Example \ref{3d_application_example} with first order implicit Newton-type sweep.  (a1) Isosurfaces of the numerical solution from $T=0.1$ to $T=1.1$ with an increment of $0.1$;  (a2) Cross-section of the numerical solution along $x=0.5$; (b1) Cross-section of the numerical solution along $y=0.5$; (b2) Cross-section of the numerical solution along $z=0.5$. } \label{figure_ex9-weno1}
\end{figure}

%weno5
\begin{figure}[!htb]
  %\centering
 % \vspace{-1.0in}
  \begin{minipage}[b]{0.3\textwidth}
    \centerline{
    \includegraphics[width=2.1in,angle=0,scale=1.4]{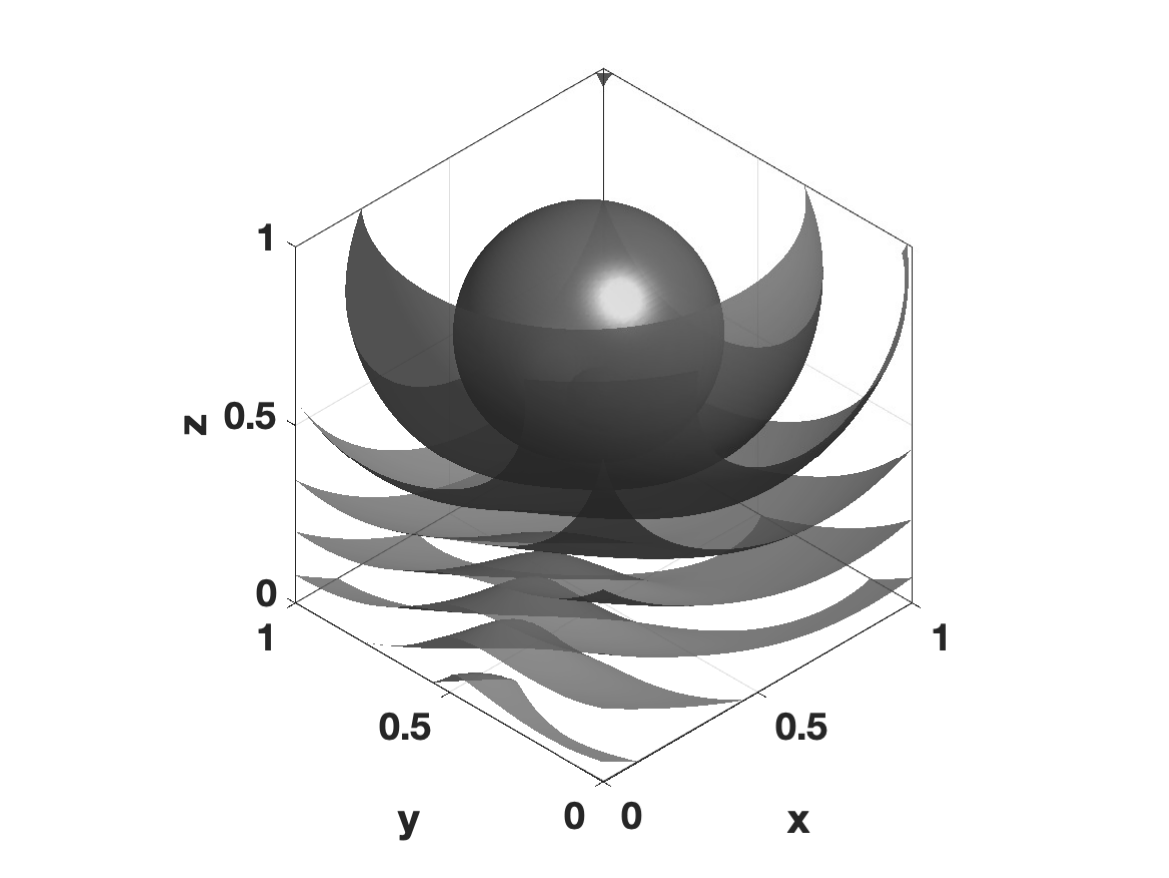}}
   % \caption{Caption 2}
 %  \vspace{-1.2in}
  \qquad     (a1) 
 %  \medskip
  \end{minipage}
  \hspace{0.2\textwidth}
 %   \vspace{-1.0in}
  \begin{minipage}[b]{0.3\textwidth}
    \centerline{
    \includegraphics[width=2.1in,angle=0,scale=1.4]{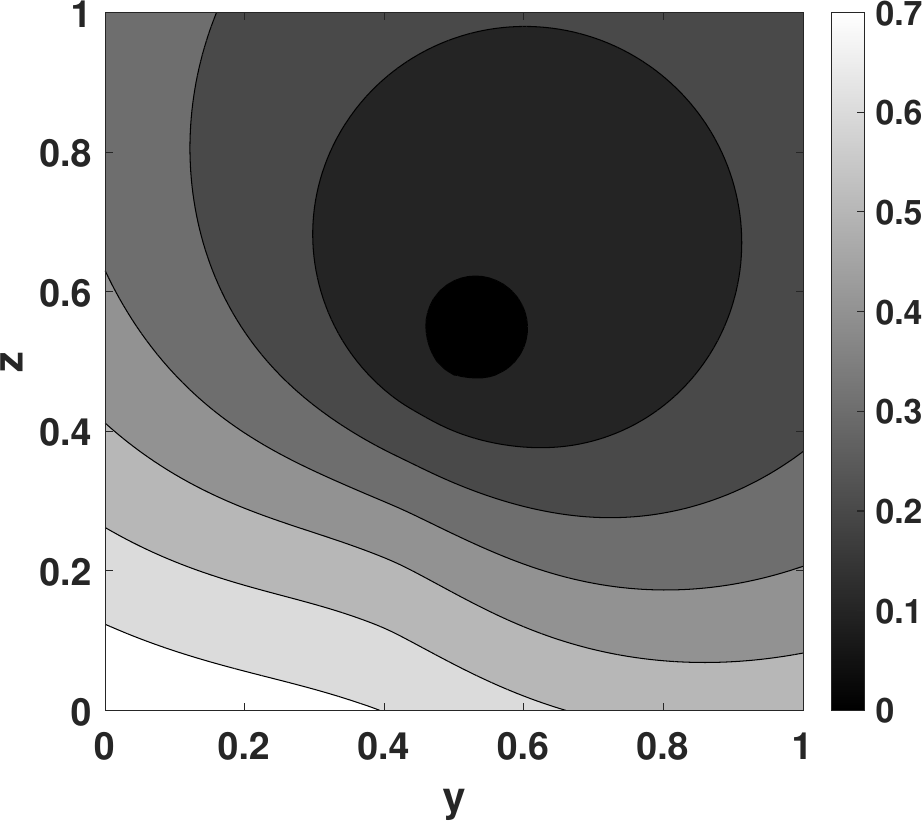}}
    %\caption{Caption 3}
   %  \vspace{-1.2in}
 \qquad   (a2) 
%    \medskip
  \end{minipage}
\begin{minipage}[b]{0.3\textwidth}
    \centerline{
    \includegraphics[width=2.1in,angle=0,scale=1.4]{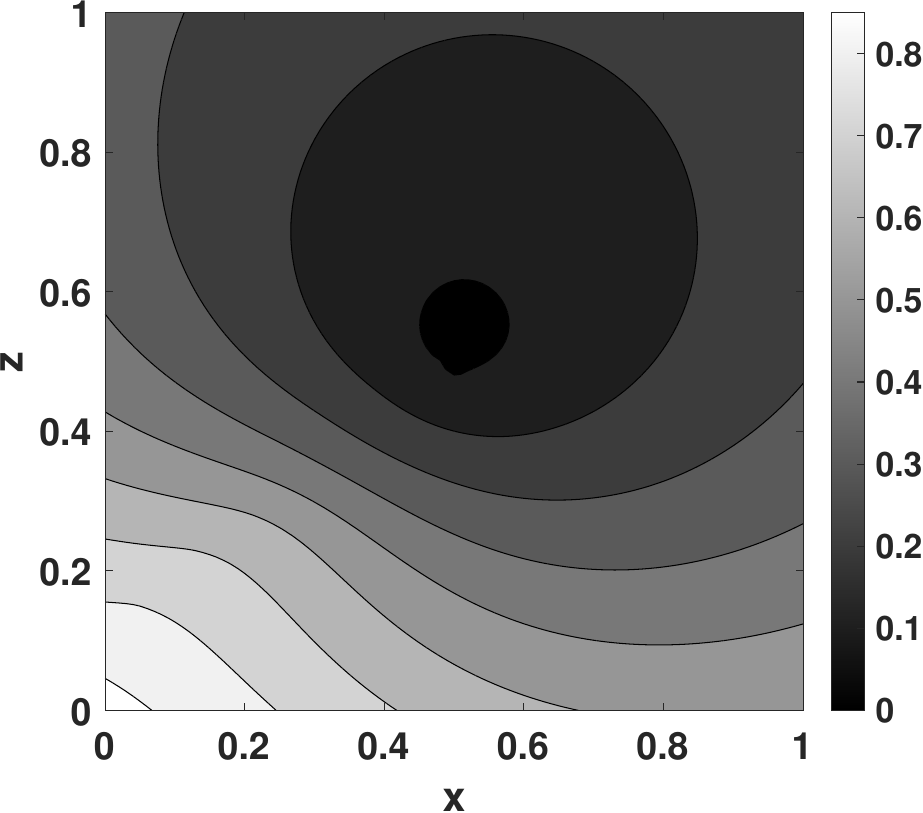}}
   % \caption{Caption 2}
 %  \vspace{-1.2in}
  \qquad      (b1) 
 %  \medskip
  \end{minipage}
  \hspace{0.2\textwidth}
 %   \vspace{-1.0in}
  \begin{minipage}[b]{0.3\textwidth}
    \centerline{
    \includegraphics[width=2.1in,angle=0,scale=1.4]{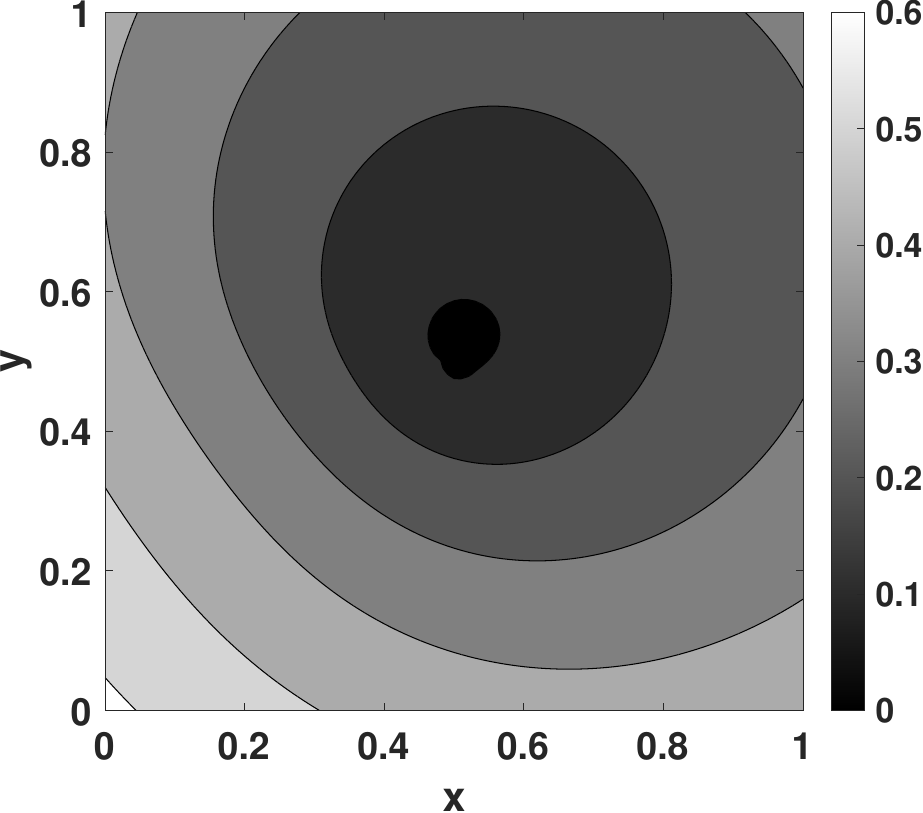}}
    %\caption{Caption 3}
   %  \vspace{-1.2in}
  \qquad      (b2)  
%    \medskip
  \end{minipage}  
    \caption{Numerical solution of Example \ref{3d_application_example} with implicit WENO5 Newton-type sweep.  (a1) Isosurfaces of the numerical solution from $T=0.1$ to $T=1.1$ with an increment of $0.1$;  (a2) Cross-section of the numerical solution along $x=0.5$; (b1) Cross-section of the numerical solution along $y=0.5$; (b2) Cross-section of the numerical solution along $z=0.5$. } \label{figure_ex9-weno5}
\end{figure}

%weno5
\begin{figure}[!htb]
  %\centering
 % \vspace{-1.0in}
  \begin{minipage}[b]{0.3\textwidth}
    \centerline{
    \includegraphics[width=2.1in,angle=0,scale=1.4]{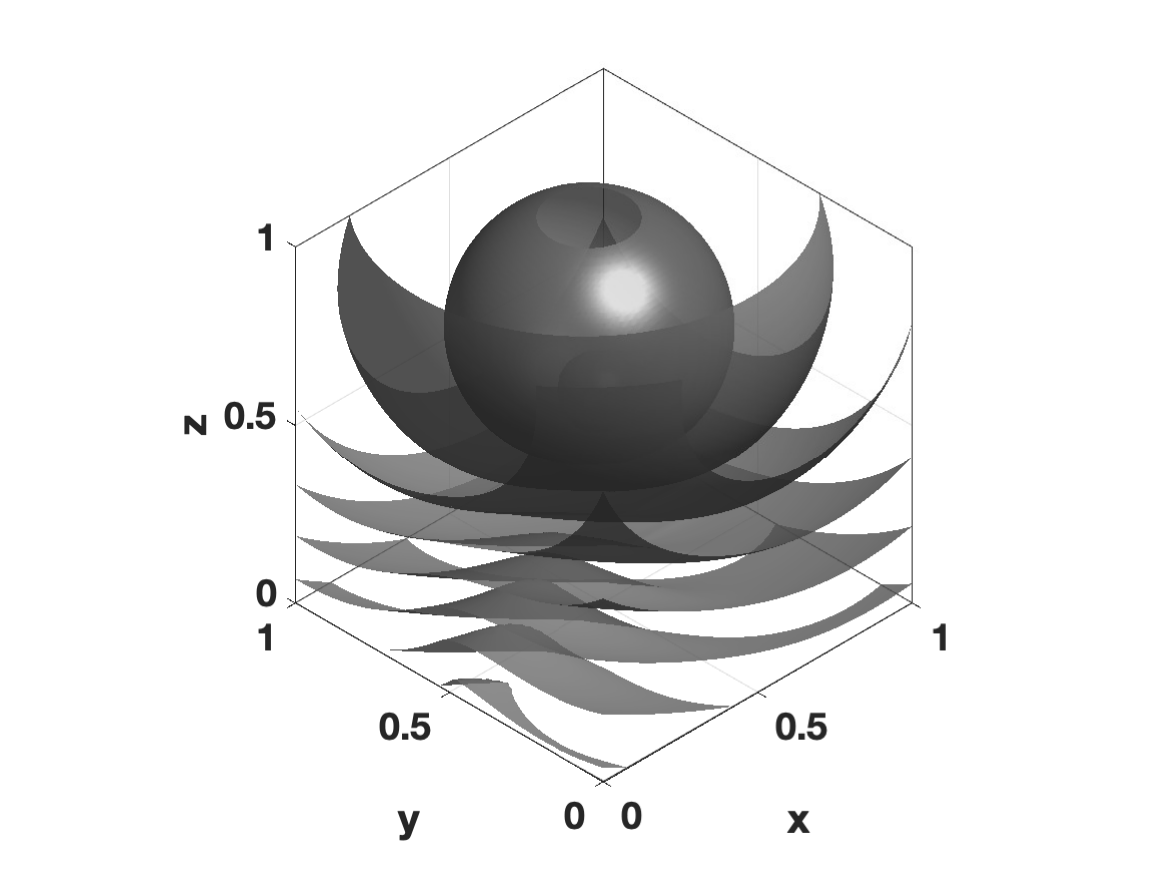}}
   % \caption{Caption 2}
 %  \vspace{-1.2in}
  \qquad     (a1) 
 %  \medskip
  \end{minipage}
  \hspace{0.2\textwidth}
 %   \vspace{-1.0in}
  \begin{minipage}[b]{0.3\textwidth}
    \centerline{
    \includegraphics[width=2.1in,angle=0,scale=1.4]{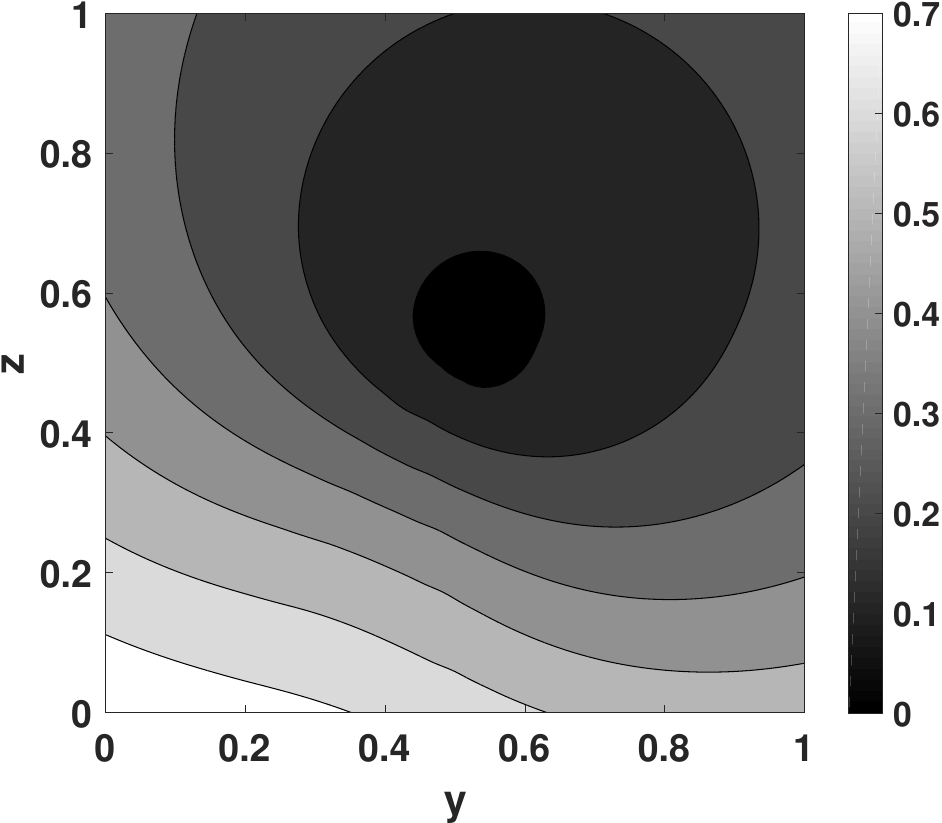}}
    %\caption{Caption 3}
   %  \vspace{-1.2in}
 \qquad   (a2) 
%    \medskip
  \end{minipage}
\begin{minipage}[b]{0.3\textwidth}
    \centerline{
    \includegraphics[width=2.1in,angle=0,scale=1.4]{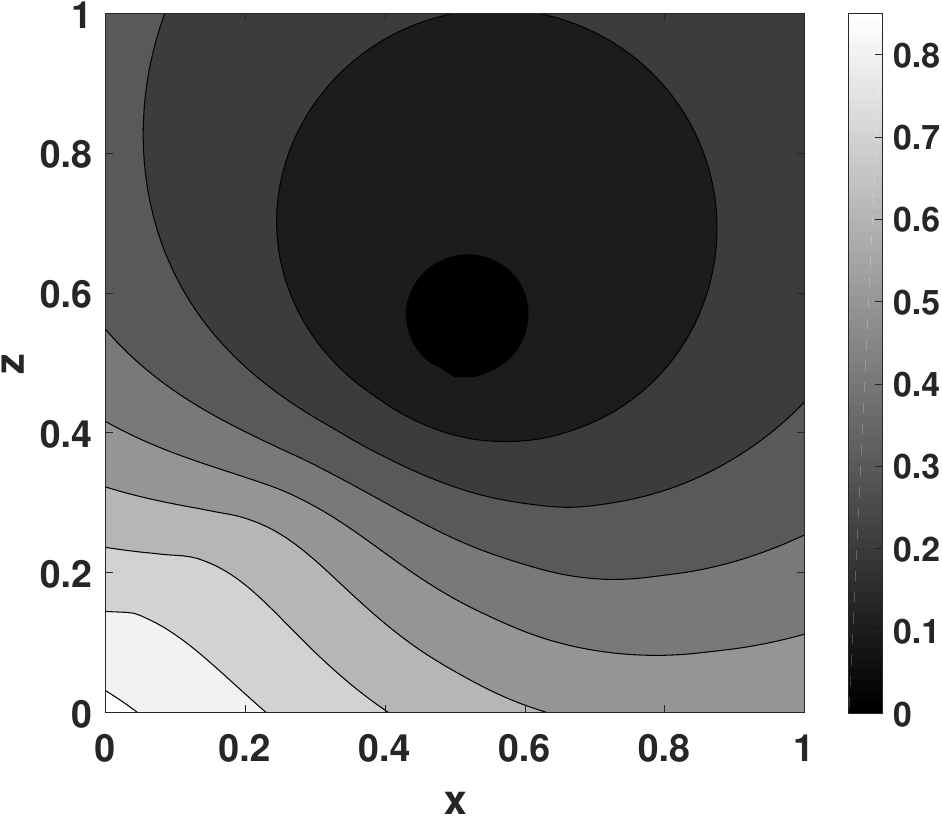}}
   % \caption{Caption 2}
 %  \vspace{-1.2in}
  \qquad      (b1) 
 %  \medskip
  \end{minipage}
  \hspace{0.2\textwidth}
 %   \vspace{-1.0in}
  \begin{minipage}[b]{0.3\textwidth}
    \centerline{
    \includegraphics[width=2.1in,angle=0,scale=1.4]{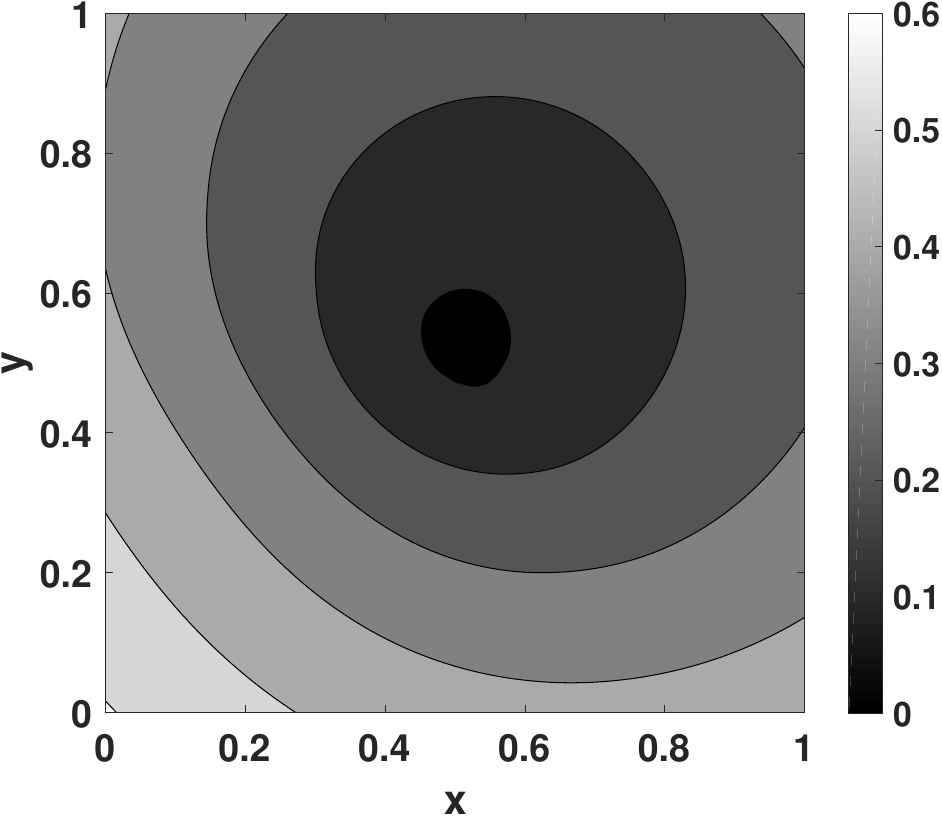}}
    %\caption{Caption 3}
   %  \vspace{-1.2in}
  \qquad      (b2)  
%    \medskip
  \end{minipage}  
    \caption{Numerical solution of Example \ref{3d_application_example} with implicit WENO7 Newton-type sweep.  (a1) Isosurfaces of the numerical solution from $T=0.1$ to $T=1.1$ with an increment of $0.1$;  (a2) Cross-section of the numerical solution along $x=0.5$; (b1) Cross-section of the numerical solution along $y=0.5$; (b2) Cross-section of the numerical solution along $z=0.5$. } \label{figure_ex9-weno7}
\end{figure}

%weno9
\begin{figure}[!htb]
  %\centering
 % \vspace{-1.0in}
  \begin{minipage}[b]{0.3\textwidth}
    \centerline{
    \includegraphics[width=2.1in,angle=0,scale=1.4]{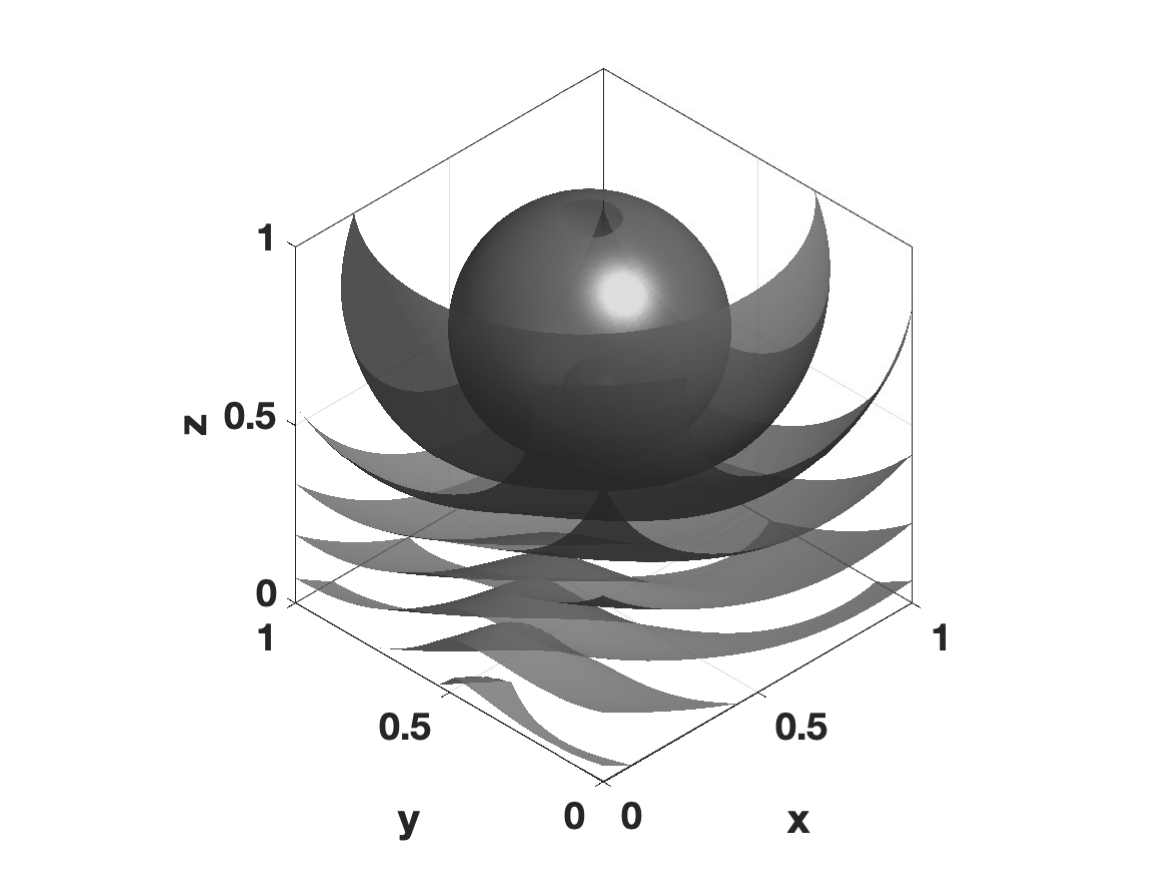}}
   % \caption{Caption 2}
 %  \vspace{-1.2in}
  \qquad     (a1) 
 %  \medskip
  \end{minipage}
  \hspace{0.2\textwidth}
 %   \vspace{-1.0in}
  \begin{minipage}[b]{0.3\textwidth}
    \centerline{
    \includegraphics[width=2.1in,angle=0,scale=1.4]{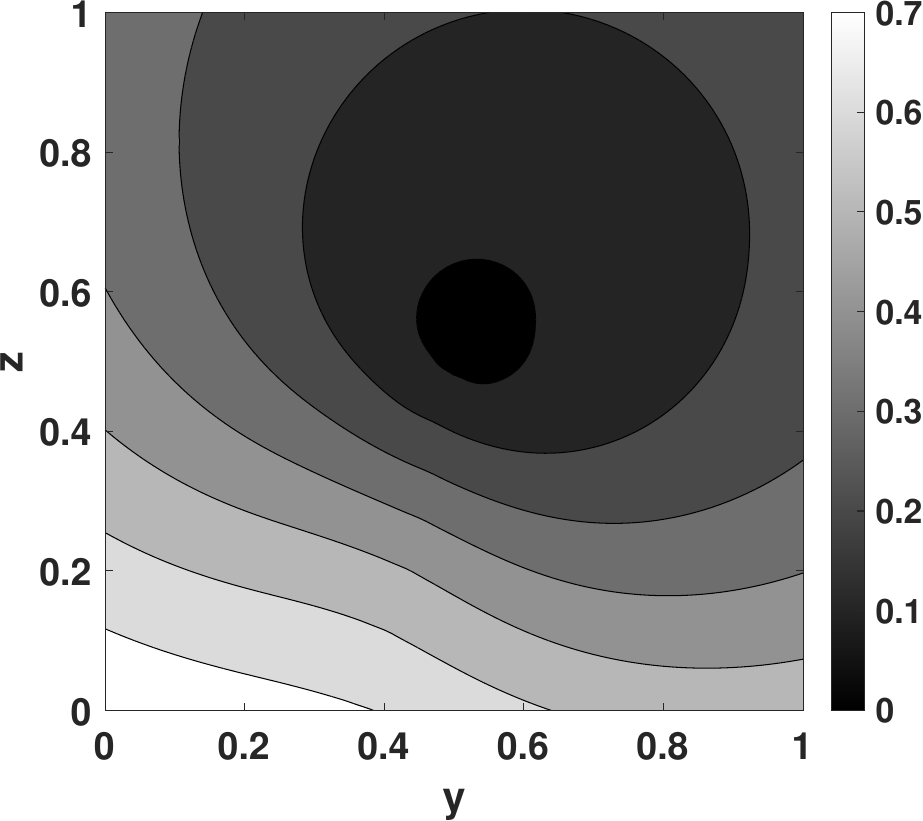}}
    %\caption{Caption 3}
   %  \vspace{-1.2in}
 \qquad   (a2) 
%    \medskip
  \end{minipage}
\begin{minipage}[b]{0.3\textwidth}
    \centerline{
    \includegraphics[width=2.1in,angle=0,scale=1.4]{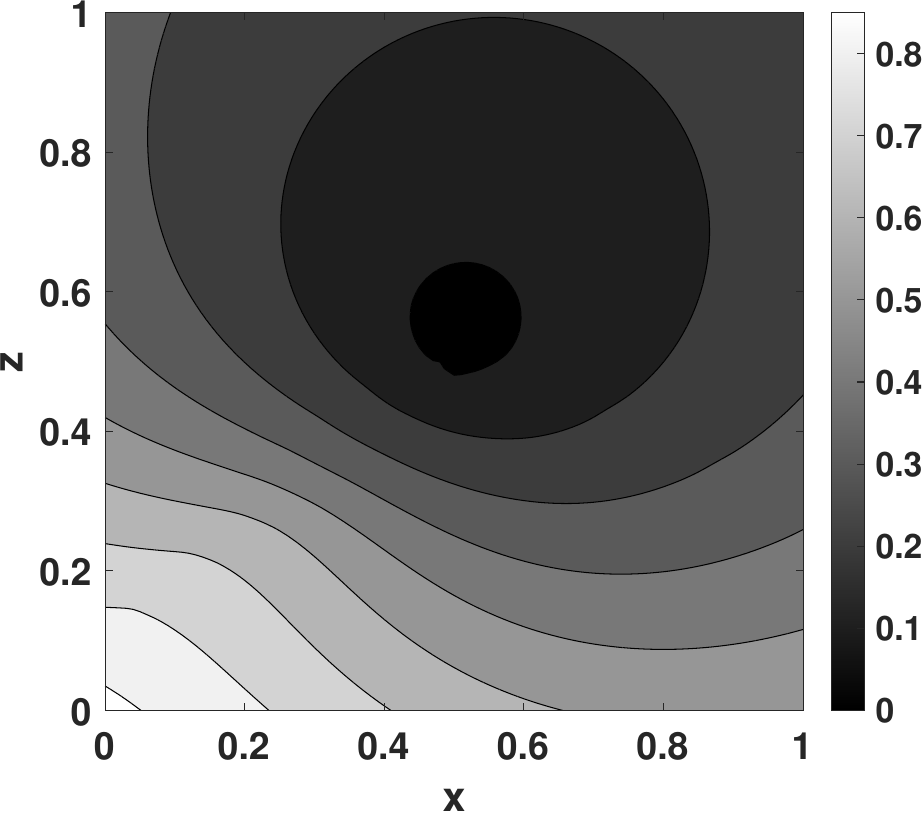}}
   % \caption{Caption 2}
 %  \vspace{-1.2in}
  \qquad      (b1) 
 %  \medskip
  \end{minipage}
  \hspace{0.2\textwidth}
 %   \vspace{-1.0in}
  \begin{minipage}[b]{0.3\textwidth}
    \centerline{
    \includegraphics[width=2.1in,angle=0,scale=1.4]{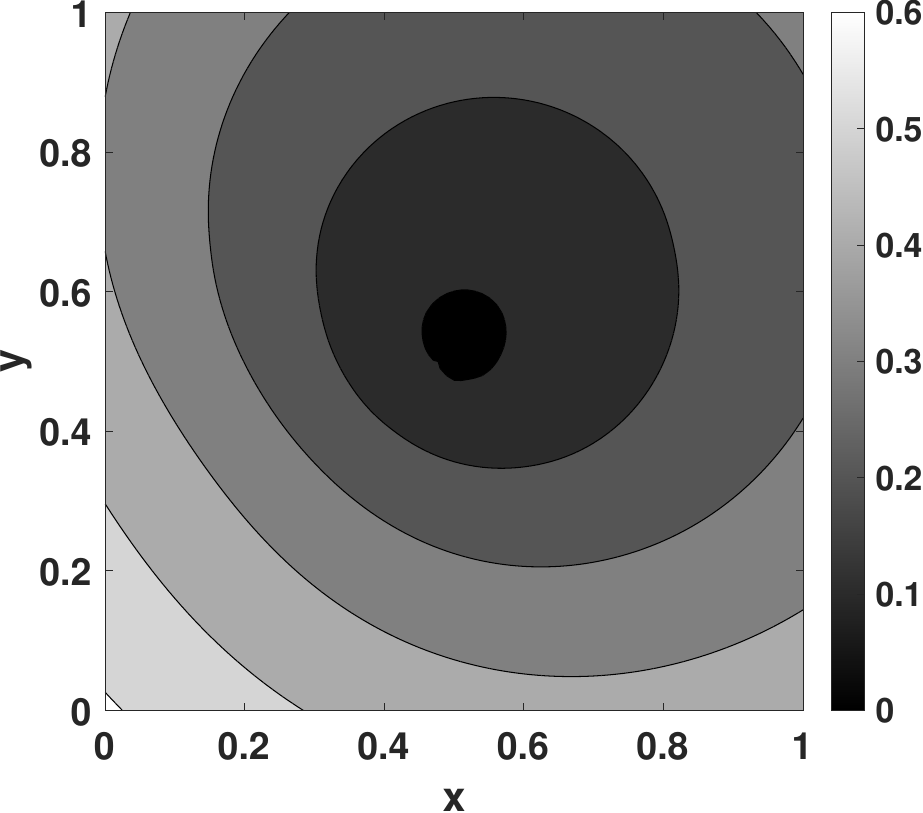}}
    %\caption{Caption 3}
   %  \vspace{-1.2in}
  \qquad      (b2)  
%    \medskip
  \end{minipage}  
    \caption{Numerical solution of Example \ref{3d_application_example} with implicit WENO9 Newton-type sweep.  (a1) Isosurfaces of the numerical solution from $T=0.1$ to $T=1.1$ with an increment of $0.1$;  (a2) Cross-section of the numerical solution along $x=0.5$; (b1) Cross-section of the numerical solution along $y=0.5$; (b2) Cross-section of the numerical solution along $z=0.5$. } \label{figure_ex9-weno9}
\end{figure}

% ===== END numerical_staticHJ.tex =====

\clearpage
% section 4
% ===== BEGIN conclusion.tex =====

\pagebreak
\section{Conclusion}
\label{numer}

In this paper, we have developed a class of high-order Newton type Gauss-Seidel Lax-Friedrichs fast sweeping methods for solving the generalized Eikonal equations arising from wave propagation in a moving fluid. By incorporating fifth-, seventh- and ninth-order WENO reconstruction into the Newton-type line-wise solver of Li and Qian \cite{li2020newton}, we achieve high-order accuracy while preserving the key structure of the Jacobian. By using the first order Lax-Friedrichs discretization, the Jacobian $\nabla \mathbf{F}_j$ or $\nabla \mathbf{F}_i$ remain tridiagonal and strictly diagonally dominant for all orders considered and ensure the Thomas algorithm's convergence with $O(N)$ complexity. Numerical examples in both two and three spatial dimensions confirm that the proposed schemes achieve their designed orders of accuracy for smooth solutions. We have also examined three alternating sweeping strategies, that include column-wise, row-wise, and combined column-row-wise, within the Gauss–Seidel iteration framework.   

Our numerical experiments also suggest that the choice of sweeping strategy may affect both the convergence behavior and the number of iterations required, and this effect appears to be related to the structure of the background velocity field $\mathbf{v}(\mathbf{x})$. The convergence of the higher-order schemes is sensitive to the choices of numerical parameters which may be attributed to the oscillatory nature of high-degree interpolating polynomials. In the future, we plan to investigate a more robust WENO reconstruction that reduces this sensitivity. Meanwhile, a systematic understanding of how the sweeping strategy should be selected based on the structure of the underlying problem and the background velocity field is an important question we plan to explore.
% ===== END conclusion.tex =====

% section 5
%\input{acknowledgment.tex}

\clearpage
%\bibliographystyle{plain}
% ===== BEGIN embedded bibliography =====

% ===== END embedded bibliography =====

\end{document}